\documentclass[11pt,oneside,a4paper]{article}
\usepackage[top=1 in, bottom=1 in, left=0.85 in, right=0.85 in]{geometry}
\usepackage{hyperref}
\usepackage{float}
\usepackage{amsmath,amssymb,graphicx,url}
\usepackage{amsthm}
\usepackage{thmtools,thm-restate,wrapfig,enumitem,mathabx}
\usepackage{xcolor}
\usepackage[noend,linesnumbered,ruled,vlined]{algorithm2e}
\usepackage{bbold}
\usepackage[utf8]{inputenc}
\usepackage[T1]{fontenc}
\usepackage{booktabs}
\usepackage{amsfonts}
\usepackage{nicefrac}
\usepackage{microtype}
\usepackage{xcolor}
\usepackage{subcaption}
\usepackage{multirow}
\usepackage{bbm}
\usepackage{mathtools}
\usepackage[compress]{cite}
\usepackage[thinc]{esdiff}

\newtheorem{proposition}{Proposition}
\newtheorem{assumption}{Assumption}
\newtheorem{theorem}{Theorem}

\newtheorem{lemma}{Lemma}
\newtheorem{corollary}{Corollary}

\definecolor{darkgreen}{RGB}{25,107,36}
\definecolor{darkteal}{RGB}{21,96,130}

\newcommand{\blue}[1]{\textcolor{blue}{#1}}

\newcommand{\orange}[1]{\textcolor{orange}{#1}}

\newcommand{\violet}[1]{\textcolor{violet}{#1}}

\newcommand{\darkgreen}[1]{\textcolor{darkgreen}{#1}}

\DeclareMathOperator*{\argmax}{arg\,max}

\newcommand{\expt}{\mathbb{E}}
\newcommand{\prob}{\mathbb{P}}
\let\oldnl\nl
\newcommand{\nonl}{\renewcommand{\nl}{\let\nl\oldnl}}

\allowdisplaybreaks

\title{Learning-Based Pricing and Matching for Two-Sided Queues}

\author{%
  Zixian Yang\\
  University of Michigan\\
  Ann Arbor, MI 48109, USA\\
  \texttt{zixian@umich.edu}\\
  \and
  Lei Ying\\
  University of Michigan\\
  Ann Arbor, MI 48109, USA\\
  \texttt{leiying@umich.edu}
}
\date{}

\begin{document}

\maketitle

\begin{abstract}
We consider a dynamic system with multiple types of customers and servers. Each type of waiting customer or server joins a separate queue, forming a bipartite graph with customer-side queues and server-side queues. The platform can match the servers and customers if their types are compatible. The matched pairs then leave the system. The platform will charge a customer a price according to their type when they arrive and will pay a server a price according to their type. The arrival rate of each queue is determined by the price according to some unknown demand or supply functions. Our goal is to design pricing and matching algorithms to maximize the profit of the platform with unknown demand and supply functions, while keeping queue lengths of both customers and servers below a predetermined threshold. This system can be used to model two-sided markets such as ride-sharing markets with passengers and drivers. The difficulties of the problem include simultaneous learning and decision making, and the tradeoff between maximizing profit and minimizing queue length. We use a longest-queue-first matching algorithm and propose a learning-based pricing algorithm, which combines gradient-free stochastic projected gradient ascent with bisection search. 
We prove that our proposed algorithm yields a sublinear regret $\tilde{O}(T^{5/6})$ and anytime queue-length bound $\tilde{O}(T^{1/6})$, where $T$ is the time horizon. We further establish a tradeoff between the regret bound and the queue-length bound: $\tilde{O}(T^{1-\gamma})$  versus $\tilde{O}(T^{\gamma})$ for $\gamma \in (0, 1/6].$
\end{abstract}

\section{Introduction}

We study pricing and matching in two-sided queueing systems with multiple types of customers and servers. Take ride sharing (Uber, Lyft, etc.) as an example. A ride order (request) can be viewed as a customer and an available driver can be viewed as a server. A ride order can have multiple types. For example, we can view different number of passengers in the order as different customer types. We can also view different pick-up locations of the order as different customer types. A driver can also have multiple types. For example, we can view different vehicle capacities as different server types. We can also view different locations of the drivers as different server types. Each type of waiting ride order or driver joins a separate queue, forming a bipartite graph with customer-side (passenger-side) queues and server-side (driver-side) queues. We consider a discrete-time system. At each time slot, the ride-sharing platform will determine a price for each customer or server type. If a passenger accepts the price, they will be charged and join the queue. Similarly, if a driver accepts the price, they will be paid and join the queue. Therefore, the arrival rate of each type of ride order (or driver) is governed by the corresponding price according to some unknown demand (or supply) function.
After that, the platform will match the drivers and ride orders in the queues if their types are compatible. For example, if the customer type is the number of passengers in the ride order and the server type is the capacity of the vehicle, then their types are compatible if the capacity of the vehicle is no less than the number of passengers in the ride order. The matched pairs then leave the system. The same process repeats in the next time slot. The profit of the platform is equal to the income from ride orders minus the cost paid to drivers. Assuming the demand and supply functions are unknown, the goal is to design pricing and matching algorithms to maximize the total profit of the platform over a finite time horizon while keeping the queue lengths of both customers and servers below a predetermined threshold.

Similar pricing and matching problems in two-sided queues have been studied in the literature~\cite{caldentey2009fcfs,adan2012exact,hu2022dynamic,aveklouris2023fluid,aveklouris2024matching,nguyen2018queueing,gurvich2015dynamic,varma2020near,chen2020pricing,vaze2022non,chen2023online,varma2023dynamic}. The works in \cite{caldentey2009fcfs,adan2012exact,hu2022dynamic,aveklouris2023fluid,aveklouris2024matching} study matching problems in two-sided queueing systems with multiple types of demand and multiple types of supply.
The work in \cite{nguyen2018queueing} studies a two-sided queueing system with on-demand servers with the goal of designing adaptive server invitation scheme to minimize waiting times.
The work in \cite{gurvich2015dynamic} studies a matching problem for multisided queues. However, these works \cite{caldentey2009fcfs}-\cite{gurvich2015dynamic} do not consider pricing. The works in \cite{varma2020near,chen2020pricing} study pricing and matching problems in two-sided queues but they consider scenarios with strategic demand side and/or supply side, which is different from our setting, where customers and servers are not strategic and their arrival rates are based on demand and supply functions. The work in \cite{vaze2022non} considers a two-sided queueing system where servers arrive at a fixed rate and the arrival rate of customers is controlled by the price. They consider one-sided pricing and single-type customers and servers, while we consider two-sided pricing and matching and multiple customer and server types.
The work in \cite{chen2023online} studies learning and pricing in a single-sided queue with unknown demand curve and service distribution.
The most relevant work in the literature is \cite{varma2023dynamic}, whose model is similar to ours. The differences between \cite{varma2023dynamic} and ours include:
\begin{itemize}[leftmargin=*]
    \item The key difference is that they assume that the platform knows the demand and supply functions while in our setting the demand and supply functions are unknown.
    \item The objective in \cite{varma2023dynamic} is the asymptotic long run expected profit of the platform subtracting a penalty which is linear in the queue length, while our objective is to maximize the total profit of the platform over a finite time horizon while keeping the queue lengths below a predetermined threshold during the entire time horizon. 
    \item In \cite{varma2023dynamic}, they consider a large-scale regime in which they scale the arrival rates and study the profit loss in steady state. We focus on finite-time analysis without scaling the arrival rates. 
\end{itemize}

The key challenge of our setting is that the demand and supply functions are unknown.  Learning the demand and supply functions explicitly will not be adaptive if the demand or supply functions change over time. Therefore, we propose to learn and make pricing and matching decisions simultaneously. One idea is to formulate the problem into a Markov decision process (MDP) and then use reinforcement learning (RL) to solve the MDP approximately. However, the state space of this MDP increases exponentially with the number of queues (i.e., customer and server types), so it is difficult to solve due to the curse of dimensionality. Instead, we propose to use a longest-queue-first matching algorithm and a learning-based pricing algorithm, which combines gradient-free (zero-order) stochastic projected gradient ascent~\cite{agarwal2010optimal} with bisection search. We also borrow the idea of fluid pricing policy from \cite{varma2023dynamic} to control the queue length 
by rejecting arrivals according to a predetermined queue-length threshold.
In this paper, we prove that our proposed algorithm yields a sublinear regret $\tilde{O}(T^{5/6})$ when being compared with a fluid-based baseline, while ensuring an anytime queue-length bound of $\tilde{O}(T^{2/3})$ and an average queue-length bound of $\tilde{O}(T^{1/2})$, where $T$ is the time horizon.
We also establish tradeoffs between the regret bound and the anytime queue-length bound, as well as between regret bound and the average queue-length bound.
Furthermore, assuming the availability of a set of balanced arrival rates and associated price intervals, we prove a regret $\tilde{O}(T^{5/6})$ and an anytime queue-length bound of $\tilde{O}(T^{1/6})$, and establish an improved tradeoff between the regret bound and the queue-length bound: $\tilde{O}(T^{1-\gamma})$ versus $\tilde{O}(T^{\gamma})$ for $\gamma\in(0, 1/6]$.

\section{Model}

\begin{figure}[tb]
    \centering
    \includegraphics[width=0.7\linewidth]{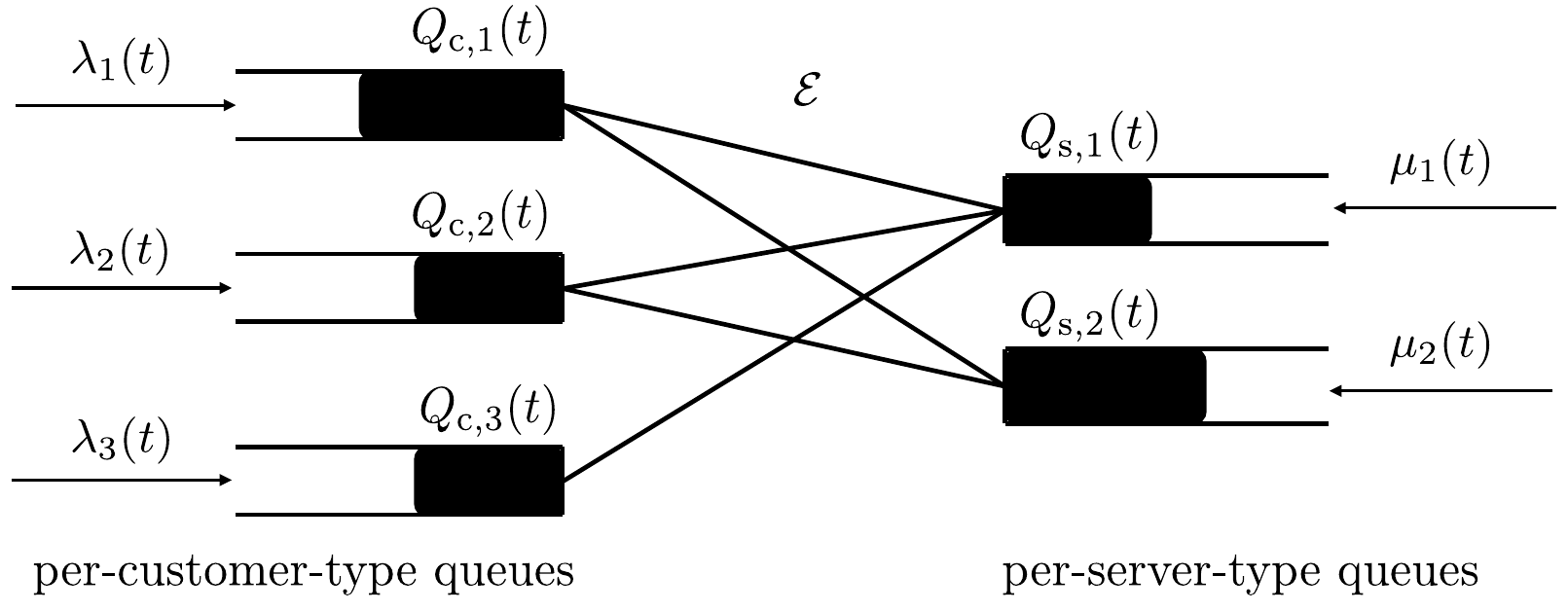}
    \caption{\small The model, an example with 3 types of customers and 2 types of servers.}
    \label{fig:model}
\end{figure}

We consider a discrete-time system with two-sided queues, a customer side and a server side. There are multiple queues on each side, representing different types of customers and different types of servers. We can also view customers as demand and servers as supply. We model the system by a bipartite graph $G({\cal I} \cup {\cal J}, {\cal E})$, where ${\cal I}=\{1,2,\ldots,I\}$ is the set of customer types and ${\cal J}=\{1,2,\ldots,J\}$ is the set of server types with $|{\cal I}|=I$ and $|{\cal J}|=J$. ${\cal E}$ is the set of all compatible links, which means that a type $i$ customer can be served by a type $j$ server if and only if $(i,j)\in {\cal E}$.
Figure~\ref{fig:model} is an example with three types of customers and two types of servers ($I=3, J=2$). 
At each time slot, there are arrivals of customers and servers.
At time slot $t$, let $A_{\mathrm{c},i}(t)$ and $A_{\mathrm{s},j}(t)$ denote the number of arrivals of type $i$ customers and type $j$ servers, respectively.
Assume that $A_{\mathrm{c},i}(t)$ and $A_{\mathrm{s},j}(t)$ follow independent Bernoulli distributions. Let $\lambda_i(t)\coloneqq \expt[A_{\mathrm{c},i}(t)]$ and $\mu_j(t)\coloneqq \expt[A_{\mathrm{s},j}(t)]$. At each time slot, after the arrival of customers and servers, the platform will match the customers and servers in the queues through compatible links. Once a customer is matched with a server, i.e., the customer is served by the server, they will leave the system immediately.
Let $Q_{\mathrm{c}, i}(t)$ and $Q_{\mathrm{s}, j}(t)$ denote the queue length of type $i$ customers and that of type $j$ servers, respectively, before arrivals and matching at time slot $t$. Assume the queues are all empty at $t=1$.

There is also a price attached to each type of customer or server. Let $p_{\mathrm{c}, i}(t)$ denote the price of the type $i$ customer and $p_{\mathrm{s}, j}(t)$ denote the price of the type $j$ server at time slot $t$. At time slot $t$, when a customer of type $i$ arrives, they will be charged $p_{\mathrm{c}, i}(t)$ units of money by the platform. When a server of type $j$ arrives, they will be paid $p_{\mathrm{s}, j}(t)$ units of money by the platform. Note that the prices determine the arrival rates. Increasing the price of type $i$ customers will reduce arrival rate $\lambda_i(t)$; and increasing the price of type $j$ servers will increase arrival rate $\mu_j(t)$. For a type $i$ customer, the demand function is denoted by $F_i:[0,1] \rightarrow [p_{\mathrm c,i,\min}, p_{\mathrm c,i,\max}]$, where the input is an arrival rate $\lambda_i\in [0,1]$, the output is the price $p_{\mathrm{c}, i}\in [p_{\mathrm c,i,\min}, p_{\mathrm c, i,\max}]$ corresponding to this arrival rate, and $p_{\mathrm c,i,\min}$ and $p_{\mathrm c, i, \max}$ are the minimum and maximum prices for type $i$ customer, respectively.
Similarly, for type $j$ server, the supply function is denoted by $G_j: [0,1] \rightarrow [p_{\mathrm s,j,\min}, p_{\mathrm s,j,\max}]$, where the input is an arrival rate $\mu_j\in [0,1]$, the output is the price $p_{\mathrm{s}, j}\in [p_{\mathrm s,j,\min}, p_{\mathrm s,j,\max}]$ corresponding to this arrival rate, and $p_{\mathrm s,j,\min}$ and $p_{\mathrm s,j,\max}$ are the minimum and maximum prices for a type $j$ server, respectively. Then we have $p_{\mathrm{c}, i}(t)=F_i(\lambda_i(t))$ and $p_{\mathrm{s}, j}(t)=G_j(\mu_j(t))$. We make the following assumption on the demand and supply functions $F_i$ and $G_j$:
\begin{assumption}\label{assum:1}
    For any customer type $i$, $F_i$ is strictly decreasing, bijective, and $L_{F_i}$-Lipschitz. Let $F_i^{-1}$ denote the inverse function of $F_i$. Assume that $F_i^{-1}$ is $L_{F_i^{-1}}$-Lipschitz.
    For any server type $j$, $G_j$ is strictly increasing, bijective, and $L_{G_j}$-Lipschitz. Let $G_j^{-1}$ denote the inverse function of $G_j$. Assume that $G_j^{-1}$ is $L_{G_j^{-1}}$-Lipschitz.
\end{assumption}
\noindent From Assumption~\ref{assum:1}, we know $F_i(1)=p_{\mathrm c,i,\min}, F_i(0) = p_{\mathrm c,i,\max}$ for any customer type $i$ and $G_j(0)=p_{\mathrm s,j,\min}, G_j(1) = p_{\mathrm s,j,\max}$ for any server type $j$.

\begin{figure}[tb]
    \centering
    \includegraphics[width=0.7\linewidth]{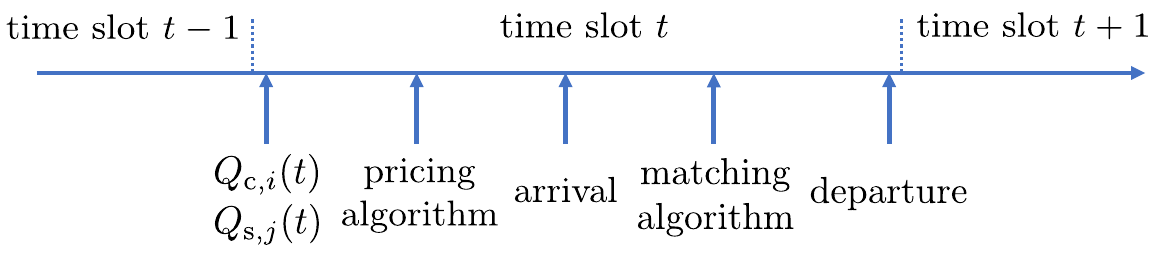}
    \caption{\small The timeline in each time slot.}
    \label{fig:timeline}
\end{figure}

The sequence of events occurring in each time slot is shown in Figure~\ref{fig:timeline}. At the beginning of time slot $t$, we first observe the lengths of all queues $Q_{\mathrm{c}, i}(t), Q_{\mathrm{s}, j}(t)$ for all $i,j$. Then the platform runs a pricing algorithm to decide the prices $p_{\mathrm{c}, i}(t), p_{\mathrm{s},j}(t),\forall i,j$. The arrival rates $\lambda_i, \mu_j, \forall i,j$ will be determined by these prices, and then the actual arrivals $A_{\mathrm{c},i}, A_{\mathrm{s},j}, \forall i,j$ occur. Next, the platform runs a matching algorithm and then those matched customers and servers leave the system.

At time slot $t$, the platform makes a profit of 
\begin{align*}
    \sum_i A_{\mathrm{c},i}(t) p_{\mathrm{c}, i}(t) - \sum_j A_{\mathrm{s},j}(t) p_{\mathrm{s}, j}(t)
    = \sum_i A_{\mathrm{c},i}(t) F_i(\lambda_i(t)) - \sum_j A_{\mathrm{s},j}(t) G_j(\mu_j(t)).
\end{align*}
Our goal is to design an online pricing and matching algorithms to maximize the profit of the platform over $T$ time slots without knowing the demand and supply functions, while keeping queue lengths of both customers and servers below a predetermined threshold, i.e., to maximize
\begin{align*}
    \sum_{t=1}^{T} \sum_i A_{\mathrm{c},i}(t) F_i(\lambda_i(t)) - \sum_j A_{\mathrm{s},j}(t) G_j(\mu_j(t)),
\end{align*}
while keeping
\begin{align*}
    \max_{t=1,\ldots,T} \max\{\max_i Q_{\mathrm{c}, i}(t), \max_j Q_{\mathrm{s}, j}(t)\}.
\end{align*}
below a predetermined threshold.

To quantify the performance of such an online algorithm, we compare it with the optimal solution to a fluid-based optimization problem~\cite{varma2023dynamic} below
\begin{align}
    \max_{\boldsymbol{\lambda}, \boldsymbol{\mu}, \boldsymbol{x}} & \sum_i \lambda_i F_i(\lambda_i) - \sum_j \mu_j G_j(\mu_j)\label{equ:fluid-opti-obj}\\
    \mathrm{s.t.} \qquad \lambda_i = &  \sum_{j: (i,j)\in {\cal E}} x_{i,j}, \quad \text{for all } i\in{\cal I} \label{equ:fluid-opti-constr1}\\
    \mu_j = &  \sum_{i: (i,j)\in {\cal E}} x_{i,j}, \quad \text{for all } j\in {\cal J}  \label{equ:fluid-opti-constr2}\\
    x_{i,j} \ge & 0, \quad \text{for all } (i,j)\in {\cal E}\label{equ:fluid-opti-constr3}\\
    \lambda_i, \mu_j \in & [0,1], \quad \text{for all } i\in{\cal I}, ~j\in {\cal J},\label{equ:fluid-opti-constr4}
\end{align}
where $\boldsymbol{x}\coloneqq (x_{i,j})_{(i,j)\in {\cal E}}$, $\boldsymbol{\lambda}\coloneqq (\lambda_i)_{i\in {\cal I}}$, and $\boldsymbol{\mu}\coloneqq (\mu_j)_{j\in {\cal J}}$.
Similarly to~\cite{varma2023dynamic}, we make the following assumption:
\begin{assumption}\label{assum:2}
    For all customer type $i$, the revenue function $r_i(\lambda_i)\coloneqq \lambda_i F_i(\lambda_i)$ is concave. For all server type $j$, the cost function $c_j(\mu_j)\coloneqq \mu_j G_j(\mu_j)$ is convex.
\end{assumption}
\noindent From Assumption~\ref{assum:2}, the objective~\eqref{equ:fluid-opti-obj} is concave.
Note that the constraints are all affine. Hence, the baseline optimization problem \eqref{equ:fluid-opti-obj}-\eqref{equ:fluid-opti-constr4} is concave and there exists a solution that achieves the optimal value, denoted as $(\boldsymbol{\lambda}^*, \boldsymbol{\mu}^*, \boldsymbol{x}^*)$. We make the following assumption on the optimal solution.
\begin{assumption}\label{assum:3}
    There exists a known positive number $a_{\min}\in (0,1)$ such that there exists an optimal solution $\boldsymbol{\lambda}^*, \boldsymbol{\mu}^*$ that satisfies $\lambda^*_i \ge a_{\min}$ and $\mu_j^* \ge a_{\min}$ for all $i$ and $j$.
\end{assumption}

In the fluid problem, $x_{i,j}$ is the asymptotic time-averaged expected number of matches on link $(i,j)$ at each time slot, i.e., $x_{i,j}\coloneqq \lim_{T\rightarrow \infty} \frac{1}{T} \sum_{t=1}^{T} \expt[X_{i,j}(t)]$, where $X_{i,j}(t)$ denote the number of matches at time $t$. $\lambda_i$ and $\mu_j$ are the asymptotic time average of $\expt[\lambda_i(t)]$ and $\expt[\mu_j(t)]$, respectively, i.e., $\lambda_i\coloneqq \lim_{T\rightarrow \infty} \frac{1}{T} \sum_{t=1}^{T} \expt[ \lambda_i(t)]$ and $\mu_j\coloneqq \lim_{T\rightarrow \infty} \frac{1}{T} \sum_{t=1}^{T} \expt[\mu_j(t)]$. Assume all the above limits exist.
We consider policies that make the queue mean rate stable~\cite{neely2022stochastic}, i.e., under the policy, for all $i,j$,
\begin{align}\label{equ:mean-rate-stable}
    \lim_{T\rightarrow \infty}\frac{1}{T} \expt [Q_{\mathrm{c},i}(T)] = 0,
    \quad \lim_{T\rightarrow \infty}\frac{1}{T} \expt [Q_{\mathrm{s},j}(T)] = 0.
\end{align}
\begin{proposition}\label{prop:1}
    Let Assumption~\ref{assum:1} and Assumption~\ref{assum:2} hold. Consider policies such that the limits $\lim_{T\rightarrow \infty} \frac{1}{T} \sum_{t=1}^{T} \expt[X_{i,j}(t)]$, $\lim_{T\rightarrow \infty} \frac{1}{T} \sum_{t=1}^{T} \expt[ \lambda_i(t)]$, and $\lim_{T\rightarrow \infty} \frac{1}{T} \sum_{t=1}^{T} \expt[\mu_j(t)]$ exist. Then any mean rate stable policy must satisfy constraints~\eqref{equ:fluid-opti-constr1} and \eqref{equ:fluid-opti-constr2}. Furthermore, the optimal value of the optimization problem~\eqref{equ:fluid-opti-obj}-\eqref{equ:fluid-opti-constr4} is an upper bound on the asymptotic time-averaged expected profit under any mean rate stable policy.
\end{proposition}
\noindent The proof of Proposition~\ref{prop:1} can be found in Appendix~\ref{app:proof-prop-1}.

For the profit maximization, we consider the following expected regret with respect to the fluid baseline:
\begin{align}\label{equ:regret-def}
    & \expt [R(T)] \coloneqq  \sum_{t=1}^{T}  \expt \left[ \biggl(\sum_i \lambda^*_i F_i(\lambda^*_i) - \sum_j \mu^*_j G_j(\mu^*_j) \biggr) 
    -  \biggl(  \sum_i A_{\mathrm{c},i}(t) F_i(\lambda_i(t)) - \sum_j A_{\mathrm{s},j}(t) G_j(\mu_j(t)) \biggr) \right]\nonumber\\
    = & \sum_{t=1}^{T} \left(  \biggl(\sum_i \lambda^*_i F_i(\lambda^*_i) - \sum_j \mu^*_j G_j(\mu^*_j) \biggr) 
    -   \expt \biggl[ \expt \Bigl[ \sum_i A_{\mathrm{c},i}(t) F_i(\lambda_i(t)) - \sum_j A_{\mathrm{s},j}(t) G_j(\mu_j(t)) \Bigl| \boldsymbol{\lambda}(t), \boldsymbol{\mu}(t) \Bigr. \Bigr] \biggr] \right)\nonumber\\
    = & \sum_{t=1}^{T} \Biggl( \biggl(\sum_i \lambda^*_i F_i(\lambda^*_i) - \sum_j \mu^*_j G_j(\mu^*_j) \biggr) 
    -   \expt \left[  \sum_i \lambda_i(t) F_i(\lambda_i(t)) - \sum_j \mu_j(t) G_j(\mu_j(t))  \right] \Biggr).
\end{align}

Assume that the functions $F_i$ and $G_j$ are {\em unknown}. The platform cannot directly control the arrival rates but they can control the prices of the customers and the servers and the matching algorithm. The platform can also observe the number of arrivals and the queue lengths of all the queues. At each time slot, the platform needs to design a price for each queue (customer or server) and to make a decision on how the customers and servers are matched. The objective is to minimize the regret $\expt [R(T)]$ while keeping the maximum queue length $\max_{t=1,\ldots,T} \max\{\max_i Q_{\mathrm{c}, i}(t), \max_j Q_{\mathrm{s}, j}(t)\}$ below a predetermined threshold.

\subsection{Notation}
We use bold symbols to represent vectors in $\mathbb{R}^{|{\cal E}|}$ or $\mathbb{R}^{I}$ or $\mathbb{R}^{J}$, such as $\boldsymbol{x}, \boldsymbol{\lambda}, \boldsymbol{\mu}$. We use subscript $i$ or $j$ to indicate an element of a vector, for example, $\lambda_i$ is the $i^{\mathrm{th}}$ element of $\boldsymbol{\lambda}$. Note that we view $\boldsymbol{x}= (x_{i,j})_{(i,j)\in {\cal E}}$ as a vector in $\mathbb{R}^{|{\cal E}|}$ rather than a matrix and the order of elements in $\boldsymbol{x}$ can be arbitrary as long as the order is fixed. We use subscript $\mathrm{c}$ to indicate the customer side and $\mathrm{s}$ for the server side. $g(T) = O(h(T))$ means that there exist positive constants $C$ and $T_0$ such that $g(T) \le C h(T)$ for all $T\ge T_0$. $g(T) = \Theta(h(T))$ means that there exists positive constants $C_1$, $C_2$ and $T_0$ such that $C_1 h(T) \le g(T) \le C_2 h(T)$ for all $T\ge T_0$. We use $\tilde{O}(\cdot)$ and $\tilde{\Theta}(\cdot)$ to suppress $\text{polylog}(T)$ factors.

\section{Summary of Our Main Results}
In this section, we provide an overview of the main results in this paper. We use a longest-queue-first matching algorithm, which is a discrete-time version of the matching algorithm in \cite{varma2023dynamic}. We propose a learning-based pricing algorithm, which combines gradient-free (zero-order) stochastic projected gradient ascent~\cite{agarwal2010optimal} with bisection search, illustrated in Section~\ref{sec:alg} in detail. The proposed pricing algorithm learns from samples and make pricing decisions simultaneously. Under the proposed matching and pricing algorithms, we establish the following profit-regret bound and queue-length bounds for sufficiently large $T$: 
\begin{align*}
    & \expt[R(T)] = \tilde{O}\Bigl(T^{\frac{5}{6}}\Bigr),\\
    & \frac{1}{T}\sum_{t=1}^{T} \expt \biggl[\sum_i Q_{\mathrm{c}, i}(t) + \sum_j Q_{\mathrm{s}, j}(t)\biggr] = \tilde{O}(T^{\frac{1}{2}}),\\
    & \max_{t=1,\ldots,T} \max\Bigl\{\max_i Q_{\mathrm{c}, i}(t), \max_j Q_{\mathrm{s}, j}(t)\Bigr\} = \tilde{O}\Bigl(T^{\frac{2}{3}}\Bigr).
\end{align*}
By changing the parameters of the proposed pricing algorithm, we can achieve the following tradeoff between regret and the maximum queue length: for any $\gamma \in (0,2/3]$, there exists a set of parameters such that the algorithm can achieve:
\begin{align*}
    \expt[R(T)] = \tilde{O}\Bigl(T^{1-\frac{\gamma}{4}}\Bigr),
    \qquad
    \max_{t=1,\ldots,T} \max\Bigl\{\max_i Q_{\mathrm{c}, i}(t), \max_j Q_{\mathrm{s}, j}(t)\Bigr\} = \tilde{O}\Bigl(T^{\gamma}\Bigr).
\end{align*}
We can also achieve the following tradeoff between regret and \emph{average} queue length: for any $\gamma \in (0,1/2]$, there exists a set of parameters such that the algorithm can achieve: 
\begin{align*}
    \expt[R(T)] = \tilde{O}\Bigl(T^{\max\{\frac{8-\gamma}{9}, 1-\gamma\}}\Bigr),
    \qquad
    \frac{1}{T}\sum_{t=1}^{T} \expt \biggl[\sum_i Q_{\mathrm{c}, i}(t) + \sum_j Q_{\mathrm{s}, j}(t)\biggr] = \tilde{O}\Bigl(T^{\gamma}\Bigr).
\end{align*}

\begin{figure}[htb]
    \centering
    \includegraphics[width=0.4\linewidth]{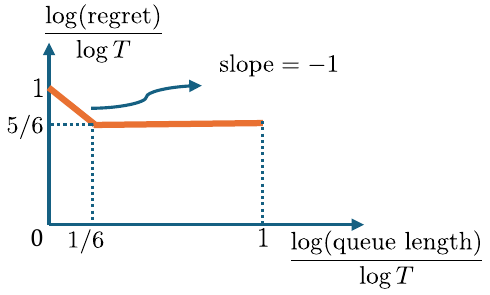}
    \caption{\small Tradeoff between regret and queue length.}
    \label{fig:tradeoff}
\end{figure}

Suppose we know a set of balanced arrival rates and associated price intervals in advance. Then we can achieve the following profit-regret bound and queue-length bounds:
\begin{align*}
    & \expt[R(T)] = \tilde{O}\Bigl(T^{\frac{5}{6}}\Bigr),\\
    & \frac{1}{T}\sum_{t=1}^{T} \expt \biggl[\sum_i Q_{\mathrm{c}, i}(t) + \sum_j Q_{\mathrm{s}, j}(t)\biggr] = \tilde{O}(T^{\frac{1}{6}}),\\
    & \max_{t=1,\ldots,T} \max\Bigl\{\max_i Q_{\mathrm{c}, i}(t), \max_j Q_{\mathrm{s}, j}(t)\Bigr\} = \tilde{O}\Bigl(T^{\frac{1}{6}}\Bigr).
\end{align*}
The tradeoff between regret and queue length can be improved as follows:
\begin{align*}
    & \expt[R(T)] = \tilde{O}\Bigl(T^{1 - \gamma}\Bigr),\\
    & \frac{1}{T}\sum_{t=1}^{T} \expt \biggl[\sum_i Q_{\mathrm{c}, i}(t) + \sum_j Q_{\mathrm{s}, j}(t)\biggr] = \tilde{O}(T^{\gamma}),\\
    & \max_{t=1,\ldots,T} \max\Bigl\{\max_i Q_{\mathrm{c}, i}(t), \max_j Q_{\mathrm{s}, j}(t)\Bigr\} = \tilde{O}\Bigl(T^{\gamma}\Bigr),
\end{align*}
for any $\gamma\in(0, 1/6]$.
The above results are shown in Figure~\ref{fig:tradeoff}. As the allowable queue length increases, the regret bound that we can achieve becomes better. However, if we allow the queue length to increase over $\tilde{\Theta}(T^{1/6})$, the regret bound cannot be further improved and remains $\tilde{O}(T^{5/6})$.

\section{Algorithm}
\label{sec:alg}
In this section, we present the proposed matching and pricing algorithms in detail. The matching algorithm is to match the customers on the demand-side queues with the servers on the supply-side queues, and the pricing algorithm is to learn and design the prices without knowing the demand and supply functions. 

\subsection{Matching Algorithm}

\begin{figure}[htbp]
    \centering
    \includegraphics[width=1.0\linewidth]{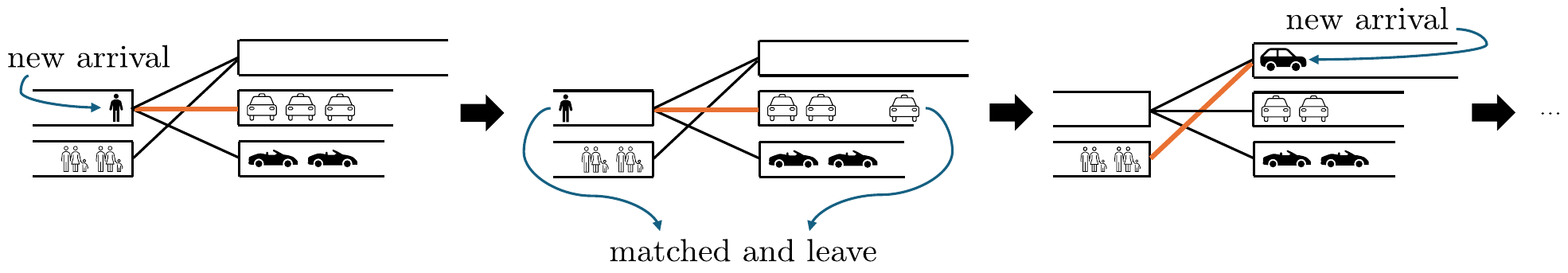}
    \caption{\small An example of longest-queue-first matching algorithm. At each time slot, for each queue, if there is a new arrival, it will be matched with one server (or customer) in the longest compatible queue on the other side. Then the matched pairs leave the system and we move on to check the next queue and repeat the same process.}
    \label{fig:matching-alg}
\end{figure}

We use a longest-queue-first matching algorithm as shown in Algorithm~\ref{alg:matching} and Figure~\ref{fig:matching-alg}, where we present the matching algorithm for one time slot and it is the same for all time slots.
In the algorithm, for each time slot, we iterate over each queue on both the customer and server sides. For each queue, we check whether there is a new arrival in the queue. If there is a new arrival (a customer or a server), the algorithm will match the arrival with one server (or customer) in the longest compatible queue on the other side, as shown in Figure~\ref{fig:matching-alg}. After matching, we decrease the queue length by one for both matched queues. Then we move on to the next queue and repeat the process.

\begin{algorithm}[htbp]
\caption{Matching algorithm }\label{alg:matching}
At each time slot $t$, without loss of generality, assume that the arrivals come in order $A_{\mathrm{c}, 1}(t), A_{\mathrm{c}, 2}(t),\ldots, A_{\mathrm{c}, I}(t), A_{\mathrm{s}, 1}(t), A_{\mathrm{s}, 2}(t),\ldots, A_{\mathrm{s}, J}(t)$.\\
{\textbf{Initialize:} $\tilde{Q}_{\mathrm{c}, i}\gets Q_{\mathrm{c}, i}(t)$ for all $i$ and $\tilde{Q}_{\mathrm{s}, j}\gets Q_{\mathrm{s}, j}(t)$ for all $j$.}\\
\For{\textnormal{$i= 1$ to $I$}}
{
    Observe $A_{\mathrm{c}, i}(t)$ and update $\tilde{Q}_{\mathrm{c}, i} \gets \tilde{Q}_{\mathrm{c}, i} + A_{\mathrm{c}, i}(t)$\;
    Let $j^*_i(t)\in\argmax_{j: (i,j)\in {\cal E}} \tilde{Q}_{\mathrm{s}, j}$ with ties broken arbitrarily\;
    \If{\textnormal{$A_{\mathrm{c}, i}(t)=1$ and $\sum_{j: (i,j)\in{\cal E}} \tilde{Q}_{\mathrm{s}, j} > 0$}}
    {
        Match the arrival $A_{\mathrm{c}, i}(t)$ with a server in queue $j^*_i(t)$\;
        Then update $\tilde{Q}_{\mathrm{c}, i} \gets \tilde{Q}_{\mathrm{c}, i} - 1$ and $\tilde{Q}_{\mathrm{s}, j^*_i(t)} \gets \tilde{Q}_{\mathrm{s}, j^*_i(t)} - 1$\;
    }
}
\For{\textnormal{$j= 1$ to $J$}}
{
    Observe $A_{\mathrm{s}, j}(t)$ and update $\tilde{Q}_{\mathrm{s}, j} \leftarrow \tilde{Q}_{\mathrm{s}, j} + A_{\mathrm{s}, j}(t)$\;
    Let $i^*_j(t)=\argmax_{i: (i,j)\in {\cal E}} \tilde{Q}_{\mathrm{c}, i}$ with ties broken arbitrarily\;
    \If{\textnormal{$A_{\mathrm{s}, j}(t)=1$ and $\sum_{i: (i,j)\in{\cal E}} \tilde{Q}_{\mathrm{c}, i} > 0$}}
    {
        Match the arrival $A_{\mathrm{s}, j}(t)$ with a customer in queue $i^*_j(t)$\;
        Then update $\tilde{Q}_{\mathrm{s}, j} \leftarrow \tilde{Q}_{\mathrm{s}, j} - 1$ and $\tilde{Q}_{\mathrm{c}, i^*_j(t)} \leftarrow \tilde{Q}_{\mathrm{c}, i^*_j(t)} - 1$\;
    }
}
\end{algorithm}

\subsection{Pricing Algorithm}

Since the functions $F_i,G_j$ are unknown and the arrival rates cannot be directly controlled, we cannot simply solve the optimization problem~\eqref{equ:fluid-opti-obj}-\eqref{equ:fluid-opti-constr4} and obtain the prices. In this section, we propose to combine the techniques of gradient-free (zero-order) stochastic projected gradient descent~\cite{agarwal2010optimal} and bisection search to solve the problem.

\begin{figure}[htbp]
    \centering
    \includegraphics[width=1.0\linewidth]{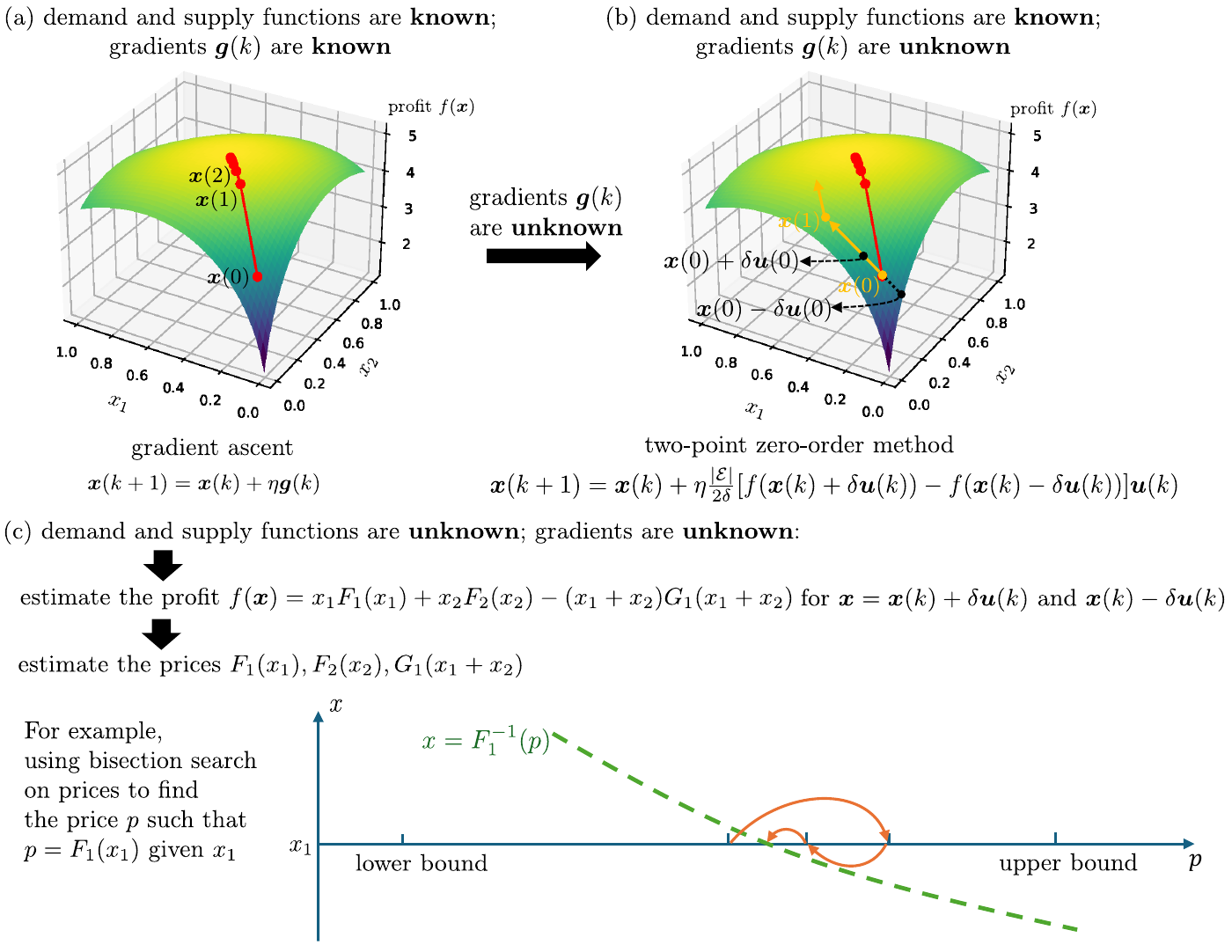}
    \caption{\small Pricing algorithm, an example with two customer queues connected to one server queue. Let $x_1$ and $x_2$ denote the arrival rates of the two customer queues, respectively. The arrival rate of the server queue should be $x_1+x_2$ because of the balance constraint. Let $\boldsymbol{x}(k)\coloneqq (x_1(k), x_2(k))$, where $k$ denote the iteration.}
    \label{fig:pricing-alg-map}
\end{figure}

Before presenting the details of our pricing algorithm, we begin by providing a simple example to understand the idea of our algorithm. We consider two customer queues, which are connected to one server queue. Let $x_1, x_2$ denote the arrival rates of customer queue 1 and customer queue 2, respectively. Then the arrival rate of the server queue should be equal to $x_1 + x_2$ due to the balance constraints~\eqref{equ:fluid-opti-constr1} and \eqref{equ:fluid-opti-constr2}.
Let $\boldsymbol{x}(k)\coloneqq (x_1(k), x_2(k))$, where $k$ denote the iteration.
Suppose the demand and supply functions are known and the gradients $\boldsymbol{g}(k)$ are known, i.e., we have first-order access to the profit function $f$. In this case, we can use the gradient ascent algorithm to solve the problem, as shown in Figure~\ref{fig:pricing-alg-map}(a), where $\eta$ is the step size. Suppose the demand and supply functions are known but the gradients are unknown, i.e., we only have zero-order access to the profit function $f$. Then we can use the two-point zero-order method~\cite{agarwal2010optimal} to solve the problem, as shown in Figure~\ref{fig:pricing-alg-map}(b). In the two point method, in each iteration $k$, we first sample a uniformly random direction $\boldsymbol{u}(k)$ ($\boldsymbol{u}(k)$ is a unit vector), and then calculate the profits using two points, $\boldsymbol{x}(k)+ \delta \boldsymbol{u}(k)$ and $\boldsymbol{x}(k)- \delta \boldsymbol{u}(k)$, where $\delta$ is small number. Then we can estimate the gradient and update $\boldsymbol{x}(k)$ using the formula in Figure~\ref{fig:pricing-alg-map}(b). However, if the demand and supply functions and the gradients are all unknown, simply using the two-point method does not work because we even do not have zero-order access to the profit function $f$. As shown in Figure~\ref{fig:pricing-alg-map}(c), for $\boldsymbol{x} = \boldsymbol{x}(k)+ \delta \boldsymbol{u}(k)$ or $\boldsymbol{x} = \boldsymbol{x}(k)- \delta \boldsymbol{u}(k)$, to estimate the profit $f(\boldsymbol{x})$, we need to estimate the prices of all queues $F_1(x_1), F_2(x_2), G_1(x_1+x_2)$. Take $F_1(x_1)$ as an example. Estimating $F_1(x_1)$ is equivalent to finding the solution to the equation $x_1=F_1^{-1}(p)$ given $x_1$.
We propose using bisection search on prices to estimate $p$ such that $p=F_1(x_1)$, as shown in Figure~\ref{fig:pricing-alg-map}(c). We first have the initial knowledge that the price $F_1(x_1)$ will be between some lower bound and some upper bound, which forms our first search interval. Then we set the price of the queue to be the midpoint of the search interval and run the system for multiple time slots to calculate average number of arrivals per time slot. Then we determine the next search interval by comparing it with $x_1$, as shown in Figure~\ref{fig:pricing-alg-map}(c). The bisection search stops when the accuracy is good enough. Note that during the process of running the system, we reject arrivals if the queue length is larger than or equal to a predetermined threshold, to control the queue length. After the bisection search, we obtain estimates of the profits $f(\boldsymbol{x}(k)+\delta \boldsymbol{u}(k))$ and $f(\boldsymbol{x}(k)-\delta \boldsymbol{u}(k))$.
Then we can substitute these estimates into the two-point method in Figure~\ref{fig:pricing-alg-map}(b), to get an estimate of the gradient and then update the arrival rate $\boldsymbol{x}(k)$.

\begin{figure}[htbp]
    \centering
    \includegraphics[width=0.9\linewidth]{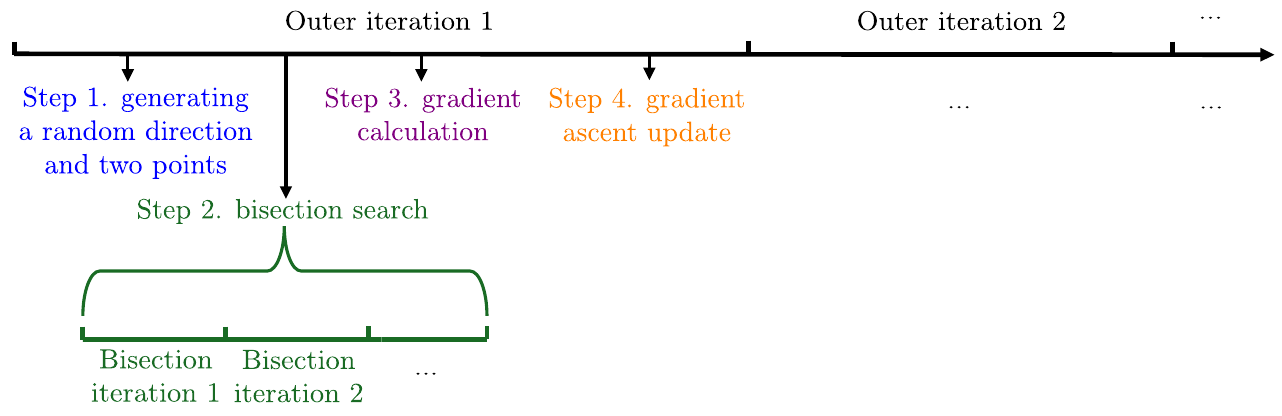}
    \caption{\small The structure of the pricing algorithm.}
    \label{fig:pricing-alg-struct}
\end{figure}

Now we present the details of the pricing algorithm in the general case. The structure of our pricing algorithm is shown in Figure~\ref{fig:pricing-alg-struct}. Before presenting the algorithm, we first define some notations and a shrunk feasible set.
Let ${\cal E}_{\mathrm{c},i}\coloneqq \{j | (i,j)\in {\cal E}\}$ and ${\cal E}_{\mathrm{s},j}\coloneqq \{i | (i,j)\in {\cal E}\}$ denote the sets of all queues that are connected to customer queue $i$ and server queue $j$, respectively. By substituting $\lambda_i, \mu_j$ and by Assumption~\ref{assum:3}, the optimization problem~\eqref{equ:fluid-opti-obj}-\eqref{equ:fluid-opti-constr4} can be equivalently rewritten as
\begin{align}
    \max_{\boldsymbol{x}} f(\boldsymbol{x}) \coloneqq & \sum_i  \Biggl(\sum_{j\in {\cal E}_{\mathrm{c},i}} x_{i,j} \Biggr) F_i\Biggl(\sum_{j\in {\cal E}_{\mathrm{c},i}} x_{i,j}\Biggr) 
    - \sum_j \Biggl(\sum_{i\in {\cal E}_{\mathrm{s},j}} x_{i,j}\Biggr) G_j\Biggl(\sum_{i\in {\cal E}_{\mathrm{s},j}} x_{i,j}\Biggr) \label{equ:eqv-opti-obj}\\
    \mathrm{s.t.}  \sum_{j\in {\cal E}_{\mathrm{c},i}} x_{i,j} \in & [a_{\min}, 1], \quad i=1,\ldots,I \label{equ:eqv-opti-constr1}\\
    \sum_{i\in {\cal E}_{\mathrm{s},j}} x_{i,j} \in & [a_{\min}, 1], \quad j=1,\ldots,J  \label{equ:eqv-opti-constr2}\\
    x_{i,j} \ge & 0, \quad (i,j)\in {\cal E}\label{equ:eqv-opti-constr3}
\end{align}
Note that we view $\boldsymbol{x}\coloneqq (x_{i,j})_{(i,j)\in {\cal E}}$ as a vector in $\mathbb{R}^{|{\cal E}|}$.
Note that the problem \eqref{equ:eqv-opti-obj}-\eqref{equ:eqv-opti-constr3} is also concave.
Let ${\cal D}\subseteq \mathbb{R}^{|{\cal E}|}$ denote the feasible set of the problem~\eqref{equ:eqv-opti-obj}-\eqref{equ:eqv-opti-constr3}, i.e.,
\begin{align}\label{equ:def-D}
    {\cal D} \coloneqq \left\{ \boldsymbol{x} \left| x_{i,j}\ge 0 \text{ for all } i,j,
    \sum_{j\in {\cal E}_{\mathrm{c},i}} x_{i,j} \in [a_{\min}, 1] \text{ for all } i,
    \sum_{i\in {\cal E}_{\mathrm{s},j}} x_{i,j} \in [a_{\min}, 1] \text{ for all } j \right.\right\}.
\end{align}
Let $N_{i,j}\coloneqq \max\{|{\cal E}_{\mathrm{c},i}|, |{\cal E}_{\mathrm{s},j}|\}$ denote the maximum cardinality of the sets ${\cal E}_{\mathrm{c},i}$ and ${\cal E}_{\mathrm{s},j}$. Define a shrunk set of ${\cal D}$ with a parameter $\delta$ as follows:
\begin{align}\label{equ:def-D'}
    & {\cal D}' \coloneqq \Biggl\{ \boldsymbol{x} \Biggl|
    \text{ for all } i,j, 
    x_{i,j}-\frac{a_{\min}+1}{2N_{i,j}} \ge - \biggl(1-\frac{\delta}{r}\biggr) \frac{a_{\min}+1}{2N_{i,j}},\nonumber\\
    & \sum_{j'\in {\cal E}_{\mathrm{c},i}} \biggl(x_{i,j'} - \frac{a_{\min}+1}{2N_{i,j'}}\biggr)\in 
    \biggl[ -\biggl(1-\frac{\delta}{r}\biggr)\biggl( \sum_{j'\in {\cal E}_{\mathrm{c},i}} \frac{a_{\min}+1}{2N_{i,j'}} - a_{\min} \biggr), 
    \biggl(1-\frac{\delta}{r}\biggr)  \biggl(  1 - \sum_{j'\in {\cal E}_{\mathrm{c},i}} \frac{a_{\min}+1}{2N_{i,j'}} \biggr) \biggr],\nonumber\\
    & \sum_{i'\in {\cal E}_{\mathrm{s},j}} \biggl(x_{i',j} - \frac{a_{\min}+1}{2N_{i',j}}\biggr)\in 
    \biggl[ -\biggl(1-\frac{\delta}{r}\biggr)\biggl( \sum_{i'\in {\cal E}_{\mathrm{s},j}} \frac{a_{\min}+1}{2N_{i',j}} - a_{\min} \biggr), 
    \biggl(1-\frac{\delta}{r}\biggr)  \biggl(  1 - \sum_{i'\in {\cal E}_{\mathrm{s},j}} \frac{a_{\min}+1}{2N_{i',j}} \biggr) \biggr]
    \Biggr.\Biggr\},
\end{align}
where
\begin{align}\label{equ:def-r}
    r \coloneqq & \min_{i,j} \Biggl\{
    \frac{1+a_{\min}}{2N_{i,j}},
    \frac{1}{|{\cal E}_{\mathrm{c}, i}|}\biggl( 1 -  \sum_{j'\in {\cal E}_{\mathrm{c},i}} \frac{a_{\min}+1}{2N_{i,j'}} \biggr),
    \frac{1}{|{\cal E}_{\mathrm{s}, j}|}\biggl( 1 -  \sum_{i'\in {\cal E}_{\mathrm{s},j}} \frac{a_{\min}+1}{2N_{i',j}} \biggr), \nonumber\\
    & \qquad \frac{1}{|{\cal E}_{\mathrm{c}, i}|} \biggl( \sum_{j'\in {\cal E}_{\mathrm{c},i}} \frac{a_{\min}+1}{2N_{i,j'}} - a_{\min} \biggr),
    \frac{1}{|{\cal E}_{\mathrm{s}, j}|} \biggl( \sum_{i'\in {\cal E}_{\mathrm{s},j}} \frac{a_{\min}+1}{2N_{i',j}} - a_{\min} \biggr)
    \Biggr\},
\end{align}
and $\delta \in (0 ,r)$. With the definition of ${\cal D}'$, we can show (Lemma~\ref{lemma:shrunk-set} in Appendix~\ref{app:theo:1-pre-lemma}) that $\boldsymbol{x} + \delta \boldsymbol{u} \in {\cal D}$ for any $\boldsymbol{x}\in {\cal D}'$ and any vector $\boldsymbol{u}$ in the unit ball.
To ensure that the definition of ${\cal D}'$ is valid, we need the following assumption.
\begin{assumption}\label{assum:4}
    $\sum_{j\in {\cal E}_{\mathrm{c},i}} \frac{a_{\min}+1}{2N_{i,j}} - a_{\min} > 0$ for all $i$ and $\sum_{i\in {\cal E}_{\mathrm{s},j}} \frac{a_{\min}+1}{2N_{i,j}} - a_{\min} > 0$ for all $j$.
\end{assumption}
As long as $a_{\min}$ is small enough, Assumption~\ref{assum:4} holds. Since $N_{i,j}\ge |{\cal E}_{\mathrm{c}, i}|$ and $a_{\min}<1$, we have $1 - \sum_{j\in {\cal E}_{\mathrm{c},i}} \frac{a_{\min}+1}{2N_{i,j}} \ge 1-\frac{a_{\min}+1}{2}>0$ for all $i$. Also, by Assumption~\ref{assum:4}, $\sum_{j\in {\cal E}_{\mathrm{c},i}} \frac{a_{\min}+1}{2N_{i,j}} - a_{\min} > 0$ for all $i$. Similarly, we have $1 - \sum_{i\in {\cal E}_{\mathrm{s},j}} \frac{a_{\min}+1}{2N_{i,j}} > 0$ and $\sum_{i\in {\cal E}_{\mathrm{s},j}} \frac{a_{\min}+1}{2N_{i,j}} - a_{\min} > 0$ for all $j$. Hence, $r>0$ and the definition of ${\cal D}'$ is valid. Let $\boldsymbol{x}_{\mathrm{ctr}} \coloneqq (\frac{a_{\min}+1}{2N_{i,j}})_{(i,j)\in{\cal E}}$. We can see that $\boldsymbol{x}_{\mathrm{ctr}} \in {\cal D}'$, which can be used as the initial point for the pricing algorithm. $\boldsymbol{x}_{\mathrm{ctr}}$ can be thought of as the ``center'' of the set ${\cal D}'$.

\begin{algorithm}[p]
\small
\caption{\small Pricing algorithm }\label{alg:pricing}
\textbf{Initialize:} Choose an exploration parameter $\delta\in(0,r)$.
Choose a step size $\eta\in (0,1)$. 
Choose an accuracy parameter $\epsilon\in (0, 1/e)$.
Choose a queue length threshold $q^{\mathrm{th}}$.
Choose $\boldsymbol{x}(1) = \boldsymbol{x}_{\mathrm{ctr}} \in {\cal D}'$ as the initial point.
Define $e_{\mathrm{c},i}$ and $e_{\mathrm{s},j}$ according to \eqref{equ:def-e-c} and \eqref{equ:def-e-s}.
Define $N\coloneqq \left\lceil\frac{\beta\ln (1/\epsilon)}{\epsilon^2}\right\rceil$ and $M \coloneqq \left\lceil \log_2 \frac{1}{\epsilon} \right\rceil$, where $\beta>0$ is a constant.\\
Time step counter $t\gets 1$;\\
Outer iteration counter $k\gets 1$;\\
\Repeat{$t > T$}
{
    \nonl \blue{//\textbf{ Step 1.~generate a random direction and two points}}\\
    Choose a unit vector $\boldsymbol{u}(k)\in \mathbb{R}^{|\cal E|}$ uniformly at random, i.e., $\|\boldsymbol{u}(k)\|_2 = 1$;\\
    Let $x^{+}_{i,j}(k)\coloneqq (\boldsymbol{x}(k) + \delta \boldsymbol{u}(k))_{i,j}$ and $x^{-}_{i,j}(k)\coloneqq (\boldsymbol{x}(k) - \delta \boldsymbol{u}(k))_{i,j}$\label{line:alg-pricing-x};\\
    Let $\lambda_i^{+}(k)\coloneqq \sum_{j\in {\cal E}_{\mathrm{c},i}} x^{+}_{i,j}(k)$,
    $\lambda_i^{-}(k)\coloneqq \sum_{j\in {\cal E}_{\mathrm{c},i}} x^{-}_{i,j}(k)$,
    $\mu_j^{+}(k)\coloneqq \sum_{i\in {\cal E}_{\mathrm{s},j}} x^{+}_{i,j}(k)$,
    $\mu_j^{-}(k)\coloneqq \sum_{i\in {\cal E}_{\mathrm{s},j}} x^{-}_{i,j}(k)$\label{line:alg-pricing-lambda-mu};\\
    Let $\boldsymbol{\lambda}^{+}(k)$ be a vector of $\lambda_i^{+}(k),i=1\ldots,I$. Define similarly $\boldsymbol{\lambda}^{-}(k)$, $\boldsymbol{\mu}^{+}(k)$, and $\boldsymbol{\mu}^{-}(k)$.\\
    \nonl \darkgreen{//\textbf{ Step 2.~bisection search to approximate $F_i(\lambda_i^{+}(k)), F_i(\lambda_i^{-}(k)), G_j(\mu_j^{+}(k)), G_j(\mu_j^{-}(k))$.}}\\
    \If{$k=1$\label{line:alg-pricing-binary-start}}
    {
        For all $i$, let $\underline{p}^{+}_{\mathrm{c}, i}(k,1)=p_{\mathrm c,i,\min}$,
        $\bar{p}^{+}_{\mathrm{c},i}(k,1)=p_{\mathrm{c},i,\max}$\label{line:alg-pricing-bound-bisection-1-i};\\
        For all $j$, let $\underline{p}^{+}_{\mathrm{s}, j}(k,1)=p_{\mathrm s,j,\min}$,
        $\bar{p}^{+}_{\mathrm{s},j}(k,1)=p_{\mathrm{s},j,\max}$\label{line:alg-pricing-bound-bisection-1-j};\\
        Let $\tilde{q}^{\mathrm{th}} = +\infty$;
        \tcp{We do not reject arrivals for $k=1$.}
    }  
    \Else
    {
        For all $i$, let $\underline{p}^{+}_{\mathrm{c}, i}(k,1)=p_{\mathrm{c},i}^+(k-1,M) - e_{\mathrm{c},i}$,
        $\bar{p}^{+}_{\mathrm{c},i}(k,1)=p_{\mathrm{c},i}^+(k-1, M) + e_{\mathrm{c},i}$\label{line:alg-pricing-bound-bisection-i};\\
        For all $j$, let $\underline{p}^{+}_{\mathrm{s}, j}(k,1)=p_{\mathrm{s},j}^+(k-1, M) - e_{\mathrm{s},j}$,
            $\bar{p}^{+}_{\mathrm{s},j}(k,1)=p_{\mathrm{s},j}^+(k-1, M) + e_{\mathrm{s},j}$\label{line:alg-pricing-bound-bisection-j};\\
        Let $\tilde{q}^{\mathrm{th}} = q^{\mathrm{th}}$;
    }
    Let $\underline{\boldsymbol{p}}^{+}_{\mathrm{c}}(k,1)$ be a vector of $\underline{p}^{+}_{\mathrm{c}, i}(k,1), i=1\ldots,I$. Define similarly $\bar{\boldsymbol{p}}^{+}_{\mathrm{c}}(k,1)$, $\underline{\boldsymbol{p}}^{+}_{\mathrm{s}}(k,1)$, $\bar{\boldsymbol{p}}^{+}_{\mathrm{s}}(k,1)$.\\
    $t$, $\boldsymbol{p}_{\mathrm{c}}^{+} (k,M)$, $\boldsymbol{p}_{\mathrm{s}}^{+} (k,M)=$ Bisection$\left(\boldsymbol{\lambda}^{+}(k), \boldsymbol{\mu}^{+}(k), \underline{\boldsymbol{p}}^{+}_{\mathrm{c}}(k,1), \bar{\boldsymbol{p}}^{+}_{\mathrm{c}}(k,1), \underline{\boldsymbol{p}}^{+}_{\mathrm{s}}(k,1), \bar{\boldsymbol{p}}^{+}_{\mathrm{s}}(k,1), \tilde{q}^{\mathrm{th}}, t, M, N, \epsilon\right)$\label{line:alg-pricing-binary-end};\\
    Do Line~\ref{line:alg-pricing-binary-start}-\ref{line:alg-pricing-binary-end} for $\boldsymbol{\lambda}^{-}(k), \boldsymbol{\mu}^{-}(k)$. Denote the counterparts of the prices by $\underline{\boldsymbol{p}}^{-}_{\mathrm{c}}(k,1)$, $\bar{\boldsymbol{p}}^{-}_{\mathrm{c}}(k,1)$,
    $\boldsymbol{p}^{-}_{\mathrm{c}}(k,M)$, $\underline{\boldsymbol{p}}^{-}_{\mathrm{s}}(k,1)$, $\bar{\boldsymbol{p}}^{-}_{\mathrm{s}}(k,1)$, $\boldsymbol{p}^{-}_{\mathrm{s}}(k,M)$\label{line:alg-pricing-bisection}.\\
    \nonl \violet{// \textbf{Step 3.~gradient calculation}}\\
    Let $\hat{\boldsymbol{g}}(k)=\frac{|{\cal E}|}{2\delta} \biggl[ \left( \sum_{i=1}^{I} \lambda_i^+(k) p_{\mathrm{c},i}^+(k, M) - \sum_{j=1}^{J} \mu_j^+(k) p_{\mathrm{s},j}^+(k, M) \right)$\label{line:alg-pricing-gradient}\\
    \nonl \qquad \qquad \qquad \quad$ - \left( \sum_{i=1}^{I} \lambda_i^-(k) p_{\mathrm{c},i}^-(k, M) - \sum_{j=1}^{J} \mu_j^-(k) p_{\mathrm{s},j}^-(k, M)
    \right) \biggr] \boldsymbol{u}(k)$;\\
    \nonl \orange{// \textbf{Step 4.~gradient ascent update}}\\
    Projected Gradient Ascent: $\boldsymbol{x}(k+1) = \Pi_{{\cal D}'}(\boldsymbol{x}(k)+\eta \hat{\boldsymbol{g}}(k) )$;\label{line:alg-pricing-pga}\\
    $k\gets k + 1$;
}
\end{algorithm}

The details of the proposed pricing algorithm are shown in Algorithm~\ref{alg:pricing} and Algorithm~\ref{alg:bisection}. 
In the notation in Algorithms~\ref{alg:pricing} and~\ref{alg:bisection}, $(k)$ denotes the $k^{\mathrm{th}}$ outer iteration; $(k,m)$ denotes the $m^{\mathrm{th}}$ bisection iteration in the $k^{\mathrm{th}}$ outer iteration. 
As shown in Algorithm~\ref{alg:pricing}, we first choose $\delta, \eta, \epsilon$, and $q^{\mathrm{th}}$. $\delta\in (0,r)$ is an exploration parameter for estimating the gradient, and it is used to construct the set ${\cal D}'$. We choose $\boldsymbol{x}(1)=\boldsymbol{x}_{\mathrm{ctr}}\in {\cal D}'$ as the initial point of the algorithm.
$\eta\in (0,1)$ is the step size for the gradient ascent step. $\epsilon\in (0,1/e)$ is the accuracy for the bisection algorithm. The threshold $q^{\mathrm{th}} > 0$ is used to control queue length. We also define $e_{\mathrm{c},i}$ and $e_{\mathrm{s},j}$, which can be later proved to be the upper bounds of the change of the prices in one outer iteration. Note that the parameters $\delta, \eta, \epsilon, q^{\mathrm{th}}$ are functions of the total number of time steps $T$. $\delta, \eta, \epsilon$ are decreasing in $T$, and $q^{\mathrm{th}}$ is increasing in $T$, and they can be optimized later to obtain a sublinear regret. In each outer iteration $k$, we have four steps:
\begin{itemize}[leftmargin=*]
    \item
    \textbf{Step 1}:
    We first generate a random direction, i.e., a unit vector $\boldsymbol{u}(k)\in \mathbb{R}^{|\cal E|}$ that is uniformly sampled from the unit sphere, and then change $\boldsymbol{x}(k)$ in the direction of $\boldsymbol{u}(k)$ and also in the opposite direction of $\boldsymbol{u}(k)$, generating $x^{+}_{i,j}(k)$ and $x^{-}_{i,j}(k)$, as shown in Line~\ref{line:alg-pricing-x} of Algorithm~\ref{alg:pricing}.
    We can then calculate the arrival rates $\lambda^{+}_i(k), \mu^{+}_j(k)$ and $\lambda^{-}_i(k), \mu^{-}_j(k)$ as shown in Line~\ref{line:alg-pricing-lambda-mu} of Algorithm~\ref{alg:pricing}.
    \item
    \textbf{Step 2}:
    Next, we find the prices that approximate $F_i(\lambda^{+}_i(k)), G_j(\mu^{+}_j(k))$ and $F_i(\lambda^{-}_i(k)), G_j(\mu^{-}_j(k))$ using the bisection search shown in Algorithm~\ref{alg:bisection}, which will be illustrated in detail later in this section. Here, we specify the initial search intervals for the bisection search.
    If $k=1$, we start with search intervals $[p_{\mathrm c,i,\min}, p_{\mathrm{c},i,\max}]$ and $[p_{\mathrm s,j,\min}, p_{\mathrm{s},j,\max}]$, as shown in Line~\ref{line:alg-pricing-bound-bisection-1-i}-\ref{line:alg-pricing-bound-bisection-1-j} of Algorithm~\ref{alg:pricing}. Note that we do not reject arrivals in the first iteration because we want to quickly learn the price that can keep the queues balanced.
    If $k>1$, we use $p_{\mathrm{c},i}^+(k-1, M_{k-1}^+)$, $p_{\mathrm{s},j}^+(k-1, M_{k-1}^+), p_{\mathrm{c},i}^-(k-1, M_{k-1}^-)$, $p_{\mathrm{s},j}^-(k-1, M_{k-1}^-), e_{\mathrm{c},i}, e_{\mathrm{s},j}$ to construct the lower and upper bounds for the bisection algorithm, as shown in Line~\ref{line:alg-pricing-bound-bisection-i}-\ref{line:alg-pricing-bound-bisection-j} of Algorithm~\ref{alg:pricing},  where
    $e_{\mathrm{c},i} = \Theta\left(\frac{\eta \epsilon}{\delta} + \eta + \delta + \epsilon\right)$ and $e_{\mathrm{s},j} = \Theta\left(\frac{\eta \epsilon}{\delta} + \eta + \delta + \epsilon\right)$ and the exact expressions can be found in Appendix~\ref{app:exact-expression}.
    \item 
    \textbf{Step 3}:
    The bisection algorithm will output approximated prices $\boldsymbol{p}_{\mathrm{c}}^{+} (k,M^+_k)$, $\boldsymbol{p}_{\mathrm{s}}^{+} (k,M^+_k)$, $\boldsymbol{p}_{\mathrm{c}}^{-} (k,M^+_k)$, $\boldsymbol{p}_{\mathrm{s}}^{-} (k,M^+_k)$ that correspond to $\boldsymbol{\lambda}^{+}(k), \boldsymbol{\mu}^{+}(k), \boldsymbol{\lambda}^{-}(k), \boldsymbol{\mu}^{-}(k)$, respectively, as shown in Line~\ref{line:alg-pricing-binary-end}-\ref{line:alg-pricing-bisection} of Algorithm~\ref{alg:pricing}. Then we can use these arrival rates and prices to calculate an estimate $\hat{\boldsymbol{g}}(k)$ of the gradient of the objective function $f(\boldsymbol{x})$ in \eqref{equ:eqv-opti-obj}, as shown in Line~\ref{line:alg-pricing-gradient} of Algorithm~\ref{alg:pricing}.
    \item 
    \textbf{Step 4}:
    With this gradient estimate, we do a projected gradient ascent as shown in Line~\ref{line:alg-pricing-pga} of Algorithm~\ref{alg:pricing}. The algorithm stops until the number of time steps $t\ge T$.
\end{itemize}

\begin{algorithm}[htbp]
\small
\caption{\small Bisection$\left(\boldsymbol{\lambda}^{+/-}(k), \boldsymbol{\mu}^{+/-}(k), \underline{\boldsymbol{p}}^{+/-}_{\mathrm{c}}(k,1), \bar{\boldsymbol{p}}^{+/-}_{\mathrm{c}}(k,1), \underline{\boldsymbol{p}}^{+/-}_{\mathrm{s}}(k,1), \bar{\boldsymbol{p}}^{+/-}_{\mathrm{s}}(k,1), \tilde{q}^{\mathrm{th}}, t, M, N, \epsilon\right)$ }\label{alg:bisection}
\For{\textnormal{$m=1$ to $M$}}
{
    Let $p^{+/-}_{\mathrm{c}, i} (k,m) = \frac{1}{2} \left(\underline{p}^{+/-}_{\mathrm{c}, i}(k,m) + \bar{p}^{+/-}_{\mathrm{c}, i}(k,m)\right)$ for all $i$; \label{line:alg-bisection-middle-i}\\
    Let $p^{+/-}_{\mathrm{s}, j} (k,m) = \frac{1}{2} \left(\underline{p}^{+/-}_{\mathrm{s}, j}(k,m) + \bar{p}^{+/-}_{\mathrm{s}, j}(k,m)\right)$ for all $j$; \label{line:alg-bisection-middle-j}\\
    Let $n_{\mathrm{c},i}(k,m)=0$ for all $i$ and $n_{\mathrm{s},j}(k,m)=0$ for all $j$; \label{line:alg-bisection-run-price-start}\\
    \Repeat{\textnormal{for all queues $i$, $n_{\mathrm{c},i}(k,m)\ge N$, i.e., the price $p^{+/-}_{\mathrm{c}, i} (k,m)$ is run for at least $N$ times and for all queues $j$, $n_{\mathrm{s},j}(k,m)\ge N$, i.e., the price $p^{+/-}_{\mathrm{s}, j} (k,m)$ is run for at least $N$ times}}
    {
        \For{\textnormal{$i=1$ to $I$}\label{line:alg-bisection-threshold-start}}
        {
            \lIf{$Q_{\mathrm{c},i}(t)\ge \tilde{q}^{\mathrm{th}}$}{set price $p_{\mathrm{c},i,\max}$ for queue $i$}
            \lElse{set price $p^{+/-}_{\mathrm{c}, i} (k,m)$ for queue $i$
            and $n_{\mathrm{c},i}(k,m) \gets n_{\mathrm{c},i}(k,m) + 1$}
        }
        \For{\textnormal{$j=1$ to $J$}}
        {
            \lIf{$Q_{\mathrm{s},j}(t)\ge \tilde{q}^{\mathrm{th}}$}{set price $p_{\mathrm{s},j,\min}$ for queue $j$}
            \lElse{set price $p^{+/-}_{\mathrm{s}, j} (k,m)$ for queue $j$
            and $n_{\mathrm{s},j}(k,m) \gets n_{\mathrm{s},j}(k,m) + 1$}
        }\label{line:alg-bisection-threshold-end}
        Use the above set of prices to run one step of the system;\\
        $t \leftarrow t + 1$;\\
        Terminate the algorithm when $t > T$;
    }\label{line:alg-bisection-run-price-end}
    Let $t^{+/-}_{\mathrm{c}, i}(k,m,n)$ denote the time slot when the price $p^{+/-}_{\mathrm{c}, i} (k,m)$ is run for the $n^{\mathrm{th}}$ time for the customer-side queue $i$; Let $t^{+/-}_{\mathrm{s}, j}(k,m,n)$ denote the time slot when the price $p^{+/-}_{\mathrm{s}, j} (k,m)$ is run for the $n^{\mathrm{th}}$ time for the server-side queue $j$;\\
    \For{\textnormal{$i=1$ to $I$}}
    {
        Let $\hat{\lambda}^{+/-}_i  (k,m)=\frac{1}{N}\sum_{n=1}^{N}A_{\mathrm{c},i}(t^{+/-}_{\mathrm{c}, i}(k,m,n))$;
        \tcp{sample average}\label{line:alg-bisection-estimate-i}
        \lIf{$\hat{\lambda}^{+/-}_i  (k,m) > \lambda^{+/-}_{i} (k)$\label{line:alg-bisection-update-i-start}}
        {$\underline{p}^{+/-}_{\mathrm{c}, i}(k,m+1) = p^{+/-}_{\mathrm{c}, i}  (k,m), 
        ~ \bar{p}^{+/-}_{\mathrm{c},i}(k,m+1) = \bar{p}^{+/-}_{\mathrm{c},i}(k,m)$}
        \lElse{$\underline{p}^{+/-}_{\mathrm{c}, i}(k,m+1) = \underline{p}  ^{+/-}_{\mathrm{c}, i}(k,m),
        ~ \bar{p}^{+/-}_{\mathrm{c},i}(k,m+1) = p^{+/-}_{\mathrm{c}, i}  (k,m)$}
        \label{line:alg-bisection-update-i-end}
    }

    \For{\textnormal{$j=1$ to $J$}}
    {
        Let $\hat{\mu}^{+/-}_j  (k,m)=\frac{1}{N}\sum_{n=1}^{N} A_{\mathrm{s},j}(t^{+/-}_{\mathrm{s}, j}(k,m,n))$;\tcp{sample average}\label{line:alg-bisection-estimate-j}
        \lIf{$\hat{\mu}^{+/-}_j  (k,m) > \mu^{+/-}_{j} (k)$\label{line:alg-bisection-update-j-start}}{$\underline{p}^{+/-}_{\mathrm{s},j}(k,m+1) = \underline{p}  ^{+/-}_{\mathrm{s},j}(k,m), 
        ~ \bar{p}^{+/-}_{\mathrm{s}, j}(k,m+1) = p^{+/-}_{\mathrm{s}, j}  (k,m)$}
        \lElse{$\underline{p}^{+/-}_{\mathrm{s}, j}(k,m+1) = p^{+/-}_{\mathrm{s}, j}  (k,m),  
        ~ \bar{p}^{+/-}_{\mathrm{s},j}(k,m+1) = \bar{p}^{+/-}_{\mathrm{s}, j}(k,m)$}
        \label{line:alg-bisection-update-j-end}
    }      
}
Let $\boldsymbol{p}^{+/-}_{\mathrm{c}} (k,M)$ be a vector of $p^{+/-}_{\mathrm{c}, i} (k,M),i=1,\ldots,I$ and define similarly $\boldsymbol{p}^{+/-}_{\mathrm{s}} (k,M)$.\\
\KwRet{$t$, $\boldsymbol{p}^{+/-}_{\mathrm{c}} (k,M)$, $\boldsymbol{p}^{+/-}_{\mathrm{s}} (k,M)$} 
\end{algorithm}

The proposed bisection algorithm is shown in Algorithm~\ref{alg:bisection}, which is to find the prices $\boldsymbol{p}^{+/-}_{\mathrm{c}}, \boldsymbol{p}^{+/-}_{\mathrm{s}}$ such that $\lambda_i^{+/-}(k)\approx F_i^{-1}(p^{+/-}_{\mathrm{c},i})$ and $\mu_j^{+/-}(k) \approx G_j^{-1}(p^{+/-}_{\mathrm{s},j})$ for some given target arrival rates $\lambda^{+/-}_i(k), \mu^{+/-}_j(k)$ with initial intervals $[\underline{p}^{+/-}_{\mathrm{c},i}(k,1), \bar{p}^{+/-}_{\mathrm{c}_i}(k,1)]$ and $[\underline{p}^{+/-}_{\mathrm{s},j}(k,1), \bar{p}^{+/-}_{\mathrm{s},j}(k,1)]$ for all $i,j$.
In each bisection iteration $m$, we first calculate the midpoints of the intervals, $p^{+/-}_{\mathrm{c}, i} (k,m)$ and $p^{+/-}_{\mathrm{s}, j} (k,m)$ for all $i,j$ just like the standard bisection search, as shown in Line~\ref{line:alg-bisection-middle-i}-\ref{line:alg-bisection-middle-j} of Algorithm~\ref{alg:bisection}. For each queue,  if the queue length is less than $\tilde{q}^{\mathrm{th}}$, we will set the price with $p^{+/-}_{\mathrm{c}, i} (k,m)$ or $p^{+/-}_{\mathrm{s}, j} (k,m)$ that was just calculated, otherwise we will temporarily reject new arrivals into the queue to control the queue length by using the prices $p_{\mathrm{c},i,\max}$ or $p_{\mathrm{s},j,\min}$. We will use this set of prices to run the system for a number of time slots such that every price has at least $N$ i.i.d. samples, as shown in Line~\ref{line:alg-bisection-run-price-start}-\ref{line:alg-bisection-run-price-end} of Algorithm~\ref{alg:bisection}. Next, with the samples we obtained, we can calculate an estimate of the arrival rate corresponding to the price we set for each queue, $\hat{\lambda}^{+/-}_i (k,m)$ or $\hat{\mu}^{+/-}_j (k,m)$, as shown in Line~\ref{line:alg-bisection-estimate-i} and Line~\ref{line:alg-bisection-estimate-j} of Algorithm~\ref{alg:bisection}. Then we can update the lower and upper bounds of the searching intervals by comparing $\hat{\lambda}^{+/-}_i (k,m)$, $\hat{\mu}^{+/-}_j (k,m)$ with $\lambda^{+/-}_i (k,m)$, $\mu^{+/-}_j (k,m)$, as shown in Line~\ref{line:alg-bisection-update-i-start}-\ref{line:alg-bisection-update-i-end} and Line~\ref{line:alg-bisection-update-j-start}-\ref{line:alg-bisection-update-j-end} of Algorithm~\ref{alg:bisection}. We conduct $M=\lceil \log_2 (1/\epsilon) \rceil$ bisection iterations to obtain sufficiently high accuracy.

\section{Main Results}

The first main result of the paper is shown in Theorem~\ref{theo:1}.
\begin{theorem}\label{theo:1}
    Let Assumptions~\ref{assum:1}-\ref{assum:4} hold. Let $\epsilon < \delta$. Under the matching algorithm~\ref{alg:matching} and the pricing algorithm~\ref{alg:pricing} as well as the bisection algorithm~\ref{alg:bisection}, for sufficiently large $T$, we have
    \begin{align}\label{equ:regret-bound}
        \expt[R(T)] = O \Biggl(
        \frac{\log^4(1/\epsilon)}{\epsilon^4 q^{\mathrm{th}}}
        + \frac{T}{q^{\mathrm{th}}}
        + T^2\epsilon^{\frac{\beta}{2}+1} + T\epsilon^{\frac{\beta}{2}-1} \log (1/\epsilon)
        + T\left(\frac{\epsilon}{\delta} + \eta + \delta\right) 
        +  \frac{\log^2 (1/\epsilon)}{\eta \epsilon^2}
        \Biggr),
    \end{align}
    and the queue length
    \begin{align}
        & \frac{1}{T}\sum_{t=1}^{T} \expt \biggl[\sum_i Q_{\mathrm{c}, i}(t) + \sum_j Q_{\mathrm{s}, j}(t)\biggr]\nonumber\\
        =  & O\left(\frac{\log^4 \frac{1}{\epsilon}}{\epsilon^4 T} +  q^{\mathrm{th}}
        +   \max\left\{\frac{\log^2 (1/\epsilon)}{\epsilon^2}, q^{\mathrm{th}}\right\} \left( T \epsilon^{\frac{\beta}{2}+1} + \epsilon^{\frac{\beta}{2}-1} \log (1/\epsilon)\right) \right), \label{equ:avg-queue-length-bound}\\
        & \max_{t=1,\ldots,T} \max\{\max_i Q_{\mathrm{c}, i}(t), \max_j Q_{\mathrm{s}, j}(t)\} 
        \le \max\left\{2\left\lceil \log_2 (1/\epsilon) \right\rceil \left\lceil \frac{\beta \ln (1/\epsilon)}{\epsilon^2} \right\rceil, q^{\mathrm{th}}\right\}.\label{equ:queue-length-bound}
    \end{align}
\end{theorem}
Based on Theorem~\ref{theo:1}, we optimize the regret bound \eqref{equ:regret-bound} by setting the parameters $\epsilon, \eta, \delta, q^{\mathrm{th}}$ to be a function of $T$. Assume that $q^{\mathrm{th}}$ and $\beta$ are large enough. Then we can first focus on minimizing the order of the following term:
\begin{align}\label{equ:regret-optimal-order}
    T\eta + T\delta + \frac{T\epsilon}{\delta} + \frac{\log^2 (1/\epsilon)}{\eta \epsilon^2}.
\end{align}
If we ignore logarithmic order, the optimal solution is $\epsilon=T^{-1/3}$, $\eta=\delta=T^{-1/6}$, and the optimal order of \eqref{equ:regret-optimal-order} is $\tilde{O}(T^{5/6})$. 
By setting $\beta=5$ and $q^{\mathrm{th}}=\Theta(T^{1/2})$, we obtain $\expt[R(T)] = \tilde{O}(T^{5/6})$. We summarize the above result as the following corollary:
\begin{corollary}\label{cor:1}
    Let Assumptions~\ref{assum:1}-\ref{assum:4} hold. Under the matching algorithm~\ref{alg:matching} and the pricing algorithm~\ref{alg:pricing} as well as the bisection algorithm~\ref{alg:bisection}, with parameters $\epsilon=T^{-1/3}$, $\eta=\delta=T^{-1/6}$, $q^{\mathrm{th}}= T^{1/2}$, $\beta=5$, for sufficiently large $T$,
    we can achieve a sublinear regret
    \begin{align*}
        \expt[R(T)] = \tilde{O}(T^{\frac{5}{6}}),
    \end{align*}
    and the queue length
    \begin{align*}
        \frac{1}{T}\sum_{t=1}^{T} \expt \biggl[\sum_i Q_{\mathrm{c}, i}(t) + \sum_j Q_{\mathrm{s}, j}(t)\biggr] = & \tilde{O} ( T^{\frac{1}{2}} ), \nonumber\\
        \max_{t=1,\ldots,T} \max\{\max_i Q_{\mathrm{c}, i}(t), \max_j Q_{\mathrm{s}, j}(t)\} = & \tilde{O}(T^{\frac{2}{3}}).
    \end{align*}
\end{corollary}
From \eqref{equ:regret-bound}, \eqref{equ:avg-queue-length-bound}, and \eqref{equ:queue-length-bound} in Theorem~\ref{theo:1}, we can see that there is a tradeoff between minimizing the regret and minimizing the average or maximum queue length. Firstly, we should not set $q^{\mathrm{th}}$ to be greater than the order of $T^{1/2}$ because it will not benefit the regret bound orderwisely since $q^{\mathrm{th}}=T^{1/2}$ is large enough to achieve the optimal order of \eqref{equ:regret-optimal-order}. However, we can reduce $q^{\mathrm{th}}$ so that the average or maximum queue-length bound can be improved, although we may suffer a larger regret. 

We first look at the tradeoff between regret and maximum queue length. Suppose we want the maximum queue length bound to be the order of $\tilde{O}(T^{\gamma})$, where $\gamma < 2/3$. Then we can set $q^{\mathrm{th}}=T^{\min\{\gamma,1/2\}}$, $\epsilon=T^{-\gamma/2}$, $\eta=\delta=T^{-\gamma/4}$, and $\beta=4/\gamma - 1$. Then by \eqref{equ:regret-bound} and \eqref{equ:queue-length-bound}, we achieve the regret bound $\tilde{O}(T^{1-\gamma/4})$ and maximum queue length bound $\tilde{O}(T^{\gamma})$.
We summarize the above tradeoff as the following corollary:
\begin{corollary}\label{cor:2}
    Let Assumptions~\ref{assum:1}-\ref{assum:4} hold. For any $\gamma \in (0, \frac{2}{3}]$, under the matching algorithm~\ref{alg:matching} and the pricing algorithm~\ref{alg:pricing} as well as the bisection algorithm~\ref{alg:bisection}, with parameters $\epsilon=T^{-\gamma/2}$, $\eta=\delta=T^{-\gamma/4}$, $q^{\mathrm{th}}= T^{\min\{\gamma,1/2\}}$, $\beta=4/\gamma-1$, for sufficiently large $T$,
    we can achieve a sublinear regret
    \begin{align*}
        \expt[R(T)] = \tilde{O}(T^{1-\gamma/4}),
    \end{align*}
    and the maximum queue length
    \begin{align*}
        \max_{t=1,\ldots,T} \max\{\max_i Q_{\mathrm{c}, i}(t), \max_j Q_{\mathrm{s}, j}(t)\} = \tilde{O}(T^{\gamma}).
    \end{align*}
\end{corollary}
From Corollary \ref{cor:1} and Corollary \ref{cor:2}, when the allowable maximum queue length increases (up to $\tilde{\Theta}(T^{2/3})$), the regret bound that we can achieve becomes better. However, even if we allow the maximum queue length to increase over $\tilde{\Theta}(T^{2/3})$, the regret bound cannot be improved and remains $\tilde{O}(T^{5/6})$.

We next look at the tradeoff between regret and average queue length. Suppose we want the average queue length bound to be the order of $\tilde{O}(T^{\gamma})$, where $\gamma < 1/2$.
Then we can set $q^{\mathrm{th}}=T^{\gamma}$, $\epsilon=T^{-2(1+\gamma)/9}$, $\eta=\delta=T^{-(1+\gamma)/9}$, and $\beta=11$. Then by \eqref{equ:regret-bound} and \eqref{equ:avg-queue-length-bound}, we achieve the regret bound $\tilde{O}(T^{\max\{(8-\gamma)/9, 1-\gamma\}})$ and average queue length bound $\tilde{O}(T^{\gamma})$. We summarize the above tradeoff as the following corollary:
\begin{corollary}\label{cor:3}
    Let Assumptions~\ref{assum:1}-\ref{assum:4} hold. For any $\gamma \in (0, \frac{1}{2}]$, under the matching algorithm~\ref{alg:matching} and the pricing algorithm~\ref{alg:pricing} as well as the bisection algorithm~\ref{alg:bisection}, with parameters $\epsilon=T^{-2(1+\gamma)/9}$, $\eta=\delta=T^{-(1+\gamma)/9}$, $q^{\mathrm{th}}= T^{\gamma}$, $\beta=11$, for sufficiently large $T$,
    we can achieve a sublinear regret
    \begin{align*}
        \expt[R(T)] = \tilde{O}(T^{\max\{(8-\gamma)/9, 1-\gamma\}}),
    \end{align*}
    and the average queue length
    \begin{align*}
        \frac{1}{T}\sum_{t=1}^{T} \expt \biggl[\sum_i Q_{\mathrm{c}, i}(t) + \sum_j Q_{\mathrm{s}, j}(t)\biggr] = \tilde{O}(T^{\gamma}).
    \end{align*}
\end{corollary}
From Corollary \ref{cor:1} and Corollary \ref{cor:3}, when the allowable average queue length increases (up to $\tilde{\Theta}(T^{1/2})$), the regret bound that we can achieve becomes better. However, even if we allow the average queue length to increase over $\tilde{\Theta}(T^{1/2})$, the regret bound cannot be improved and remains $\tilde{O}(T^{5/6})$.

Next, we make the following assumption:
\begin{assumption}\label{assum:5}
    Assume that we know some $\boldsymbol{x}(1)\in \mathcal{D}'$ and a set of prices
    \[
    \left\{\underline{p}^{+/-}_{\mathrm{c}, i}(k,1),\bar{p}^{+/-}_{\mathrm{c},i}(k,1), \underline{p}^{+/-}_{\mathrm{s}, j}(k,1), \bar{p}^{+/-}_{\mathrm{s},j}(k,1)\right\}
    _{i\in {\cal I}, j\in {\cal J}}
    \]
    such that
    \begin{align*}
        \sum_{j\in {\cal E}_{\mathrm{c}, i}} x_{i,j}(1) \in & \left[F^{-1}_i(\bar{p}^{+/-}_{\mathrm{c},i}(k,1))-\epsilon + \sqrt{|{\cal E}_{\mathrm{c}, i}|} \delta,
        F^{-1}_i(\underline{p}^{+/-}_{\mathrm{c},i}(k,1)) + \epsilon - \sqrt{|{\cal E}_{\mathrm{c}, i}|} \delta\right], \mbox{for all } i\nonumber\\
        \sum_{i\in {\cal E}_{\mathrm{s}, j}} x_{i,j}(1) \in & \left[G^{-1}_j(\underline{p}^{+/-}_{\mathrm{s},j}(k,1))-\epsilon + \sqrt{|{\cal E}_{\mathrm{s}, j}|} \delta, 
        G^{-1}_j(\bar{p}^{+/-}_{\mathrm{s},j}(k,1)) + \epsilon - \sqrt{|{\cal E}_{\mathrm{s}, j}|} \delta \right], \mbox{for all } j,
    \end{align*}
    and
    \begin{align*}
        \bar{p}^{+/-}_{\mathrm{c},i}(k,1) - \underline{p}^{+/-}_{\mathrm{c},i}(k,1) \le 2e_{\mathrm{c}, i}, \mbox{for all } i\nonumber\\
        \bar{p}^{+/-}_{\mathrm{s},j}(k,1) - \underline{p}^{+/-}_{\mathrm{s},j}(k,1) \le 2 e_{\mathrm{s}, j}, \mbox{for all } j.
    \end{align*}
\end{assumption}
\noindent Assumption~\ref{assum:5} means that we know a set of balanced arrival rates and the bounds for the corresponding prices. Consider a pricing algorithm slightly modified from Algorithm~\ref{alg:pricing} and \ref{alg:bisection}. The differences are as follows. We choose the initial vector $\boldsymbol{x}(1)$ and the initial price searching intervals $[\underline{p}^{+/-}_{\mathrm{c}, i}(k,1),\bar{p}^{+/-}_{\mathrm{c},i}(k,1)], [\underline{p}^{+/-}_{\mathrm{s}, j}(k,1), \bar{p}^{+/-}_{\mathrm{s},j}(k,1)]$ to be the ones in Assumption~\ref{assum:5} and for each queue we reject arrivals if the queue length is larger than or equal to $q^{\mathrm{th}}$ in \textit{all iterations including the first one}. We refer to this modified pricing algorithm as \textit{balanced pricing algorithm}.

Adding Assumption~\ref{assum:5} and using this \textit{balanced pricing algorithm}, we have the following result.
\begin{theorem}\label{theo:2}
    Let Assumptions \ref{assum:1}-\ref{assum:5} hold. Let $\epsilon < \delta$. Under the matching algorithm~\ref{alg:matching} and the balanced pricing algorithm described above, for sufficiently large $T$, we have
    \begin{align}\label{equ:regret-bound-2}
        \expt[R(T)] = O \Biggl(
        \frac{T}{q^{\mathrm{th}}}
        + T^2\epsilon^{\frac{\beta}{2}+1} + T\epsilon^{\frac{\beta}{2}-1} \log (1/\epsilon)
        + T\left(\frac{\epsilon}{\delta} + \eta + \delta\right) 
        +  \frac{\log^2 (1/\epsilon)}{\eta \epsilon^2}
        \Biggr),
    \end{align}
    and the queue length
    \begin{align}
        \frac{1}{T}\sum_{t=1}^{T} \expt \biggl[\sum_i Q_{\mathrm{c}, i}(t) + \sum_j Q_{\mathrm{s}, j}(t)\biggr] \le & q^{\mathrm{th}}, \label{equ:avg-queue-length-bound-2}\\
        \max_{t=1,\ldots,T} \max\{\max_i Q_{\mathrm{c}, i}(t), \max_j Q_{\mathrm{s}, j}(t)\} \le & q^{\mathrm{th}}.\label{equ:queue-length-bound-2}
    \end{align}
\end{theorem}
Compare with Theorem~\ref{theo:1}, the regret bound and the queue length bounds in Theorem~\ref{theo:2} are improved. The reason is that we have a set of nearly balanced arrival rates at the beginning of the algorithm and therefore avoid the rapid increase of queue length in the first iteration due to imbalanced arrival rates and avoid the additional regret induced by rejecting arrivals to reduce the queue length after the first iteration.

Similar to Corollary~\ref{cor:1}, we can optimize the regret bound~\eqref{equ:regret-bound-2} in Theorem~\ref{theo:2} by setting the parameters $\epsilon, \eta, \delta, \beta, q^{\mathrm{th}}$. 
\begin{corollary}\label{cor:4}
    Let Assumptions~\ref{assum:1}-\ref{assum:5} hold. Under the matching algorithm~\ref{alg:matching} and the balanced pricing algorithm, with parameters $\epsilon=T^{-1/3}$, $\eta=\delta=T^{-1/6}$, $q^{\mathrm{th}}= T^{1/6}$, $\beta=5$, for sufficiently large $T$,
    we can achieve a sublinear regret
    \begin{align*}
        \expt[R(T)] = \tilde{O}(T^{\frac{5}{6}}),
    \end{align*}
    and the queue length
    \begin{align*}
        \frac{1}{T}\sum_{t=1}^{T} \expt \biggl[\sum_i Q_{\mathrm{c}, i}(t) + \sum_j Q_{\mathrm{s}, j}(t)\biggr] = & \tilde{O} ( T^{\frac{1}{6}} ), \nonumber\\
        \max_{t=1,\ldots,T} \max\{\max_i Q_{\mathrm{c}, i}(t), \max_j Q_{\mathrm{s}, j}(t)\} = & \tilde{O}(T^{\frac{1}{6}}).
    \end{align*}
\end{corollary}
From Corollary~\ref{cor:4}, we can see that both the average and the maximum queue length bounds are significantly reduced while the optimal regret order is the same as that in Corollary~\ref{cor:1}.

The following Corollary~\ref{cor:5} illustrates the tradeoff among regret, average queue length, and maximum queue length.
\begin{corollary}\label{cor:5}
    Let Assumptions~\ref{assum:1}-\ref{assum:5} hold. For any $\gamma \in (0, \frac{1}{6}]$, under the matching algorithm~\ref{alg:matching} and the balanced pricing algorithm, with parameters $\epsilon=T^{-2\gamma}$, $\eta=\delta=T^{-\gamma}$, $q^{\mathrm{th}}= T^{\gamma}$, $\beta=1/\gamma-1$, for sufficiently large $T$, we can achieve a sublinear regret
    \begin{align*}
        \expt[R(T)] = \tilde{O}(T^{1-\gamma}),
    \end{align*}
    and the queue length
    \begin{align*}
        \frac{1}{T}\sum_{t=1}^{T} \expt \biggl[\sum_i Q_{\mathrm{c}, i}(t) + \sum_j Q_{\mathrm{s}, j}(t)\biggr] = & \tilde{O} ( T^{\gamma} ), \nonumber\\
        \max_{t=1,\ldots,T} \max\{\max_i Q_{\mathrm{c}, i}(t), \max_j Q_{\mathrm{s}, j}(t)\} = & \tilde{O}(T^{\gamma}).
    \end{align*}
\end{corollary}
From Corollary \ref{cor:4} and Corollary \ref{cor:5}, when the allowable maximum and average queue length increases (up to $\tilde{\Theta}(T^{1/6})$), the regret bound that we can achieve becomes better. However, even if we allow the maximum and average queue length to increase over $\tilde{\Theta}(T^{1/6})$, the regret bound cannot be improved and remains $\tilde{O}(T^{5/6})$.

\section{Proof Roadmap}

In this section, we present proof roadmaps for Theorem~\ref{theo:1} and Theorem~\ref{theo:2}.

We first present the proof roadmap for Theorem~\ref{theo:1}. The complete proof of Theorem~\ref{theo:1} can be found in Appendix~\ref{app:theo:1}.
In the proof roadmap, we will consider the choice of parameters as in Corollary~\ref{cor:1} to make the proof idea clear, i.e.,
$\epsilon=T^{-1/3}$, $\eta=\delta=T^{-1/6}$, $q^{\mathrm{th}}= T^{1/2}$, $\beta=5$.

For the first outer iteration, we have no control of the queue length, so the queue length will possibly increase up to $2MN=\tilde{\Theta}(T^{2/3})$. After the first outer iteration, the queue length is controlled by the threshold $q^{\mathrm{th}}= T^{1/2}$. As a result, the anytime maximum queue-length bound is $\tilde{O}(T^{2/3})$ and the average queue-length bound is $\tilde{O}(T^{1/2})$.

\begin{figure}[htbp]
    \centering
    \includegraphics[width=1.0\linewidth]{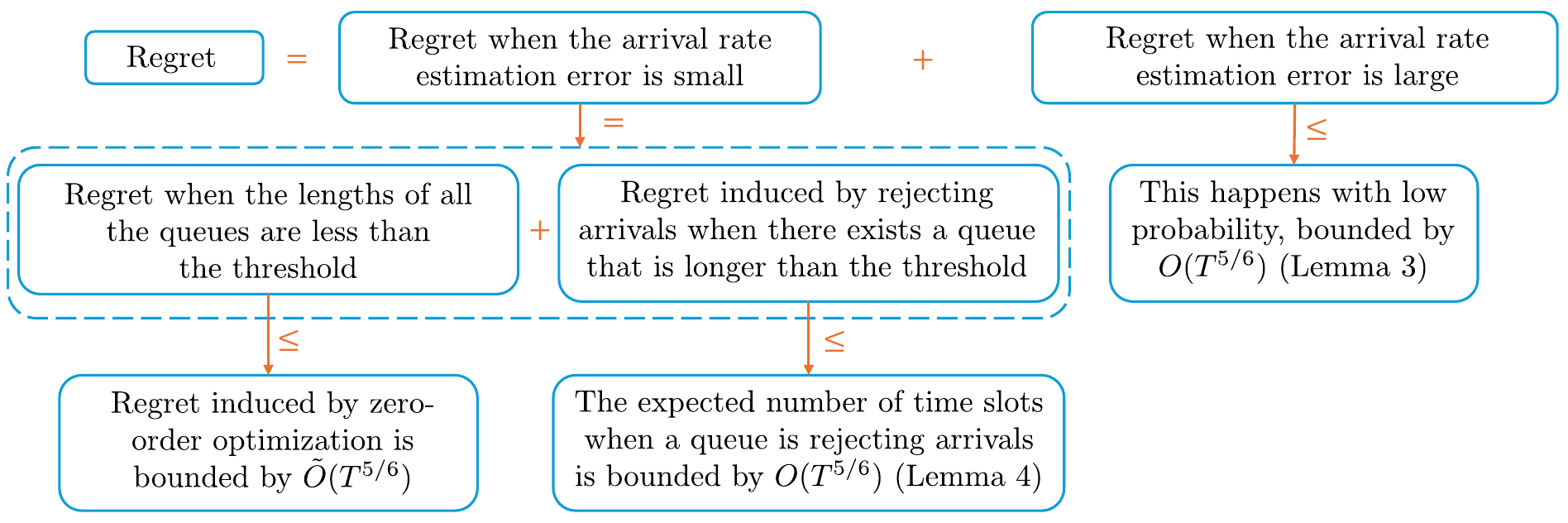}
    \caption{\small Proof Roadmap.}
    \label{fig:roadmap}
\end{figure}

For the regret bound, the proof roadmap is shown in Figure~\ref{fig:roadmap}. We first define an event ${\cal C}$, which means that the estimations of all the arrival rates are $\epsilon$-accurate. Then the regret can be divided into the regret when the event ${\cal C}$ occurs and the regret when the event ${\cal C}$ does not occur. We show in Lemma~\ref{lemma:concentration} in Appendix~\ref{sec:roadmap-add-concentration} that the complement of the event, ${\cal C}^{\mathrm{c}}$, happens with low probability. Intuitively, since for each estimation we use $N=\tilde{O}(1/\epsilon^2)$ samples, we can obtain $\epsilon$ accuracy with high probability by Hoeffding's inequality. Note that Lemma~\ref{lemma:concentration} cannot be proved by simply applying Hoeffding's inequality and we will discuss the challenge in Section~\ref{sec:challenge}.

Therefore, the regret when the event ${\cal C}$ does not happen can be bounded. For the regret when the event ${\cal C}$ happens, i.e., the estimations are $\epsilon$-accurate, we further divide it into two parts. One part is the regret when the lengths of all the queues are less than the threshold $q^{\mathrm{th}}$ so that there is no queue that is rejecting arrivals, and the other part is the regret induced by rejecting arrivals when there exists a queue whose length is greater than or equal to the threshold $q^{\mathrm{th}}$. For the latter part, we can bound it by the expected number of time slots when a queue is rejecting arrivals. We show in Lemma~\ref{lemma:bound-time-large-queue} in Appendix~\ref{app:bound-regret-reject-arr} that this can be bounded by $O(T^{5/6})$. The proof is using Lyapunov drift analysis with a Lyapunov function $V_{\mathrm{c}}(t)\coloneqq \sum_i Q_{\mathrm{c},i}^2(t)$ for the customer side and $V_{\mathrm{s}}(t)\coloneqq \sum_j Q_{\mathrm{s},j}^2(t)$ for the server side. The challenge of the proof will be discussed in Section~\ref{sec:challenge}.

For the part of regret when there is no queue that is rejecting arrivals, regret is induced by zero-order optimization process. We show in Appendix~\ref{app:regret-zero-order-opt} that this part of regret can be bounded by $\tilde{O}(T^{5/6})$. The intuition is as follows. If there is only one run of the system in each outer iteration, then from the literature of two-point method of zero-order optimization~\cite{agarwal2010optimal}, the regret will be bounded by $O(\sqrt{K})$, where $K$ is the number of outer iterations. In our algorithm, there are $\tilde{\Theta}(T^{2/3})$ runs of the systems in each outer iterations and there are $\tilde{\Theta}(T^{1/3})$ outer iterations. Hence, the regret will be bounded by $\tilde{\Theta}(T^{2/3}) O(\sqrt{T^{1/3}})=\tilde{O}(T^{5/6})$. However, the proof needs additional steps because even in the same outer iteration, the prices and arrival rates are changing. In fact, we need to connect the actual prices to the prices corresponding to $\boldsymbol{x}(k)$.
Moreover, there can be dependence between the arrival rate estimation and the randomness of the optimization algorithm, which we will discuss in Section~\ref{sec:challenge}.

For Theorem~\ref{theo:2}, the complete proof can be found in Appendix~\ref{app:theo:2}.
The proof roadmap of Theorem~\ref{theo:2} is similar to that of Theorem~\ref{theo:1}. The outline is the same as that in Figure~\ref{fig:roadmap} but the difference is that the regret induced by rejecting arrivals will be significantly reduced.
Under Assumption~\ref{assum:5}, we know that there is a set of nearly balanced arrival rates at the beginning of the algorithm. Then under the \textit{balanced pricing algorithm}, we show in Lemma~\ref{lemma:improved-bisection-error} that the arrival rates will remain nearly balanced for all time. This avoids the rapid increase of queue length in the first iteration due to imbalanced arrival rates and therefore avoids the additional regret induced by rejecting arrivals to reduce the
queue length after the first iteration. Hence, the regret induced by rejecting arrivals when there exists a queue whose length is greater than or equal to the threshold is significantly reduced. Therefore, compared to Theorem~\ref{theo:1}, Theorem~\ref{theo:2} has a better tradeoff between regret and queue length, as shown in Corollary~\ref{cor:5}.

\section{Challenges in the Proof}
\label{sec:challenge}

In this section, we discuss the challenges in the proof and the methods we use to overcome them.

For the proof of Theorem~\ref{theo:1}, The first challenge is to prove that the event ${\cal C}^{\mathrm{c}}$ occurs with low probability (Lemma~\ref{lemma:concentration} in Appendix~\ref{sec:roadmap-add-concentration}). The event ${\cal C}$ is defined as
\begin{align*}
    & {\cal C} \coloneqq \nonumber\\
    \Biggl\{
    & \text{for all outer iteration $k=1,\ldots,\lceil T/(2MN) \rceil$, all bisection iteration $m=1,\ldots, M$, all $i,j$,}\nonumber\\ 
    & \text{and all arrival rates $\lambda_i(t^{+}_{\mathrm{c}, i}(k,m,n)), \lambda_i(t^{-}_{\mathrm{c}, i}(k,m,n)), \mu_j(t^{+}_{\mathrm{s}, j}(k,m,n)), \mu_j(t^{-}_{\mathrm{s}, j}(k,m,n)) \in[0,1]$}, \nonumber\\
    & \biggl| \frac{1}{N} \sum_{n=1}^{N} A_{\mathrm{c},i} (t^{+}_{\mathrm{c}, i}(k,m,n)) - \lambda_i(t^{+}_{\mathrm{c}, i}(k,m,n)) \biggr| < \epsilon, \biggl| \frac{1}{N} \sum_{n=1}^{N} A_{\mathrm{c},i} (t^{-}_{\mathrm{c}, i}(k,m,n)) - \lambda_i(t^{-}_{\mathrm{c}, i}(k,m,n)) \biggr| < \epsilon, \nonumber\\
    & \biggl| \frac{1}{N} \sum_{n=1}^{N} A_{\mathrm{s},j} (t^{+}_{\mathrm{s}, j}(k,m,n)) - \mu_j(t^{+}_{\mathrm{s}, j}(k,m,n)) \biggr| < \epsilon, \biggl| \frac{1}{N} \sum_{n=1}^{N} A_{\mathrm{s},j} (t^{-}_{\mathrm{s}, j}(k,m,n)) - \mu_j(t^{-}_{\mathrm{s}, j}(k,m,n)) \biggr| < \epsilon
    \Biggr\},
\end{align*}
where $t^{+}_{\mathrm{c}, i}(k,m,n)$ is the time slot when the price $p^{+}_{\mathrm{c}, i} (k,m)$ is run for the $n^{\mathrm{th}}$ time for the customer-side queue $i$, and similarly we define $t^{-}_{\mathrm{c}, i}(k,m,n)$, $t^{+}_{\mathrm{s}, j}(k,m,n)$, and $t^{-}_{\mathrm{s}, j}(k,m,n)$. 
In the definition, we require that the accurate estimation holds for any arrival rates within $[0,1]$ because we want to ensure the independence between the event ${\cal C}$ and the randomness of $\boldsymbol{u}(k)$ in Algorithm~\ref{alg:pricing}, which we will discuss later in this section. To prove that the event ${\cal C}$ holds with high probability, simply using Hoeffding's inequality with the union bound does not work because arrival rates are continuous random variables within $[0,1]$. Therefore, we use a discretization idea that is similar to the covering argument for linear bandits~\cite{lattimore2020bandit}. We define a set ${\cal S}_{[0,1]}$, a discretized version of $[0,1]$ with resolution $\epsilon/2$ and rewrite the first inequality in the definition of ${\cal C}$ as
$|\frac{1}{N} \sum_{n=1}^{N} \mathbb{1} \{U^+_{\mathrm{c},i}(k,m,n) \le \lambda \} - \lambda | < \epsilon$ where $U^+_{\mathrm{c},i}(k,m,n)$ is a random variable uniformly distributed in $[0,1]$ that is generated independently at the very beginning of the process, and similarly we can rewrite other inequalities in the definition of ${\cal C}$.
We can bound the quantization error and then use Hoeffding's inequality and the union bound over $\lambda\in {\cal S}_{[0,1]}$ to bound the probability of ${\cal C}^{\mathrm{c}}$.

The second challenge is to prove the upper bound on the expected number of time slots when the length of one of the queues is greater than or equal to the threshold $q^{\mathrm{th}}$ (Lemma~\ref{lemma:bound-time-large-queue} in Appendix~\ref{app:bound-regret-reject-arr}). As mentioned in the proof roadmap in the previous section, we use Lyapunov drift analysis with quadratic Lyapunov functions. The challenge is how we make sure that the actual arrival rates are balanced or almost balanced. We know that $\boldsymbol{x}(k) \in {\cal D}'$ so the corresponding arrival rates are balanced. Since $\boldsymbol{x}^{+}(k)$ and $\boldsymbol{x}^{-}(k)$ are $\delta$-close to $\boldsymbol{x}(k)$, their corresponding arrival rates are almost balanced.
We want to prove that the actual arrival rates are close enough to the almost balanced arrival rates corresponding to $\boldsymbol{x}^{+}(k)$ or $\boldsymbol{x}^{-}(k)$. To show this, we use an important property that the true price corresponding to $\boldsymbol{x}^{+}(k)$ or $\boldsymbol{x}^{-}(k)$ is within the search interval of the bisection algorithm with high probability and an important fact that the width of the search interval is small enough for $k\ge 2$. 
Therefore, the actual arrival rates are also almost balanced for $k\ge 2$.

The third challenge is bounding the regret induced by zero-order optimization conditioned on $\epsilon$-accurate arrival rate estimations (Lemma~\ref{lemma:optimization} in Appendix~\ref{app:regret-zero-order-opt}). As mentioned in the proof roadmap in the previous section, the challenge is that there could be dependence between the estimation of arrival rates and the randomness ($\boldsymbol{u}(k)$) in the optimization algorithm. By defining the event ${\cal C}$ such that the $\epsilon$-accuracy is required for all arrival rates in $[0,1]$, we can show that the event ${\cal C}$ is independent of $\boldsymbol{u}(k)$. Note that if we change the definition of ${\cal C}$ such that only the estimations of the actual arrival rates are $\epsilon$-accurate, this independence does not hold.
This independence between ${\cal C}$ and $\boldsymbol{u}(k)$ is essential so that conditioned on the event ${\cal C}$, we can use the proof idea of two-point method of zero-order optimization to bound the regret.

Although the proof of Theorem~\ref{theo:2} is similar to that of Theorem~\ref{theo:1}, the key point in the proof of Theorem~\ref{theo:2} is to show that the arrival rates will remain balanced all the time (Lemma~\ref{lemma:improved-bisection-error}). The key idea is to use an induction argument and verify the base case with Assumption~\ref{assum:5}.

\section{Conclusions}

We studied the pricing and matching problem for two-sided queueing systems with unknown demand and supply functions. We proposed a novel learning-based pricing algorithm that combines zero-order stochastic projected gradient ascent with bisection search. 
We proved that the proposed algorithm yields a sublinear regret $\tilde{O}(T^{5/6})$, while ensuring an anytime queue-length bound of $\tilde{O}(T^{2/3})$ and an average queue-length bound of $\tilde{O}(T^{1/2})$. We also established tradeoffs between the regret bound and the anytime queue-length bound, as well as between the regret bound and the average queue-length bound.
Furthermore, assuming the availability of a set of balanced arrival rates and associated price intervals, we established an improved tradeoff between the regret bound and the queue-length bound: $\tilde{O}(T^{1-\gamma})$ versus $\tilde{O}(T^{\gamma})$ for $\gamma\in(0, 1/6]$.


\clearpage
\newpage
\bibliographystyle{unsrt}
\bibliography{refs}


\clearpage
\newpage
\appendix

\clearpage

\section{Proof of Proposition~\ref{prop:1}}
\label{app:proof-prop-1}

We first show that any mean rate stable policy must satisfy the constraints \eqref{equ:fluid-opti-constr1} and \eqref{equ:fluid-opti-constr2}.

Recall 
\begin{align*}
    x_{i,j}\coloneqq & \lim_{T\rightarrow \infty} \frac{1}{T} \sum_{t=1}^{T} \expt[X_{i,j}(t)]\\
    \lambda_i\coloneqq & \lim_{T\rightarrow \infty} \frac{1}{T} \sum_{t=1}^{T} \expt[ \lambda_i(t)]\\
    \mu_j\coloneqq & \lim_{T\rightarrow \infty} \frac{1}{T} \sum_{t=1}^{T} \expt[\mu_j(t)].
\end{align*}
Note that the dynamics of the queue length of type $i$ customers is
\begin{align}\label{equ:dynamics}
    Q_{\mathrm{c}, i}(t+1) = Q_{\mathrm{c}, i}(t) + A_{\mathrm{c}, i}(t) - \sum_{j: (i,j)\in {\cal E}} X_{i,j}(t).
\end{align}
Taking expectation on both sides, we have
\begin{align*}
    \expt [Q_{\mathrm{c}, i}(t+1)] =  \expt [Q_{\mathrm{c}, i}(t) ] + \expt[\lambda_{i}(t)] - \sum_{j: (i,j)\in {\cal E}} \expt[X_{i,j}(t)].
\end{align*}
Taking time average on both sides, we have
\begin{align*}
     \frac{1}{T} \sum_{t=1}^{T} \expt [Q_{\mathrm{c}, i}(t+1)] = \frac{1}{T} \sum_{t=1}^{T} \expt [Q_{\mathrm{c}, i}(t) ] +  \frac{1}{T} \sum_{t=1}^{T} \expt[\lambda_{i}(t)] - \sum_{j: (i,j)\in {\cal E}} \frac{1}{T} \sum_{t=1}^{T} \expt[X_{i,j}(t)],
\end{align*}
which implies that 
\begin{align*}
    \frac{1}{T} \sum_{t=1}^{T} \expt[\lambda_{i}(t)] - \sum_{j: (i,j)\in {\cal E}} \frac{1}{T} \sum_{t=1}^{T} \expt[X_{i,j}(t)] = \frac{1}{T} \sum_{t=1}^{T} \left(\expt [Q_{\mathrm{c}, i}(t+1)] -  \expt [Q_{\mathrm{c}, i}(t) ]\right) = \frac{1}{T} \expt [Q_{\mathrm{c}, i}(T+1)].
\end{align*}
Taking limit on both sides and by the definition of mean rate stability~\eqref{equ:mean-rate-stable}, we have $\lambda_i - \sum_{j: (i,j)\in {\cal E}} x_{i,j} = 0$.
Similar argument holds for the server side and hence, $\mu_j = \sum_{i: (i,j)\in {\cal E}} x_{i,j}$.

Now we show that the optimal value of the optimization problem~\eqref{equ:fluid-opti-obj}-\eqref{equ:fluid-opti-constr4} is an upper bound of the asymptotic time-averaged profit under any mean rate stable policy, i.e., we will show
\begin{align*}
    \lim_{t\rightarrow \infty}\frac{1}{T} \sum_{t=1}^{T} \expt \left[ \sum_i A_{\mathrm{c},i}(t) F_i(\lambda_i(t)) - \sum_j A_{\mathrm{s},j}(t) G_j(\mu_j(t)) \right]
    \le \sum_i \lambda_i^* F_i(\lambda_i^*) - \sum_j \mu_j^* G_j(\mu_j^*),
\end{align*}
where $(\boldsymbol{\lambda}^*, \boldsymbol{\mu}^*)$ denotes an optimal solution to the optimization problem~\eqref{equ:fluid-opti-obj}-\eqref{equ:fluid-opti-constr4}. First note that
\begin{align}\label{equ:expected-profit}
    & \expt \left[ \sum_i A_{\mathrm{c},i}(t) F_i(\lambda_i(t)) - \sum_j A_{\mathrm{s},j}(t) G_j(\mu_j(t)) \right]\nonumber\\
    = & \expt \biggl[ \expt \Bigl[ \sum_i A_{\mathrm{c},i}(t) F_i(\lambda_i(t)) - \sum_j A_{\mathrm{s},j}(t) G_j(\mu_j(t)) \Bigl| \boldsymbol{\lambda}(t), \boldsymbol{\mu}(t) \Bigr. \Bigr] \biggr]\nonumber\\
    = & \expt \biggl[ \sum_i \lambda_{i}(t) F_i(\lambda_i(t)) - \sum_j \mu_{j}(t) G_j(\mu_j(t))\biggr].
\end{align}
By Assumption~\ref{assum:2}, we have
\begin{align}\label{equ:upper-bound-convexity}
    & \frac{1}{T} \sum_{t=1}^{T} \expt \biggl[ \sum_i \lambda_{i}(t) F_i(\lambda_i(t)) - \sum_j \mu_{j}(t) G_j(\mu_j(t))\biggr]\nonumber\\
    \le & \expt \biggl[ \sum_i \biggl(\frac{1}{T} \sum_{t=1}^{T} \lambda_{i}(t) \biggr)
    F_i\biggl(\frac{1}{T} \sum_{t=1}^{T} \lambda_i(t)\biggr)
    - \sum_j \biggl(\frac{1}{T} \sum_{t=1}^{T} \mu_{j}(t) \biggr)
    G_j\biggl( \frac{1}{T} \sum_{t=1}^{T} \mu_j(t)\biggr)\biggr].
\end{align}
By Assumption~\ref{assum:2} and the Jensen's inequality, we have
\begin{align}\label{equ:upper-bound-jensen}
    & \expt \biggl[ \sum_i \biggl(\frac{1}{T} \sum_{t=1}^{T} \lambda_{i}(t) \biggr)
    F_i\biggl(\frac{1}{T} \sum_{t=1}^{T} \lambda_i(t)\biggr)
    - \sum_j \biggl(\frac{1}{T} \sum_{t=1}^{T} \mu_{j}(t) \biggr)
    G_j\biggl( \frac{1}{T} \sum_{t=1}^{T} \mu_j(t)\biggr)\biggr]\nonumber\\
    \le &  \sum_i \biggl( \frac{1}{T} \sum_{t=1}^{T} \expt[\lambda_{i}(t)] \biggr)
    F_i\biggl(\frac{1}{T} \sum_{t=1}^{T} \expt[\lambda_i(t)]\biggr)
    - \sum_j \biggl(\frac{1}{T} \sum_{t=1}^{T} \expt[\mu_{j}(t)] \biggr)
    G_j\biggl( \frac{1}{T} \sum_{t=1}^{T} \expt[\mu_j(t)]\biggr).
\end{align}
Taking limit over $T$, by Assumption~\ref{assum:1}, we have
\begin{align}\label{equ:dct-continuity}
    & \lim_{T\rightarrow \infty} \sum_i \biggl( \frac{1}{T} \sum_{t=1}^{T} \expt[\lambda_{i}(t)] \biggr)
    F_i\biggl(\frac{1}{T} \sum_{t=1}^{T} \expt[\lambda_i(t)]\biggr)
    - \lim_{T\rightarrow \infty} \sum_j \biggl(\frac{1}{T} \sum_{t=1}^{T} \expt[\mu_{j}(t)] \biggr)
    G_j\biggl( \frac{1}{T} \sum_{t=1}^{T} \expt[\mu_j(t)]\biggr)\nonumber\\
    = & \sum_i \biggl( \lim_{T\rightarrow \infty} \frac{1}{T} \sum_{t=1}^{T} \expt[\lambda_{i}(t)] \biggr)
    F_i\biggl(\lim_{T\rightarrow \infty} \frac{1}{T} \sum_{t=1}^{T} \expt[\lambda_i(t)]\biggr)\nonumber\\
    & - \sum_j \biggl( \lim_{T\rightarrow \infty}\frac{1}{T} \sum_{t=1}^{T} \expt[\mu_{j}(t)] \biggr)
    G_j\biggl( \lim_{T\rightarrow \infty} \frac{1}{T} \sum_{t=1}^{T} \expt[\mu_j(t)]\biggr).
\end{align}
At the beginning of this section, we show that for any mean rate stable policy,
\[
\lim_{T\rightarrow \infty} \frac{1}{T} \sum_{t=1}^{T} \expt[\lambda_i(t)] =  \sum_{j:(i,j)\in{\cal E}} \lim_{T\rightarrow \infty} \frac{1}{T} \sum_{t=1}^{T} \expt[X_{i,j}(t)]
\]
and
\[
\lim_{T\rightarrow \infty} \frac{1}{T} \sum_{t=1}^{T} \expt[\mu_j(t)] =  \sum_{i:(i,j)\in{\cal E}} \lim_{T\rightarrow \infty} \frac{1}{T} \sum_{t=1}^{T} \expt[X_{i,j}(t)]
\]
for all $i,j$. 
Note that
\[
\lim_{T\rightarrow \infty} \frac{1}{T} \sum_{t=1}^{T} \expt[X_{i,j}(t)] \ge 0
\]
for all $i,j$, and
\[
\lim_{T\rightarrow \infty} \frac{1}{T} \sum_{t=1}^{T} \expt[\lambda_i(t)] \in [0,1], \lim_{T\rightarrow \infty} \frac{1}{T} \sum_{t=1}^{T} \expt[\mu_j(t)] \in [0,1].
\]
for all $i,j$. Hence,
\[
\left(\lim_{T\rightarrow \infty} \frac{1}{T} \sum_{t=1}^{T} \expt[\boldsymbol{\lambda}(t)], \lim_{T\rightarrow \infty} \frac{1}{T} \sum_{t=1}^{T} \expt[\boldsymbol{\mu}(t)],
\lim_{T\rightarrow \infty} \frac{1}{T} \sum_{t=1}^{T} \expt[\boldsymbol{X}(t)]\right)
\]
is a feasible solution to the optimization problem~\eqref{equ:fluid-opti-obj}-\eqref{equ:fluid-opti-constr4}, where $\boldsymbol{X}(t)\coloneqq (X_{i,j}(t))_{(i,j)\in{\cal E}}$. Hence,
\begin{align}\label{equ:fluid-upper-bound-final}
    & \sum_i \biggl( \lim_{T\rightarrow \infty} \frac{1}{T} \sum_{t=1}^{T} \expt[\lambda_{i}(t)] \biggr)
    F_i\biggl(\lim_{T\rightarrow \infty} \frac{1}{T} \sum_{t=1}^{T} \expt[\lambda_i(t)]\biggr)\nonumber\\
    & - \sum_j \biggl( \lim_{T\rightarrow \infty}\frac{1}{T} \sum_{t=1}^{T} \expt[\mu_{j}(t)] \biggr)
    G_j\biggl( \lim_{T\rightarrow \infty} \frac{1}{T} \sum_{t=1}^{T} \expt[\mu_j(t)]\biggr)\nonumber\\
    \le & \sum_i \lambda_i^* F_i(\lambda_i^*) - \sum_j \mu_j^* G_j(\mu_j^*).
\end{align}
From \eqref{equ:expected-profit}, \eqref{equ:upper-bound-convexity}, \eqref{equ:upper-bound-jensen}, \eqref{equ:dct-continuity}, and \eqref{equ:fluid-upper-bound-final}, we have
\begin{align*}
    \lim_{t\rightarrow \infty}\frac{1}{T} \sum_{t=1}^{T} \expt \left[ \sum_i A_{\mathrm{c},i}(t) F_i(\lambda_i(t)) - \sum_j A_{\mathrm{s},j}(t) G_j(\mu_j(t)) \right]
    \le \sum_i \lambda_i^* F_i(\lambda_i^*) - \sum_j \mu_j^* G_j(\mu_j^*).
\end{align*}
Proposition~\ref{prop:1} is proved.

\section{Exact Expressions of \texorpdfstring{$e_{\mathrm{c},i}$ and $e_{\mathrm{s},j}$}{e\_\{c,i\} and e\_\{s,j\}} in Algorithm~\ref{alg:pricing}}
\label{app:exact-expression}

\begin{align}
    e_{\mathrm{c},i}
    =&\frac{2 \eta \epsilon |{\cal E}|^{3/2} L_{F_i} }{\delta} 
     \Biggl[ \sum_{i'=1}^{I} L_{F_{i'}}
     \left(1 + L_{F^{-1}_{i'}} \left(p_{\mathrm{c},i',\max} - p_{\mathrm{c},i',\min}\right)\right)\nonumber\\
     & + \sum_{j'=1}^{J} L_{G_{j'}}
    \left( 1 + L_{G^{-1}_{j'}}
    \left(p_{\mathrm{s},j',\max} - p_{\mathrm{s},j',\min}\right)\right)\Biggr]
    + 2\epsilon L_{F_i} \left(1 + L_{F^{-1}_i} \left(p_{\mathrm{c},i,\max} - p_{\mathrm{c},i,\min}\right)\right)\nonumber\\
    & + \eta |{\cal E}|^{3/2} L_{F_i} \left( \sum_{i'=1}^{I} |{\cal E}_{\mathrm{c},i'}| ( L_{F_{i'}} + p_{\mathrm{c}, i', \max} ) + \sum_{j'=1}^{J}  |{\cal E}_{\mathrm{s},j'}| ( L_{G_{j'}} + p_{\mathrm{s}, j', \max} ) \right) + 2 \delta |{\cal E}|^{1/2} L_{F_i}\label{equ:def-e-c}\\
    = & \Theta\left(\frac{\eta \epsilon}{\delta} + \eta + \delta + \epsilon\right),\nonumber\\
    e_{\mathrm{s},j}
    =&\frac{2 \eta \epsilon |{\cal E}|^{3/2} L_{G_j} }{\delta} 
     \Biggl[ \sum_{i'=1}^{I} L_{F_{i'}}
     \left(1 + L_{F^{-1}_{i'}} \left(p_{\mathrm{c},i',\max} - p_{\mathrm{c},i',\min}\right)\right)\nonumber\\
    & + \sum_{j'=1}^{J} L_{G_{j'}}
    \left( 1 + L_{G^{-1}_{j'}}
    \left(p_{\mathrm{s},j',\max} - p_{\mathrm{s},j',\min}\right)\right)\Biggr]
    + 2\epsilon L_{G_j} \left(1 + L_{G^{-1}_j} \left( p_{\mathrm{s},j,\max} - p_{\mathrm{s},j,\min}\right)\right)\nonumber\\
    & + \eta |{\cal E}|^{3/2} L_{G_j} \left( \sum_{i'=1}^{I} |{\cal E}_{\mathrm{c},i'}| ( L_{F_{i'}} + p_{\mathrm{c}, i', \max} ) + \sum_{j'=1}^{J}  |{\cal E}_{\mathrm{s},j'}| ( L_{G_{j'}} + p_{\mathrm{s}, j', \max} ) \right) + 2 \delta |{\cal E}|^{1/2} L_{G_j}\label{equ:def-e-s}\\
    = & \Theta\left(\frac{\eta \epsilon}{\delta} + \eta + \delta + \epsilon\right)\nonumber.
\end{align}

\section{Proof of Theorem~\ref{theo:1}}
\label{app:theo:1}

Before proving Theorem~\ref{theo:1}, we first present some preliminary lemmas.

\subsection{Preliminary Lemmas}
\label{app:theo:1-pre-lemma}

Under the matching algorithm \ref{alg:matching}, we have the following lemma.
\begin{lemma}\label{lemma:empty-queues}
    Assume all queues are empty in the beginning. Under Algorithm~\ref{alg:matching} and any pricing algorithm, at any time slot $t$, if a queue $i$ (or $j$) is not empty before the arrivals, then all the queues on the other side that are linked to this queue are all empty, i.e., $\sum_{j: (i,j)\in{\cal E}} Q_{\mathrm{s}, j}(t) = 0$ (or $\sum_{i: (i,j)\in{\cal E}} Q_{\mathrm{c}, i}(t) = 0$).
\end{lemma}
From Lemma~\ref{lemma:empty-queues}, we know that the number of matches on a compatible link $(i,j)$ at time $t$ , $X_{i,j}(t)$, is either zero or one because it is impossible that $Q_{\mathrm{c}, i}(t) + A_{\mathrm{c}, i}(t) \ge 2$ and $Q_{\mathrm{s}, j}(t) + A_{\mathrm{s}, j}(t) \ge 2$ hold simultaneously for Bernoulli arrivals. Proof of Lemma~\ref{lemma:empty-queues} can be found in Appendix~\ref{app:proof-lemma-empty-queues}.

Next, we present a lemma characterizing the relation between the feasible set ${\cal D}$ and the shrunk feasible set ${\cal D}'$.
\begin{lemma}\label{lemma:shrunk-set}
    Under Assumption~\ref{assum:4}, we have ${\cal D}' = (1-\frac{\delta}{r})({\cal D} - \boldsymbol{x}_{\mathrm{ctr}})+ \boldsymbol{x}_{\mathrm{ctr}}$, where $\boldsymbol{x}_{\mathrm{ctr}} = \left(\frac{a_{\min}+1}{2N_{i,j}}\right)_{(i,j)\in {\cal E}}$. And we have $\boldsymbol{x} + \delta \boldsymbol{u} \in {\cal D}$ for any $\boldsymbol{x}\in {\cal D}'$ and any vector $\boldsymbol{u}\in {\mathbb B}$, where ${\mathbb B}$ is the unit ball in $\mathbb{R}^{|\cal E|}$.
\end{lemma}
By Lemma~\ref{lemma:shrunk-set}, for any iteration $k$, we have $\boldsymbol{x}^{+}(k)\in {\cal D}$ and $\boldsymbol{x}^{-}(k)\in {\cal D}$, which implies that the arrival rates $\lambda_i^{+}(k)$, $\lambda_i^{-}(k)$, $\mu_j^{+}(k)$, and $\mu_j^{-}(k)$ in Algorithm~\ref{alg:pricing} are in $[a_{\min}, 1]$ for all $i,j$, i.e., they are valid arrival rates.
Proof of Lemma~\ref{lemma:shrunk-set} can be found in Appendix~\ref{app:proof-lemma-shrunk-set}. 

In the following subsections, we present the complete proof of Theorem \ref{theo:1}.

\subsection{Adding the Event of Concentration}
\label{sec:roadmap-add-concentration}
Before bounding the regret, we first define a ``good'' event ${\cal C}$.
With a fixed total number of time slots $T$, the number of outer iterations is a random variable and we know that it is less than or equal to $T/(2MN)$ because there are at least $N$ time slots in each bisection iteration and there are $2M$ bisection iterations in each outer iteration.
Define the event ${\cal C}$ as
\begin{align*}
    & {\cal C} \coloneqq \nonumber\\
    \Biggl\{
    & \text{for all outer iteration $k=1,\ldots,\lceil T/(2MN) \rceil$, all bisection iteration $m=1,\ldots, M$, all $i,j$,}\nonumber\\ 
    & \text{and all arrival rates $\lambda_i(t^{+}_{\mathrm{c}, i}(k,m,n)), \lambda_i(t^{-}_{\mathrm{c}, i}(k,m,n)), \mu_j(t^{+}_{\mathrm{s}, j}(k,m,n)), \mu_j(t^{-}_{\mathrm{s}, j}(k,m,n)) \in[0,1]$}, \nonumber\\
    & \biggl| \frac{1}{N} \sum_{n=1}^{N} A_{\mathrm{c},i} (t^{+}_{\mathrm{c}, i}(k,m,n)) - \lambda_i(t^{+}_{\mathrm{c}, i}(k,m,n)) \biggr| < \epsilon, \biggl| \frac{1}{N} \sum_{n=1}^{N} A_{\mathrm{c},i} (t^{-}_{\mathrm{c}, i}(k,m,n)) - \lambda_i(t^{-}_{\mathrm{c}, i}(k,m,n)) \biggr| < \epsilon, \nonumber\\
    & \biggl| \frac{1}{N} \sum_{n=1}^{N} A_{\mathrm{s},j} (t^{+}_{\mathrm{s}, j}(k,m,n)) - \mu_j(t^{+}_{\mathrm{s}, j}(k,m,n)) \biggr| < \epsilon, \biggl| \frac{1}{N} \sum_{n=1}^{N} A_{\mathrm{s},j} (t^{-}_{\mathrm{s}, j}(k,m,n)) - \mu_j(t^{-}_{\mathrm{s}, j}(k,m,n)) \biggr| < \epsilon
    \Biggr\},
\end{align*}
where $t^{+}_{\mathrm{c}, i}(k,m,n)$ is the time slot when the price $p^{+}_{\mathrm{c}, i} (k,m)$ is run for the $n^{\mathrm{th}}$ time for the customer-side queue $i$, $t^{-}_{\mathrm{c}, i}(k,m,n)$ is the time slot when the price $p^{-}_{\mathrm{s}, j} (k,m)$ is run for the $n^{\mathrm{th}}$ time for the customer-side queue $i$, $t^{+}_{\mathrm{s}, j}(k,m,n)$ is the time slot when the price $p^{+}_{\mathrm{s}, j} (k,m)$ is run for the $n^{\mathrm{th}}$ time for the server-side queue $j$, and $t^{-}_{\mathrm{s}, j}(k,m,n)$ is the time slot when the price $p^{-}_{\mathrm{s}, j} (k,m)$ is run for the $n^{\mathrm{th}}$ time for the server-side queue $j$. The ``good'' event ${\cal C}$ means that all the estimated arrival rates are close to (concentrate to) the corresponding true arrival rates. Note that if $K< \lceil T/(2MN) \rceil$, we can create extra imaginary outer iterations that continue the same process so the definition of ${\cal C}$ is valid.

We first add the ``good'' event ${\cal C}$ to the regret definition \eqref{equ:regret-def} by the law of total expectation as follows:
\begin{align}
    & \expt [R(T)] 
    = \sum_{t=1}^{T}  \Biggl[ \biggl(\sum_i \lambda^*_i F_i(\lambda^*_i) - \sum_j \mu^*_j G_j(\mu^*_j) \biggr) 
    -   \expt \left[  \sum_i \lambda_i(t) F_i(\lambda_i(t)) - \sum_j \mu_j(t) G_j(\mu_j(t))  \right]  \Biggr]\nonumber\\
    = & \prob({\cal C})  \expt \Biggl[ \sum_{t=1}^{T} \biggl[ \biggl(\sum_i \lambda^*_i F_i(\lambda^*_i) - \sum_j \mu^*_j G_j(\mu^*_j) \biggr) 
    -    \biggl(  \sum_i \lambda_i(t) F_i(\lambda_i(t)) - \sum_j \mu_j(t) G_j(\mu_j(t))  \biggr) \biggr] \Biggl|\Biggr. {\cal C} \Biggr]\label{equ:regret-term-1-2}\\
    & + \prob({\cal C}^{\mathrm{c}}) \expt \Biggl[ \sum_{t=1}^{T} \biggl[ \biggl(\sum_i \lambda^*_i F_i(\lambda^*_i) - \sum_j \mu^*_j G_j(\mu^*_j) \biggr) 
    -    \biggl(  \sum_i \lambda_i(t) F_i(\lambda_i(t)) - \sum_j \mu_j(t) G_j(\mu_j(t))  \biggr) \biggr] \Biggl|\Biggr. {\cal C}^{\mathrm{c}} \Biggr].\label{equ:regret-term-1-3}
\end{align}
For the event ${\cal C}$, we have the following high probability lemma.
\begin{lemma}\label{lemma:concentration}
    $\prob \left( {\cal C}^{\mathrm{c}} \right) \le \Theta\left(T\epsilon^{\frac{\beta}{2}+1} + \epsilon^{\frac{\beta}{2}-1} \log (1/\epsilon) \right),$ where in the notation $\Theta(\cdot)$ we ignore the variables that do not depend on $T$.
\end{lemma}
Lemma~\ref{lemma:concentration} means that the ``good'' event ${\cal C}$ happens with high probability as long as $\epsilon$ is small.
Lemma~\ref{lemma:concentration} is meaningful when $\epsilon < T^{-\frac{2}{\beta+2}}$ and $\beta > 2$. 
Proof of Lemma~\ref{lemma:concentration} can be found in Appendix~\ref{app:proof-lemma-concentration}.

Since $\sum_i \lambda_i(t) F_i(\lambda_i(t)) \ge 0$ and $\sum_j \mu_j(t) G_j(\mu_j(t)) \le \sum_j G_j(1)$, we have
\begin{align}\label{equ:regret-term-1-1}
    & \biggl[\sum_i \lambda^*_i F_i(\lambda^*_i) - \sum_j \mu^*_j G_j(\mu^*_j)\biggr] - \biggl[\sum_i \lambda_i(t) F_i(\lambda_i(t)) - \sum_j \mu_j(t) G_j(\mu_j(t)) \biggr]\nonumber\\
    \le & \sum_i \lambda^*_i F_i(\lambda^*_i) - \sum_j \mu^*_j G_j(\mu^*_j) + \sum_j G_j(1).
\end{align}
Note that the term \eqref{equ:regret-term-1-1} is nonnegative since $\sum_j G_j(1) \ge \sum_j \mu^*_j G_j(\mu^*_j)$.
From Lemma~\ref{lemma:concentration} and \eqref{equ:regret-term-1-1}, the term \eqref{equ:regret-term-1-3} can be bounded by
\begin{align}\label{equ:regret-term-1-3-bound}
    \eqref{equ:regret-term-1-3} \le \Theta\left(T^2\epsilon^{\frac{\beta}{2}+1} + T\epsilon^{\frac{\beta}{2}-1} \log (1/\epsilon) \right).
\end{align}

\subsection{Bounding the Regret Caused by Rejecting Arrivals for Queue Length Control}
\label{app:bound-regret-reject-arr}

In this subsection, we will divide the term \eqref{equ:regret-term-1-2} into two cases according to whether the queue length is over the threshold $q^{\mathrm{th}}$ or not and then bound the regret caused by queue length control.

We define an event ${\cal H}_t$ as
\begin{align*}
    {\cal H}_t \coloneqq \left\{ Q_{\mathrm{c},i}(t) < q^{\mathrm{th}} \mbox{ for all } i, \mbox{ and } Q_{\mathrm{s},j}(t) < q^{\mathrm{th}} \mbox{ for all } j\right\}.
\end{align*}
This event ${\cal H}_t$ indicates that there is no queue length control at time $t$. We add an indicator of the event ${\cal H}_t$ into \eqref{equ:regret-term-1-2} as follows:
\begin{align}
    \eqref{equ:regret-term-1-2} & =  \prob({\cal C})  \expt \Biggl[ \sum_{t=1}^{T} \mathbb{1}_{{\cal H}_t^{\mathrm{c}}} \biggl[ \biggl(\sum_i \lambda^*_i F_i(\lambda^*_i) - \sum_j \mu^*_j G_j(\mu^*_j) \biggr) 
    -    \biggl(  \sum_i \lambda_i(t) F_i(\lambda_i(t)) - \sum_j \mu_j(t) G_j(\mu_j(t))  \biggr) \biggr] \Biggl|\Biggr. {\cal C} \Biggr]\nonumber\\
    & + \prob({\cal C})  \expt \Biggl[ \sum_{t=1}^{T} \mathbb{1}_{{\cal H}_t} \biggl[ \biggl(\sum_i \lambda^*_i F_i(\lambda^*_i) - \sum_j \mu^*_j G_j(\mu^*_j) \biggr) 
    -    \biggl(  \sum_i \lambda_i(t) F_i(\lambda_i(t)) - \sum_j \mu_j(t) G_j(\mu_j(t))  \biggr) \biggr] \Biggl|\Biggr. {\cal C} \Biggr]\nonumber\\
    & \le \biggl[\sum_i \lambda^*_i F_i(\lambda^*_i) - \sum_j \mu^*_j G_j(\mu^*_j) + \sum_j G_j(1)\biggr] \prob({\cal C})  \expt \Biggl[ \sum_{t=1}^{T} \mathbb{1}_{{\cal H}_t^{\mathrm{c}}} \Biggl|\Biggr. {\cal C} \Biggr]\label{equ:regret-term-1-4}\\
    & +  \expt \Biggl[ \sum_{t=1}^{T} \mathbb{1}_{{\cal H}_t} \biggl[ \biggl(\sum_i \lambda^*_i F_i(\lambda^*_i) - \sum_j \mu^*_j G_j(\mu^*_j) \biggr) 
    -    \biggl(  \sum_i \lambda_i(t) F_i(\lambda_i(t)) - \sum_j \mu_j(t) G_j(\mu_j(t))  \biggr) \biggr] \Biggl|\Biggr. {\cal C} \Biggr],\label{equ:regret-term-1-5}
\end{align}
where in the last inequality we use \eqref{equ:regret-term-1-1} and the fact that $\prob({\cal C})\le 1$. 
For the term \eqref{equ:regret-term-1-4}, by the law of total expectation and the union bound, we have
\begin{align}\label{equ:regret-term-1-6}
    & \prob({\cal C})  \expt \Biggl[ \sum_{t=1}^{T} \mathbb{1}_{{\cal H}_t^{\mathrm{c}}} \Biggl|\Biggr. {\cal C} \Biggr] 
    = \expt \Biggl[ \mathbb{1}_{\cal C} \sum_{t=1}^{T} \mathbb{1}_{{\cal H}_t^{\mathrm{c}}} \Biggr]
    - \prob({\cal C}^{\mathrm{c}})  \expt \Biggl[ \mathbb{1}_{\cal C} \sum_{t=1}^{T} \mathbb{1}_{{\cal H}_t^{\mathrm{c}}} \Biggl|\Biggr. {\cal C}^{\mathrm{c}} \Biggr] 
    = \expt \Biggl[ \mathbb{1}_{\cal C} \sum_{t=1}^{T} \mathbb{1}_{{\cal H}_t^{\mathrm{c}}} \Biggr] \nonumber\\
    = & \expt \Biggl[ \sum_{t=1}^{T}  \mathbb{1}_{\cal C} \mathbb{1} \left\{ \mbox{there exists } i \mbox{ such that } Q_{\mathrm{c},i}(t) \ge q^{\mathrm{th}} \mbox{ or } 
     \mbox{there exists } j \mbox{ such that } Q_{\mathrm{s},j}(t) \ge q^{\mathrm{th}}\right\}   \Biggr]\nonumber\\
    \le & \sum_{t=1}^{T} \sum_i \expt \biggl[  \mathbb{1}_{\cal C} \mathbb{1} \left\{ Q_{\mathrm{c},i}(t) \ge q^{\mathrm{th}}\right\} \biggr] 
    + \sum_{t=1}^{T} \sum_j \expt \biggl[  \mathbb{1}_{\cal C} \mathbb{1} \left\{ Q_{\mathrm{s},j}(t) \ge q^{\mathrm{th}}\right\}   \biggr]\nonumber\\
    = & \sum_{t=1}^{2MN} \sum_i \expt \biggl[  \mathbb{1}_{\cal C} \mathbb{1} \left\{ Q_{\mathrm{c},i}(t) \ge q^{\mathrm{th}}\right\} \biggr]
    + \sum_{t=2MN+1}^{T} \sum_i \expt \biggl[  \mathbb{1}_{\cal C} \mathbb{1} \left\{ Q_{\mathrm{c},i}(t) \ge q^{\mathrm{th}}\right\} \biggr]\nonumber\\
    & + \sum_{t=1}^{2MN} \sum_j \expt \biggl[  \mathbb{1}_{\cal C} \mathbb{1} \left\{ Q_{\mathrm{s},j}(t) \ge q^{\mathrm{th}}\right\}   \biggr]
    +\sum_{t=2MN+1}^{T} \sum_j \expt \biggl[  \mathbb{1}_{\cal C} \mathbb{1} \left\{ Q_{\mathrm{s},j}(t) \ge q^{\mathrm{th}}\right\}   \biggr],
\end{align}
where in the last equality, we divide each sum into two parts, corresponding to $k=1$ and $k>1$. 
In fact, there are $2MN$ time slots in the first outer iteration $k=1$ because there is no queue length control in the first outer iteration.
Since the length of each queue can increase at most $1$ at each time slot, for the first outer iteration, we have
\begin{align}\label{equ:regret-term-1-7}
    \sum_{t=1}^{2MN} \sum_i \expt \biggl[  \mathbb{1}_{\cal C} \mathbb{1} \left\{ Q_{\mathrm{c},i}(t) \ge q^{\mathrm{th}}\right\} \biggr] 
    \le I \sum_{t=1}^{2MN} \mathbb{1} \bigl\{ t \ge q^{\mathrm{th}}\bigr\} 
    = I \max \bigl\{0, 2MN - q^{\mathrm{th}}\bigr\},\nonumber\\
    \sum_{t=1}^{2MN} \sum_j \expt \biggl[  \mathbb{1}_{\cal C} \mathbb{1} \left\{ Q_{\mathrm{s},j}(t) \ge q^{\mathrm{th}}\right\}   \biggr]
    \le J \sum_{t=1}^{2MN} \mathbb{1} \bigl\{ t \ge q^{\mathrm{th}}\bigr\} 
    = J \max \bigl\{0, 2MN - q^{\mathrm{th}}\bigr\}.
\end{align}
For the terms $\sum_{t=2MN+1}^{T} \sum_i \expt \left[  \mathbb{1}_{\cal C} \mathbb{1} \left\{ Q_{\mathrm{c},i}(t) \ge q^{\mathrm{th}}\right\} \right]$ and $\sum_{t=2MN+1}^{T} \sum_j \expt \left[  \mathbb{1}_{\cal C} \mathbb{1} \left\{ Q_{\mathrm{s},j}(t) \ge q^{\mathrm{th}}\right\}   \right]$ in \eqref{equ:regret-term-1-6}, we have the following lemma.
\begin{lemma}\label{lemma:bound-time-large-queue}
    Let Assumption~\ref{assum:1} and Assumption~\ref{assum:4} hold. Suppose $\epsilon < \delta$ and $T$ is sufficiently large. Then for the customer side, we have
    \begin{align*}
    &  \sum_{t=2MN+1}^{T} \sum_i \expt \Bigl[ 
    \mathbb{1}_{{\cal C}}
    \mathbb{1} \Bigl\{ Q_{\mathrm{c},i}(t) \ge q^{\mathrm{th}} \Bigr\}
    \Bigr]\nonumber\\
    \le & \Theta \left(\frac{\log^4(1/\epsilon)}{\epsilon^4 q^{\mathrm{th}}}
    + \frac{T}{q^{\mathrm{th}}}
    + T^2\epsilon^{\frac{\beta}{2}+1} + T\epsilon^{\frac{\beta}{2}-1} \log (1/\epsilon) 
    + T \left(\frac{\eta \epsilon}{\delta} + \eta + \delta + \epsilon\right) \right),
    \end{align*}
    and for the server side, we have
    \begin{align*}
    &  \sum_{t=2MN+1}^{T} \sum_j \expt \Bigl[ 
    \mathbb{1}_{{\cal C}}
    \mathbb{1} \Bigl\{ Q_{\mathrm{s},j}(t) \ge q^{\mathrm{th}} \Bigr\}
    \Bigr]\nonumber\\
    \le & \Theta \left(\frac{\log^4(1/\epsilon)}{\epsilon^4 q^{\mathrm{th}}}
    + \frac{T}{q^{\mathrm{th}}}
    + T^2\epsilon^{\frac{\beta}{2}+1} + T\epsilon^{\frac{\beta}{2}-1} \log (1/\epsilon) 
    + T \left(\frac{\eta \epsilon}{\delta} + \eta + \delta + \epsilon\right) \right),
    \end{align*}
    where in the notation $\Theta(\cdot)$ we ignore the variables that do not depend on $T$.
\end{lemma}
\noindent Lemma~\ref{lemma:bound-time-large-queue} provides upper bounds for the expected number of time slots when the queue length is greater than or equal to $q^{\mathrm{th}}$ and the ``good'' event ${\cal C}$ holds. 
Proof of Lemma~\ref{lemma:bound-time-large-queue} can be found in Appendix~\ref{app:proof-lemma-bound-time-large-queue}.
From \eqref{equ:regret-term-1-6}, \eqref{equ:regret-term-1-7}, and Lemma~\ref{lemma:bound-time-large-queue}, we have
\begin{align}\label{equ:regret-term-1-8}
    & \prob({\cal C})  \expt \Biggl[ \sum_{t=1}^{T} \mathbb{1}_{{\cal H}_t^{\mathrm{c}}} \Biggl|\Biggr. {\cal C} \Biggr]
    \le \Theta \left(\max \biggl\{0, 2\left\lceil \log_2 (1/\epsilon) \right\rceil \left\lceil \frac{\beta \ln (1/\epsilon)}{\epsilon^2} \right\rceil - q^{\mathrm{th}}\biggr\} \right)\nonumber\\
    & + \Theta \left(\frac{\log^4(1/\epsilon)}{\epsilon^4 q^{\mathrm{th}}}
    + \frac{T}{q^{\mathrm{th}}}
    + T^2\epsilon^{\frac{\beta}{2}+1} + T\epsilon^{\frac{\beta}{2}-1} \log (1/\epsilon) 
    + T \left(\frac{\eta \epsilon}{\delta} + \eta + \delta + \epsilon\right) \right).
\end{align}

\subsection{Bounding the Regret during Optimization}
\label{app:regret-zero-order-opt}

\begin{figure}[tb]
    \centering
    \includegraphics[width=0.6\linewidth]{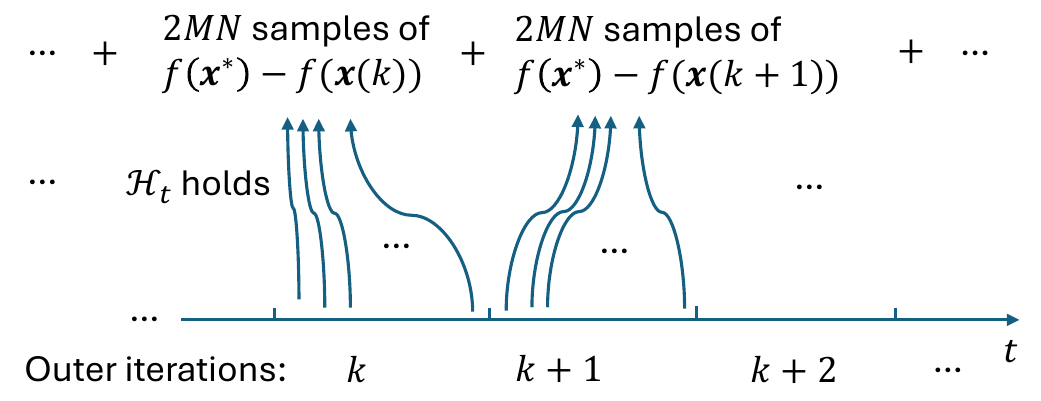}
    \caption{Transforming the summation.}
    \label{fig:transform}
\end{figure}

In this subsection, we will bound the term \eqref{equ:regret-term-1-5}. \eqref{equ:regret-term-1-5} can be divided into two terms as follows:
\begin{align*}
    \eqref{equ:regret-term-1-5} = \expt \Biggl[ \sum_{t=1}^{2MN} \mathbb{1}_{{\cal H}_t} \biggl[ \biggl(\sum_i \lambda^*_i F_i(\lambda^*_i) - \sum_j \mu^*_j G_j(\mu^*_j) \biggr) 
    -    \biggl(  \sum_i \lambda_i(t) F_i(\lambda_i(t)) - \sum_j \mu_j(t) G_j(\mu_j(t))  \biggr) \biggr] \Biggl|\Biggr. {\cal C} \Biggr]\nonumber\\
    + \expt \Biggl[ \sum_{t=2MN+1}^{T} \mathbb{1}_{{\cal H}_t} \biggl[ \biggl(\sum_i \lambda^*_i F_i(\lambda^*_i) - \sum_j \mu^*_j G_j(\mu^*_j) \biggr) 
    -    \biggl(  \sum_i \lambda_i(t) F_i(\lambda_i(t)) - \sum_j \mu_j(t) G_j(\mu_j(t))  \biggr) \biggr] \Biggl|\Biggr. {\cal C} \Biggr].
\end{align*}
By \eqref{equ:regret-term-1-1}, we further obtain
\begin{align}\label{equ:regret-term-1-9-temp}
    & \eqref{equ:regret-term-1-5} \le 2MN \biggl[ \sum_i \lambda^*_i F_i(\lambda^*_i) - \sum_j \mu^*_j G_j(\mu^*_j) + \sum_j G_j(1) \biggr]\nonumber\\
    & + \expt \Biggl[ \sum_{t=2MN+1}^{T} \mathbb{1}_{{\cal H}_t} \biggl[ \biggl(\sum_i \lambda^*_i F_i(\lambda^*_i) - \sum_j \mu^*_j G_j(\mu^*_j) \biggr) 
    - \biggl(  \sum_i \lambda_i(t) F_i(\lambda_i(t)) - \sum_j \mu_j(t) G_j(\mu_j(t))  \biggr) \biggr] \Biggl|\Biggr. {\cal C} \Biggr].
\end{align}
Recall that $\lambda_i(t),\mu_j(t)$ are the actual arrival rates at time $t$. Define $\lambda'_i(t)\coloneqq F^{-1}_i(p^{+/-}_{\mathrm{c},i}(k(t),m(t)))$, where $k(t),m(t)$ denote the outer iteration and bisection iteration at time $t$, respectively. Similarly define $\mu'_j(t)\coloneqq G^{-1}_j(p^{+/-}_{\mathrm{s},j}(k(t),m(t)))$. In fact, $\lambda'_i(t),\mu'_j(t)$ are the arrival rates as if there is no queue length control.
Under the event ${\cal H}_t$, there is no queue length control, which implies that $\lambda_i(t) = \lambda'_i(t)$ and $\mu_j(t)=\mu'_j(t)$. Hence, we can change $\lambda(t),\mu(t)$ in \eqref{equ:regret-term-1-9-temp} to $\lambda'(t),\mu'(t)$ as follows:
\begin{align}\label{equ:regret-term-1-9}
    & \eqref{equ:regret-term-1-5} \le 2MN \biggl[ \sum_i \lambda^*_i F_i(\lambda^*_i) - \sum_j \mu^*_j G_j(\mu^*_j) + \sum_j G_j(1) \biggr]\nonumber\\
    & + \expt \Biggl[ \sum_{t=2MN+1}^{T} \mathbb{1}_{{\cal H}_t} \biggl[ \biggl(\sum_i \lambda^*_i F_i(\lambda^*_i) - \sum_j \mu^*_j G_j(\mu^*_j) \biggr) 
    - \biggl(  \sum_i \lambda'_i(t) F_i(\lambda'_i(t)) - \sum_j \mu'_j(t) G_j(\mu'_j(t))  \biggr) \biggr] \Biggl|\Biggr. {\cal C} \Biggr].
\end{align}
Next, we will relate $\lambda'_i(t)$ and $\mu'_j(t)$ to the arrival rates $\sum_{j\in {\cal E}_{\mathrm{c},i}} x_{i,j}(k(t))$ and $\sum_{i\in {\cal E}_{\mathrm{s},j}} x_{i,j}(k(t))$, which are the target arrival rates generated in the optimization process.
Notice that by Assumption~\ref{assum:1}, we have
\begin{align}\label{equ:diff-bound-f-1}
    & \lambda_1 F_i(\lambda_1) - \lambda_2 F_i(\lambda_2)
    = (\lambda_1 F_i(\lambda_1)  - \lambda_1 F_i(\lambda_2)) + (\lambda_1 F_i(\lambda_2) - \lambda_2 F_i(\lambda_2))\nonumber\\
    \le & \lambda_1 L_{F_i}|\lambda_1 - \lambda_2| + |\lambda_1 - \lambda_2| F_i(\lambda_2)
    \le (L_{F_i} + p_{\mathrm{c},i,\max})|\lambda_1 - \lambda_2|,
\end{align}
for any $\lambda_1, \lambda_2\in[0,1]$. Similarly, for any $\mu_1, \mu_2\in[0,1]$, we have
\begin{align}\label{equ:diff-bound-f-2}
    \mu_1 G_j(\mu_1) - \mu_2 G_j(\mu_2)
    \le (L_{G_j} + p_{\mathrm{s},j,\max})|\mu_1 - \mu_2|.
\end{align}
Then by \eqref{equ:diff-bound-f-1}, we have
\begin{align}\label{equ:diff-bound-f-3-temp}
     & \sum_i \biggl(\sum_{j\in {\cal E}_{\mathrm{c},i}} x_{i,j}(k(t))\biggr)
     F_i\biggl(\sum_{j\in {\cal E}_{\mathrm{c},i}} x_{i,j}(k(t))\biggr) 
     - \sum_i \lambda'_i(t) F_i(\lambda'_i(t))\nonumber\\
     \le & \sum_i (L_{F_i} + p_{\mathrm{c},i,\max})
     \biggl|\lambda'_i(t) - \sum_{j\in {\cal E}_{\mathrm{c},i}} x_{i,j}(k(t))\biggr|.
\end{align}
Lemma~\ref{lemma:bisection-error} in Appendix~\ref{app:proof-lemma-bound-time-large-queue} bounds the difference between the actual arrival rate of any queue without queue length control and the target arrival rate generated in the optimization process, i.e., $|\lambda'_i(t) - \sum_{j\in {\cal E}_{\mathrm{c},i}} x_{i,j}(k(t))|$ is bounded and small. Hence, by \eqref{equ:diff-bound-f-3-temp} and Lemma~\ref{lemma:bisection-error}, we have
\begin{align}\label{equ:diff-bound-f-3}
     & \sum_i \biggl(\sum_{j\in {\cal E}_{\mathrm{c},i}} x_{i,j}(k(t))\biggr)
     F_i\biggl(\sum_{j\in {\cal E}_{\mathrm{c},i}} x_{i,j}(k(t))\biggr) 
     - \sum_i \lambda'_i(t) F_i(\lambda'_i(t))\nonumber\\
     \le & \sum_i (L_{F_i} + p_{\mathrm{c},i,\max})
     \biggl[ 2 L_{F^{-1}_i} e_{\mathrm{c},i} + 2\epsilon \left[1 + L_{F^{-1}_i} \left( p_{\mathrm{c},i,\max} - p_{\mathrm{c},i,\min}\right)\right]  + \delta \sqrt{|{\cal E}_{\mathrm{c},i}|} \biggr]\nonumber\\
     = & \Theta\biggl(\frac{\eta \epsilon}{\delta} + \eta + \delta + \epsilon\biggr).
\end{align}
Similarly, by \eqref{equ:diff-bound-f-2} and Lemma~\ref{lemma:bisection-error}, we have
\begin{align}\label{equ:diff-bound-f-4}
    \sum_j \mu'_j(t) G_j(\mu'_j(t))
    - \sum_j \biggl( \sum_{i\in {\cal E}_{\mathrm{s},j}} x_{i,j}(k(t)) \biggr)
    G_j\biggl( \sum_{i\in {\cal E}_{\mathrm{s},j}} x_{i,j}(k(t)) \biggr)
    \le \Theta\biggl(\frac{\eta \epsilon}{\delta} + \eta + \delta + \epsilon\biggr).
\end{align}
By \eqref{equ:regret-term-1-9}, \eqref{equ:diff-bound-f-3}, and \eqref{equ:diff-bound-f-4}, we have
\begin{align}\label{equ:regret-term-1-10}
    \eqref{equ:regret-term-1-5} \le 
    \Theta\biggl( \frac{\log^2 (1/\epsilon)}{\epsilon^2} + T \biggl(\frac{\eta \epsilon}{\delta} + \eta + \delta + \epsilon \biggr) \biggr)
    + \expt \Biggl[ \sum_{t=2MN+1}^{T} \mathbb{1}_{{\cal H}_t}
    \bigl[f(\boldsymbol{x}^*) - f(\boldsymbol{x}(k(t)))\bigr]
    \Biggl|\Biggr. {\cal C} \Biggr],
\end{align}
where we recall the function $f: {\cal D} \rightarrow \mathbb{R}$ is defined as follows:
\begin{align*}
    f(\boldsymbol{x}) \coloneqq \sum_i \biggl(\sum_{j\in {\cal E}_{\mathrm{c},i}} x_{i,j}\biggr) F_i\biggl(\sum_{j\in {\cal E}_{\mathrm{c},i}} x_{i,j}\biggr) 
    - \sum_j \biggl( \sum_{i\in {\cal E}_{\mathrm{s},j}} x_{i,j} \biggr) G_j\biggl( \sum_{i\in {\cal E}_{\mathrm{s},j}} x_{i,j} \biggr).
\end{align*}
The term $\expt [ \sum_{t=2MN+1}^{T} \mathbb{1}_{{\cal H}_t}[f(\boldsymbol{x}^*) - f(\boldsymbol{x}(k(t)))] | {\cal C} ]$ in \eqref{equ:regret-term-1-10} can be viewed as the regret induced during the optimization process.
By the definition of ${\cal H}_t$, the summation in the term $\expt [ \sum_{t=2MN+1}^{T} \mathbb{1}_{{\cal H}_t}[f(\boldsymbol{x}^*) - f(\boldsymbol{x}(k(t)))] | {\cal C} ]$ is taken over the time slots when the queue length is less than $q^{\mathrm{th}}$ for all queues. We know that the number of these time slots in one outer iteration should be $2MN$. Also notice that $f(\boldsymbol{x}^*) - f(\boldsymbol{x}(k(t)))$ has the same value for all $t$ in the same outer iteration. There are $k(T)$ outer iterations in total and in each outer iteration there are $2MN$ samples that count into the summation. Therefore, as shown in Figure~\ref{fig:transform}, by taking $2MN$ samples from each outer iteration and then summing over all outer iterations, we can transform the term $\expt [ \sum_{t=2MN+1}^{T} \mathbb{1}_{{\cal H}_t}[f(\boldsymbol{x}^*) - f(\boldsymbol{x}(k(t)))] | {\cal C} ]$ into
\begin{align}\label{equ:regret-term-1-11}
    \expt \Biggl[ \sum_{t=2MN+1}^{T} \mathbb{1}_{{\cal H}_t}
    \bigl[f(\boldsymbol{x}^*) - f(\boldsymbol{x}(k(t)))\bigr]
    \Biggl|\Biggr. {\cal C} \Biggr]
    = & 2MN \expt \Biggl[ \sum_{k=2}^{k(T)}
    \bigl[f(\boldsymbol{x}^*) - f(\boldsymbol{x}(k))\bigr]
    \Biggl|\Biggr. {\cal C} \Biggr].
\end{align}
Note that $\boldsymbol{x}^* \in \argmax_{\boldsymbol{x}\in {\cal D}} f(\boldsymbol{x})$ and $\boldsymbol{x}(k)\in {\cal D}' \subseteq {\cal D}$ by Lemma~\ref{lemma:shrunk-set}. Hence, $f(\boldsymbol{x}^*) \ge f(\boldsymbol{x}(k))$. Since there are at least $2MN$ time slots in each outer iteration, $k(T)\le T/(2MN)$. Hence, \eqref{equ:regret-term-1-11} can be bounded by
\begin{align}\label{equ:regret-term-1-12}
    \expt \Biggl[ \sum_{t=2MN+1}^{T} \mathbb{1}_{{\cal H}_t}
    \bigl[f(\boldsymbol{x}^*) - f(\boldsymbol{x}(k(t)))\bigr]
    \Biggl|\Biggr. {\cal C} \Biggr]    \le & 2MN  \sum_{k=1}^{\lceil T/(2MN) \rceil} \expt \Bigl[
    f(\boldsymbol{x}^*) - f(\boldsymbol{x}(k))
    \Bigl|\Bigr. {\cal C} \Bigr].
\end{align}
For the term $\sum_{k=1}^{\lceil T/(2MN) \rceil} \expt [f(\boldsymbol{x}^*) - f(\boldsymbol{x}(k))| {\cal C} ]$ in \eqref{equ:regret-term-1-12}, we have the following lemma:
\begin{lemma}\label{lemma:optimization}
    Let Assumption~\ref{assum:1}--\ref{assum:4} hold. Then for any fixed integer $K$, we have
    \begin{align*}
        \sum_{k=1}^{K} \expt \Bigl[f(\boldsymbol{x}^*) - f(\boldsymbol{x}(k))
        \Bigl|\Bigr. {\cal C} \Bigr] 
        \le  \Theta \left( \frac{1}{\eta} +  K\eta + \frac{K\eta \epsilon^2}{\delta^2} + \frac{K\epsilon}{\delta} + K\delta\right),
    \end{align*}
    where in the notation $\Theta(\cdot)$ we ignore the variables that do not depend on $T$ or $K$.
\end{lemma}
Proof of Lemma~\ref{lemma:optimization} can be found in Appendix~\ref{app:proof-lemma-optimization}. The proof uses the techniques of gradient-free stochastic projected gradient descent \cite{agarwal2010optimal}. One difference is that our estimated gradient is biased and we need to bound the bias. Another difference is that the expectation is conditioned on the event ${\cal C}$. Since the event ${\cal C}$ is carefully designed, ${\cal C}$ is independent of $\boldsymbol{u}(k)$, which is essential in our proof.
By Lemma~\ref{lemma:optimization}, \eqref{equ:regret-term-1-12}, $M=\left\lceil\log_2(1/\epsilon)\right\rceil$, $N=\left\lceil\frac{\beta \ln (1/\epsilon)}{\epsilon^2}\right\rceil$, $\eta<1$, and $\delta<1$, we have
\begin{align}\label{equ:regret-term-1-13}
    & \expt \Biggl[ \sum_{t=2MN+1}^{T} \mathbb{1}_{{\cal H}_t}
    \bigl[f(\boldsymbol{x}^*) - f(\boldsymbol{x}(k(t)))\bigr]
    \Biggl|\Biggr. {\cal C} \Biggr]\nonumber\\    
    \le & \Theta \left(  T\eta + \frac{T\eta \epsilon^2}{\delta^2} + \frac{T\epsilon}{\delta} + T\delta + \left( \frac{1}{\eta \epsilon^2 } + \frac{\eta}{\delta^2}
    + \frac{1}{\epsilon \delta} \right) \log^2(1/\epsilon) \right).
\end{align}
By \eqref{equ:regret-term-1-10}, \eqref{equ:regret-term-1-13}, $\eta<1$, and $\delta<1$, we have
\begin{align}\label{equ:regret-term-1-14}
    \eqref{equ:regret-term-1-5} 
    \le & \Theta\left(  \frac{T\epsilon}{\delta} + T\eta + T\delta + \frac{T\eta \epsilon^2}{\delta^2} 
    + \log^2 (1/\epsilon) \left(\frac{1}{\eta \epsilon^2} + \frac{\eta}{\delta^2} + \frac{1}{\epsilon \delta}\right) \right),
\end{align}

Combining \eqref{equ:regret-term-1-4}, \eqref{equ:regret-term-1-5}, \eqref{equ:regret-term-1-8}, and \eqref{equ:regret-term-1-14}, we obtain
\begin{align}\label{equ:regret-term-1-15}
    \eqref{equ:regret-term-1-2}
    \le &   
    \Theta \Biggl(
    \frac{\log^4(1/\epsilon)}{\epsilon^4 q^{\mathrm{th}}}
    + \frac{T}{q^{\mathrm{th}}}
    + T^2\epsilon^{\frac{\beta}{2}+1} + T\epsilon^{\frac{\beta}{2}-1} \log (1/\epsilon) 
    + T \left(\frac{\eta \epsilon}{\delta} + \eta + \delta + \epsilon\right)
    \nonumber\\
    & \quad  + \frac{T\epsilon}{\delta} + \frac{T\eta \epsilon^2}{\delta^2}
    + \log^2 (1/\epsilon) \left(\frac{1}{\eta \epsilon^2} + \frac{\eta}{\delta^2} + \frac{1}{\epsilon \delta}\right)
    \Biggr).
\end{align}
From \eqref{equ:regret-term-1-2}, \eqref{equ:regret-term-1-3}, \eqref{equ:regret-term-1-3-bound}, and \eqref{equ:regret-term-1-15}, we have
\begin{align}\label{equ:regret-term-1-final}
    \expt [R(T)]
    \le &
    \Theta \Biggl(
    \frac{\log^4(1/\epsilon)}{\epsilon^4 q^{\mathrm{th}}}
    + \frac{T}{q^{\mathrm{th}}}
    + T^2\epsilon^{\frac{\beta}{2}+1} + T\epsilon^{\frac{\beta}{2}-1} \log (1/\epsilon)
    + T \left(\frac{\eta \epsilon}{\delta} + \eta + \delta + \epsilon\right)
    \nonumber\\
    & \quad + \frac{T\epsilon}{\delta}  + \frac{T\eta \epsilon^2}{\delta^2}
    + \log^2 (1/\epsilon) \left(\frac{1}{\eta \epsilon^2} + \frac{\eta}{\delta^2} + \frac{1}{\epsilon \delta}\right)
    \Biggr).
\end{align}
Notice that if $\epsilon\ge \delta$, the regret bound will be at least linear in $T$. Hence, we will choose $\epsilon, \delta$ such that $\epsilon < \delta$. Therefore, the bound \eqref{equ:regret-term-1-final} can be reduced to \eqref{equ:regret-bound} in Theorem~\ref{theo:1}.

\subsection{Bounding the Maximum Queue Length}

Since the length of each queue can increase at most $1$ at each time slot and there are $2MN$ time slots in the first outer iteration ($k=1$), the length of each queue can increase to at most $2MN$ after the first outer iteration. When $k>1$, since for each queue we set the arrival rate to be zero when the queue length is greater than or equal to $q^{\mathrm{th}}$, the queue length must be less than or equal to 
$\max\{2MN,q^{\mathrm{th}}\}$. Substituting $M$ and $N$ with their values, we obtain the maximum queue-length bound \eqref{equ:queue-length-bound} in Theorem~\ref{theo:1}.

\subsection{Bounding the Average Queue Length}
To bounding the average queue length, we will use an intermediate result~\eqref{equ:drift-intermediate} in the proof of Lemma~\ref{lemma:bound-time-large-queue} in Appendix~\ref{app:proof-lemma-bound-time-large-queue}. 

From \eqref{equ:drift-intermediate}, we have
\begin{align*}
    & \expt [V_{\mathrm{c}}(t+1) - V_{\mathrm{c}}(t) ] \nonumber\\
    \le & 
    - a_{\min} \sum_i \expt \Bigl[ 
    \mathbb{1}_{{\cal C}}
    Q_{\mathrm{c},i}(t)
    \mathbb{1} \Bigl\{ Q_{\mathrm{c},i}(t) \ge q^{\mathrm{th}} \Bigr\}
    \Bigr]\nonumber\\
    & + q^{\mathrm{th}} \Theta\left(\frac{\eta \epsilon}{\delta} + \eta + \delta + \epsilon\right)
    + 2 I q^{\mathrm{th}} \expt \bigl[ \mathbb{1}_{{\cal C}^{\mathrm{c}}} \bigr] + \Theta(1)
\end{align*}
for all $t\ge 2MN+1$.
Note that $\expt [ \mathbb{1}_{{\cal C}^{\mathrm{c}}} ] = \prob ( {\cal C}^{\mathrm{c}}) \le \Theta\left(T\epsilon^{\frac{\beta}{2}+1} + \epsilon^{\frac{\beta}{2}-1} \log (1/\epsilon) \right) $ by Lemma~\ref{lemma:concentration}. Then we have
\begin{align*}
    & \expt [V_{\mathrm{c}}(t+1) - V_{\mathrm{c}}(t) ] \nonumber\\
    \le & 
    - a_{\min} \sum_i \expt \Bigl[ 
    \mathbb{1}_{{\cal C}}
    Q_{\mathrm{c},i}(t)
    \mathbb{1} \Bigl\{ Q_{\mathrm{c},i}(t) \ge q^{\mathrm{th}} \Bigr\}
    \Bigr]\nonumber\\
    & + q^{\mathrm{th}} \Theta\left(\frac{\eta \epsilon}{\delta} + \eta + \delta + \epsilon\right)
    + q^{\mathrm{th}} \Theta\left(T\epsilon^{\frac{\beta}{2}+1} + \epsilon^{\frac{\beta}{2}-1} \log (1/\epsilon) \right) + \Theta(1)
\end{align*}
for all $t\ge 2MN+1$.
Taking summation over $t$ from $2MN+1$ to $T$, we have
\begin{align*}
    & \sum_{t=2MN+1}^{T} a_{\min} \sum_i \expt \Bigl[ 
    \mathbb{1}_{{\cal C}}
    Q_{\mathrm{c},i}(t)
    \mathbb{1} \Bigl\{ Q_{\mathrm{c},i}(t) \ge q^{\mathrm{th}} \Bigr\}
    \Bigr]\nonumber\\
    \le & \expt [V_{\mathrm{c}}(2MN+1) ]
    + T q^{\mathrm{th}} \Theta\left(\frac{\eta \epsilon}{\delta} + \eta + \delta + \epsilon\right)
    + q^{\mathrm{th}} \Theta\left(T^2\epsilon^{\frac{\beta}{2}+1} + T\epsilon^{\frac{\beta}{2}-1} \log (1/\epsilon) \right) + \Theta(T).
\end{align*}
Note that we do not reject arrivals to control the queue lengths during the first outer iteration ($k=1$) and control it with threshold $q^{\mathrm{th}}$ when $k\ge 2$. Therefore, we have $V_{\mathrm{c}}(2MN+1) = \sum_i Q^2_{\mathrm{c},i}(2MN+1) \le 4 I M^2 N^2$. Then we have
\begin{align}\label{equ:avg-queue-length-interm-1}
    & \sum_{t=2MN+1}^{T} \sum_i \expt \Bigl[ 
    \mathbb{1}_{{\cal C}}
    Q_{\mathrm{c},i}(t)
    \mathbb{1} \Bigl\{ Q_{\mathrm{c},i}(t) \ge q^{\mathrm{th}} \Bigr\}
    \Bigr]\nonumber\\
    \le & \Theta (M^2 N^2)
    + T q^{\mathrm{th}} \Theta\left(\frac{\eta \epsilon}{\delta} + \eta + \delta + \epsilon\right)
    + q^{\mathrm{th}} \Theta\left(T^2 \epsilon^{\frac{\beta}{2}+1} + T \epsilon^{\frac{\beta}{2}-1} \log (1/\epsilon) \right) + \Theta(T).
\end{align}
Note that
\begin{align*}
    & \sum_{t=1}^{T} \expt \biggl[\sum_i Q_{\mathrm{c},i}(t)\biggr]\nonumber\\
    = & \sum_{t=1}^{2MN} \expt \biggl[\sum_i Q_{\mathrm{c},i}(t)\biggr]
    + \sum_{t=2MN+1}^{T} \expt \biggl[\sum_i Q_{\mathrm{c},i}(t)\biggr]\nonumber\\
    \le & \sum_{t=1}^{2MN} It
    + \sum_{t=2MN+1}^{T} \expt \biggl[\sum_i Q_{\mathrm{c},i}(t)\biggr]\nonumber\\
    = & \Theta(M^2 N^2) 
    + \sum_{t=2MN+1}^{T} \expt \biggl[\sum_i Q_{\mathrm{c},i}(t) \mathbb{1} \Bigl\{ Q_{\mathrm{c},i}(t) \ge q^{\mathrm{th}} \Bigr\} \biggr]
    + \sum_{t=2MN+1}^{T} \expt \biggl[\sum_i Q_{\mathrm{c},i}(t) \mathbb{1} \Bigl\{ Q_{\mathrm{c},i}(t) < q^{\mathrm{th}} \Bigr\} \biggr]
    \nonumber\\
    \le & \Theta(M^2 N^2) + \Theta(T q^{\mathrm{th}})
    + \sum_{t=2MN+1}^{T} \expt \biggl[\sum_i Q_{\mathrm{c},i}(t) \mathbb{1} \Bigl\{ Q_{\mathrm{c},i}(t) \ge q^{\mathrm{th}} \Bigr\} \biggr].
\end{align*}
Adding the concentration event ${\cal C}$, we have
\begin{align*}
    \sum_{t=1}^{T} \expt \biggl[\sum_i Q_{\mathrm{c},i}(t)\biggr]
    \le & \Theta(M^2 N^2) + \Theta(T q^{\mathrm{th}})
    + \sum_{t=2MN+1}^{T} \expt \biggl[
    \mathbb{1}_{{\cal C}}
    \sum_i Q_{\mathrm{c},i}(t) \mathbb{1} \Bigl\{ Q_{\mathrm{c},i}(t) \ge q^{\mathrm{th}} \Bigr\} \biggr]\nonumber\\
    & + \sum_{t=2MN+1}^{T} \expt \biggl[
    \mathbb{1}_{{\cal C}^\mathrm{c}}
    \sum_i Q_{\mathrm{c},i}(t) \mathbb{1} \Bigl\{ Q_{\mathrm{c},i}(t) \ge q^{\mathrm{th}} \Bigr\} \biggr].
\end{align*}
By the maximum queue-length bound \eqref{equ:queue-length-bound}, we have
\begin{align*}
    \sum_{t=1}^{T} \expt \biggl[\sum_i Q_{\mathrm{c},i}(t)\biggr]
    \le & \Theta(M^2 N^2) + \Theta(T q^{\mathrm{th}})
    + \sum_{t=2MN+1}^{T} \expt \biggl[
    \mathbb{1}_{{\cal C}}
    \sum_i Q_{\mathrm{c},i}(t) \mathbb{1} \Bigl\{ Q_{\mathrm{c},i}(t) \ge q^{\mathrm{th}} \Bigr\} \biggr]\nonumber\\
    & + T I \max\{2MN, q^{\mathrm{th}}\} \expt \bigl[
    \mathbb{1}_{{\cal C}^\mathrm{c}} \bigr].
\end{align*}
By Lemma~\ref{lemma:concentration}, we have
\begin{align}\label{equ:avg-queue-length-interm-2}
    \sum_{t=1}^{T} \expt \biggl[\sum_i Q_{\mathrm{c},i}(t)\biggr]
    \le & \Theta(M^2 N^2) + \Theta(T q^{\mathrm{th}})
    + \sum_{t=2MN+1}^{T} \expt \biggl[
    \mathbb{1}_{{\cal C}}
    \sum_i Q_{\mathrm{c},i}(t) \mathbb{1} \Bigl\{ Q_{\mathrm{c},i}(t) \ge q^{\mathrm{th}} \Bigr\} \biggr]\nonumber\\
    & + \max\{2MN, q^{\mathrm{th}}\} \Theta\left(T^2 \epsilon^{\frac{\beta}{2}+1} + T \epsilon^{\frac{\beta}{2}-1} \log (1/\epsilon) \right).
\end{align}
From \eqref{equ:avg-queue-length-interm-1} and \eqref{equ:avg-queue-length-interm-2}, we have
\begin{align*}
    \sum_{t=1}^{T} \expt \biggl[\sum_i Q_{\mathrm{c},i}(t)\biggr]
    \le & \Theta(M^2 N^2) + \Theta(T q^{\mathrm{th}})
    + T q^{\mathrm{th}} \Theta\left(\frac{\eta \epsilon}{\delta} + \eta + \delta + \epsilon\right)
    + \Theta(T)
    \nonumber\\
    & + \max\{2MN, q^{\mathrm{th}}\} \Theta\left(T^2 \epsilon^{\frac{\beta}{2}+1} + T \epsilon^{\frac{\beta}{2}-1} \log (1/\epsilon) \right).
\end{align*}
Divided by $T$, we have
\begin{align*}
    \frac{1}{T}\sum_{t=1}^{T} \expt \biggl[\sum_i Q_{\mathrm{c},i}(t)\biggr]
    \le & \Theta(\frac{M^2 N^2}{T}) + \Theta(q^{\mathrm{th}})
    + q^{\mathrm{th}} \Theta\left(\frac{\eta \epsilon}{\delta} + \eta + \delta + \epsilon\right)
    \nonumber\\
    & + \max\{2MN, q^{\mathrm{th}}\} \Theta\left(T \epsilon^{\frac{\beta}{2}+1} + \epsilon^{\frac{\beta}{2}-1} \log (1/\epsilon) \right)\nonumber\\
    \le & \Theta(\frac{M^2 N^2}{T}) + \Theta(q^{\mathrm{th}})
    + \max\{2MN, q^{\mathrm{th}}\} \Theta\left(T \epsilon^{\frac{\beta}{2}+1} + \epsilon^{\frac{\beta}{2}-1} \log (1/\epsilon) \right),
\end{align*}
where the last inequality is due to the fact that $\epsilon < \delta$ and $\epsilon, \delta,\eta < 1$. Following the same argument, we can obtain a bound of the same order for $\frac{1}{T}\sum_{t=1}^{T} \expt \biggl[\sum_j Q_{\mathrm{s},j}(t)\biggr]$.
Substituting the definitions of $M$ and $N$ into the above bounds, we obtain the average queue-length bound \eqref{equ:avg-queue-length-bound} in Theorem~\ref{theo:1}.

Theorem~\ref{theo:1} is proved.

\section{Proof of Theorem~\ref{theo:2}}
\label{app:theo:2}

In this section, we present the complete proof of Theorem~\ref{theo:2}.
The proof of Theorem~\ref{theo:2} is similar to that of Theorem~\ref{theo:1}, but the main difference is bounding the regret caused by rejecting arrivals for queue length control (Lemma~\ref{lemma:bound-time-large-queue}). In fact, with Assumption~\ref{assum:5}, we can obtain a better result, where the regret caused by rejecting arrivals in the first outer iteration will be much smaller and the queue length at the end of the first outer iteration will also be much smaller since we have a set of arrival rates that remain nearly balanced during the first outer iteration.

\subsection{Adding the Event of Concentration}
Recall the definition of the ``good'' event $\mathcal{C}$, which means that the average of the arrival samples is approximately equal to the arrival rate for all arrival rates. Adding the event $\mathcal{C}$ to the regret definition~\eqref{equ:regret-def} and following the same argument as that in Section~\ref{sec:roadmap-add-concentration}, we have
\begin{align}
    & \expt [R(T)] 
    = \sum_{t=1}^{T}  \Biggl[ \biggl(\sum_i \lambda^*_i F_i(\lambda^*_i) - \sum_j \mu^*_j G_j(\mu^*_j) \biggr) 
    -   \expt \left[  \sum_i \lambda_i(t) F_i(\lambda_i(t)) - \sum_j \mu_j(t) G_j(\mu_j(t))  \right]  \Biggr]\nonumber\\
    = & \prob({\cal C})  \expt \Biggl[ \sum_{t=1}^{T} \biggl[ \biggl(\sum_i \lambda^*_i F_i(\lambda^*_i) - \sum_j \mu^*_j G_j(\mu^*_j) \biggr) 
    -    \biggl(  \sum_i \lambda_i(t) F_i(\lambda_i(t)) - \sum_j \mu_j(t) G_j(\mu_j(t))  \biggr) \biggr] \Biggl|\Biggr. {\cal C} \Biggr]\label{equ:improved-regret-term-1-2}\\
    & + \Theta\left(T^2\epsilon^{\frac{\beta}{2}+1} + T\epsilon^{\frac{\beta}{2}-1} \log (1/\epsilon) \right).\label{equ:improved-regret-term-1-3}
\end{align}

\subsection{Bounding the Regret Caused by Rejecting Arrivals}
In this subsection, we will divide the term \eqref{equ:improved-regret-term-1-2} into two cases by comparing the queue lengths with the threshold $q^{\mathrm{th}}$, and then bound the regret induced by rejecting arrivals when the queue length exceeds the threshold $q^{\mathrm{th}}$.

Recall the definition of ${\cal H}_t$:
\[
{\cal H}_t \coloneqq \left\{ Q_{\mathrm{c},i}(t) < q^{\mathrm{th}} \mbox{ for all } i, \mbox{ and } Q_{\mathrm{s},j}(t) < q^{\mathrm{th}} \mbox{ for all } j\right\}.
\]
Following the same argument as that in Section~\ref{app:bound-regret-reject-arr}, we have
\begin{align}
    \eqref{equ:improved-regret-term-1-2} & =  \prob({\cal C})  \expt \Biggl[ \sum_{t=1}^{T} \mathbb{1}_{{\cal H}_t^{\mathrm{c}}} \biggl[ \biggl(\sum_i \lambda^*_i F_i(\lambda^*_i) - \sum_j \mu^*_j G_j(\mu^*_j) \biggr) 
    -    \biggl(  \sum_i \lambda_i(t) F_i(\lambda_i(t)) - \sum_j \mu_j(t) G_j(\mu_j(t))  \biggr) \biggr] \Biggl|\Biggr. {\cal C} \Biggr]\nonumber\\
    & + \prob({\cal C})  \expt \Biggl[ \sum_{t=1}^{T} \mathbb{1}_{{\cal H}_t} \biggl[ \biggl(\sum_i \lambda^*_i F_i(\lambda^*_i) - \sum_j \mu^*_j G_j(\mu^*_j) \biggr) 
    -    \biggl(  \sum_i \lambda_i(t) F_i(\lambda_i(t)) - \sum_j \mu_j(t) G_j(\mu_j(t))  \biggr) \biggr] \Biggl|\Biggr. {\cal C} \Biggr]\nonumber\\
    & \le \biggl[\sum_i \lambda^*_i F_i(\lambda^*_i) - \sum_j \mu^*_j G_j(\mu^*_j) + \sum_j G_j(1)\biggr]\nonumber\\
    & \qquad \sum_{t=1}^{T} \sum_i \expt \biggl[  \mathbb{1}_{\cal C} \mathbb{1} \left\{ Q_{\mathrm{c},i}(t) \ge q^{\mathrm{th}}\right\} \biggr] 
    + \sum_{t=1}^{T} \sum_j \expt \biggl[  \mathbb{1}_{\cal C} \mathbb{1} \left\{ Q_{\mathrm{s},j}(t) \ge q^{\mathrm{th}}\right\}   \biggr]\label{equ:improved-regret-term-1-4}\\
    & +  \expt \Biggl[ \sum_{t=1}^{T} \mathbb{1}_{{\cal H}_t} \biggl[ \biggl(\sum_i \lambda^*_i F_i(\lambda^*_i) - \sum_j \mu^*_j G_j(\mu^*_j) \biggr) 
    -    \biggl(  \sum_i \lambda_i(t) F_i(\lambda_i(t)) - \sum_j \mu_j(t) G_j(\mu_j(t))  \biggr) \biggr] \Biggl|\Biggr. {\cal C} \Biggr],\label{equ:improved-regret-term-1-5}
\end{align}
To bound the term \eqref{equ:improved-regret-term-1-4}, unlike the proof of Theorem~\ref{theo:1}, we do not need to consider whether $t$ is in the first outer iteration ($t\le 2MN$) or not. The reason is that we know from Assumption~\ref{assum:5} that we have a set of arrival rates that are nearly balanced at the beginning of the algorithm so that the arrival rates will remain nearly balanced all the time, which will be shown in Lemma~\ref{lemma:improved-bisection-error}. Also, we use the threshold $q^{\mathrm{th}}$ to reject arrivals all the time under the balanced pricing algorithm. In fact, we have the following lemma.
\begin{lemma}\label{lemma:improved-bound-time-large-queue}
    Let Assumption~\ref{assum:1}, Assumption~\ref{assum:4}, and Assumption~\ref{assum:5} hold. Suppose $\epsilon < \delta$ and $T$ is sufficiently large. Under the matching algorithm~\ref{alg:matching} and the balanced pricing algorithm, we have
    \begin{align*}
    &  \sum_{t=1}^{T} \sum_i \expt \Bigl[ 
    \mathbb{1}_{{\cal C}}
    \mathbb{1} \Bigl\{ Q_{\mathrm{c},i}(t) \ge q^{\mathrm{th}} \Bigr\}
    \Bigr] + \sum_{t=1}^{T} \sum_j \expt \Bigl[ 
    \mathbb{1}_{{\cal C}}
    \mathbb{1} \Bigl\{ Q_{\mathrm{s},j}(t) \ge q^{\mathrm{th}} \Bigr\}
    \Bigr] \nonumber\\
    \le & \Theta \left(
    \frac{T}{q^{\mathrm{th}}}
    + T^2\epsilon^{\frac{\beta}{2}+1} + T\epsilon^{\frac{\beta}{2}-1} \log (1/\epsilon) 
    + T \left(\frac{\eta \epsilon}{\delta} + \eta + \delta + \epsilon\right)  \right),
    \end{align*}
    where in the notation $\Theta(\cdot)$ we ignore the variables that do not depend on $T$.
\end{lemma}
\noindent Lemma~\ref{lemma:improved-bound-time-large-queue} provides upper bounds for the expected number of time slots in which the queue length is greater than or equal to $q^{\mathrm{th}}$ and the ``good'' event ${\cal C}$ holds. The result in Lemma~\ref{lemma:improved-bound-time-large-queue} improves that in Lemma~\ref{lemma:bound-time-large-queue} under the additional Assumption~\ref{assum:5} and the \textit{balanced pricing algorithm}. Proof of Lemma~\ref{lemma:improved-bound-time-large-queue} can be found in Appendix~\ref{app:proof-lemma-improved-bound-time-large-queue}.

From \eqref{equ:improved-regret-term-1-4} and Lemma~\ref{lemma:improved-bound-time-large-queue}, the regret caused by rejecting arrivals can be bounded by
\begin{align}\label{equ:improved-regret-term-1-4-final}
    \eqref{equ:improved-regret-term-1-4}
    \le &  \Theta \left(
    \frac{T}{q^{\mathrm{th}}}
    + T^2\epsilon^{\frac{\beta}{2}+1} + T\epsilon^{\frac{\beta}{2}-1} \log (1/\epsilon) 
    + T \left(\frac{\eta \epsilon}{\delta} + \eta + \delta + \epsilon\right)  \right).
\end{align}

\subsection{Bounding the Regret during Optimization}
In this subsection, we will bound the term \eqref{equ:improved-regret-term-1-5}, which are the regret during optimization. The process of bounding \eqref{equ:improved-regret-term-1-5} is similar to that in Section~\ref{app:regret-zero-order-opt}. The difference is that we do not need to divide the regret into two parts, one part in the first outer iteration and the other part after the first outer iteration. We can treat them together because we show in Lemma~\ref{lemma:improved-bisection-error} that the arrival rates are close to the target arrival rates generated in the optimization process for all time.

Recall that $\lambda_i(t),\mu_j(t)$ are the actual arrival rates at time $t$, and $\lambda'_i(t)\coloneqq F^{-1}_i(p^{+/-}_{\mathrm{c},i}(k(t),m(t)))$ and $\mu'_j(t)\coloneqq G^{-1}_j(p^{+/-}_{\mathrm{s},j}(k(t),m(t)))$, which are the arrival rates as if there is no queue control. Following the same argument in Section~\ref{app:regret-zero-order-opt}, we have
\begin{align}\label{equ:improved-regret-term-1-9}
    & \eqref{equ:improved-regret-term-1-5} = \expt \Biggl[ \sum_{t=1}^{T} \mathbb{1}_{{\cal H}_t} \biggl[ \biggl(\sum_i \lambda^*_i F_i(\lambda^*_i) - \sum_j \mu^*_j G_j(\mu^*_j) \biggr) 
    - \biggl(  \sum_i \lambda'_i(t) F_i(\lambda'_i(t)) - \sum_j \mu'_j(t) G_j(\mu'_j(t))  \biggr) \biggr] \Biggl|\Biggr. {\cal C} \Biggr].
\end{align}
By Lemma~\ref{lemma:improved-bisection-error}, we can relate $\lambda'_i(t)$ and $\mu'_j(t)$ to the target arrival rates generated in the optimization process, i.e., $\sum_{j\in {\cal E}_{\mathrm{c},i}} x_{i,j}(k(t))$ and $\sum_{i\in {\cal E}_{\mathrm{s},j}} x_{i,j}(k(t))$. Following the same argument as that in Section~\ref{app:regret-zero-order-opt}, we have
\begin{align}\label{equ:improved-regret-term-1-10}
    \eqref{equ:improved-regret-term-1-5} \le &
    \Theta\biggl( T \biggl(\frac{\eta \epsilon}{\delta} + \eta + \delta + \epsilon \biggr) \biggr)
    + \expt \Biggl[ \sum_{t=1}^{T} \mathbb{1}_{{\cal H}_t}
    \bigl[f(\boldsymbol{x}^*) - f(\boldsymbol{x}(k(t)))\bigr]
    \Biggl|\Biggr. {\cal C} \Biggr]\nonumber\\
    \le & \Theta\biggl( T \biggl(\frac{\eta \epsilon}{\delta} + \eta + \delta + \epsilon \biggr) \biggr) + 2MN  \sum_{k=1}^{\lceil T/(2MN) \rceil} \expt \Bigl[
    f(\boldsymbol{x}^*) - f(\boldsymbol{x}(k))
    \Bigl|\Bigr. {\cal C} \Bigr],
\end{align}
where we recall the function $f: {\cal D} \rightarrow \mathbb{R}$ is defined as follows:
\begin{align*}
    f(\boldsymbol{x}) \coloneqq \sum_i \biggl(\sum_{j\in {\cal E}_{\mathrm{c},i}} x_{i,j}\biggr) F_i\biggl(\sum_{j\in {\cal E}_{\mathrm{c},i}} x_{i,j}\biggr) 
    - \sum_j \biggl( \sum_{i\in {\cal E}_{\mathrm{s},j}} x_{i,j} \biggr) G_j\biggl( \sum_{i\in {\cal E}_{\mathrm{s},j}} x_{i,j} \biggr).
\end{align*}
Note that Lemma~\ref{lemma:optimization} also holds for the \textit{balanced pricing algorithm}. We omit the proof because it is the almost the same except that we should use Lemma~\ref{lemma:improved-bisection-error} instead of Lemma~\ref{lemma:bisection-error}, but both lemmas have the same result.
By \eqref{equ:improved-regret-term-1-10}, Lemma~\ref{lemma:optimization}, and the definitions of $M$ and $N$, we have
\begin{align}\label{equ:improved-regret-term-1-5-final}
    \eqref{equ:improved-regret-term-1-5}
    \le & \Theta\biggl( T \biggl(\frac{\eta \epsilon}{\delta} + \eta + \delta + \epsilon \biggr) \biggr) + \Theta \left( T \left(\eta + \frac{\eta \epsilon^2}{\delta^2} + \frac{\epsilon}{\delta} + \delta \right) + \left( \frac{1}{\eta \epsilon^2 } + \frac{\eta}{\delta^2}
    + \frac{1}{\epsilon \delta} \right) \log^2(1/\epsilon) \right)\nonumber\\
    \le & \Theta \left( T\left(\frac{\epsilon}{\delta} + \eta + \delta\right) 
        +  \frac{\log^2 (1/\epsilon)}{\eta \epsilon^2} \right)
\end{align}
where the last inequality is by $\epsilon < \delta$.

Combining \eqref{equ:improved-regret-term-1-4}, \eqref{equ:improved-regret-term-1-5}, \eqref{equ:improved-regret-term-1-4-final}, and \eqref{equ:improved-regret-term-1-5-final}, we have
\begin{align}\label{equ:improved-regret-term-1-2-final}
\eqref{equ:improved-regret-term-1-2} \le  \Theta \left(
    \frac{T}{q^{\mathrm{th}}}
    + T^2\epsilon^{\frac{\beta}{2}+1} + T\epsilon^{\frac{\beta}{2}-1} \log (1/\epsilon) 
    +  T\left(\frac{\epsilon}{\delta} + \eta + \delta\right) 
    +  \frac{\log^2 (1/\epsilon)}{\eta \epsilon^2} \right).
\end{align}
Combining \eqref{equ:improved-regret-term-1-2}, \eqref{equ:improved-regret-term-1-3}, and \eqref{equ:improved-regret-term-1-2-final}, we have
\begin{align*}
    \expt[R(T)] \le \Theta \left(
    \frac{T}{q^{\mathrm{th}}}
    + T^2\epsilon^{\frac{\beta}{2}+1} + T\epsilon^{\frac{\beta}{2}-1} \log (1/\epsilon) 
    +  T\left(\frac{\epsilon}{\delta} + \eta + \delta\right) 
    +  \frac{\log^2 (1/\epsilon)}{\eta \epsilon^2} \right).
\end{align*}
Therefore, the regret bound \eqref{equ:regret-bound-2} in Theorem~\ref{theo:2} is proved.

\subsection{Bounding the Maximum and Average Queue Length}

Under \textit{the balanced pricing algorithm}, since for each queue we set the arrival rate to be zero when the queue length is greater than or equal to $q^{\mathrm{th}}$ for all iterations and all the queues are initially empty, the queue length at any time must be less than or equal to $q^{\mathrm{th}}$. Hence, the maximum queue-length bound \eqref{equ:queue-length-bound-2} and the average queue-length bound \eqref{equ:avg-queue-length-bound-2} in Theorem~\ref{theo:2} are proved.

\section{Proof of Lemma~\ref{lemma:empty-queues}}
\label{app:proof-lemma-empty-queues}

Prove by induction.

First at time slot $t=1$, all queues are empty, so the claimed property holds. 

Assume that at time slot $t$, the claimed property holds, i.e., at time slot $t$, for any demand-side queue $i$, if $Q_{\mathrm{c},i}(t)>0$, then $\sum_{j: (i,j)\in{\cal E}} Q_{\mathrm{s}, j}(t) = 0$; for any supply-side queue $j$, if $Q_{\mathrm{s},j}(t)>0$, then $\sum_{i: (i,j)\in{\cal E}} Q_{\mathrm{c}, i}(t) = 0$. We need to show that this also holds for $t+1$.

We will prove this by contradiction. 
Suppose this does not hold for $t+1$. Then there exists a demand-side queue $i$, $Q_{\mathrm{c},i}(t+1)>0$ and $\sum_{j: (i,j)\in{\cal E}} Q_{\mathrm{s}, j}(t+1) > 0$, or there exists or a supply side queue $j$, $Q_{\mathrm{s},j}(t+1)>0$ and $\sum_{i: (i,j)\in{\cal E}} Q_{\mathrm{c}, i}(t+1) > 0$.

First consider the case where there exists a demand-side queue $i$, $Q_{\mathrm{c},i}(t+1)>0$ and $\sum_{j: (i,j)\in{\cal E}} Q_{\mathrm{s}, j}(t+1) > 0$. There are two possibilities:
\begin{enumerate}
    \item $Q_{\mathrm{c},i}(t)>0$ and $\sum_{j: (i,j)\in{\cal E}} Q_{\mathrm{s}, j}(t) = 0$\\
    In this case, there must be arrivals in at least one of the queues $\{j: (i,j)\in{\cal E}\}$ at time $t$. However, under Algorithm~\ref{alg:matching}, these arrivals must be all matched at time slot $t$ because $Q_{\mathrm{c},i}(t)>0$ and $Q_{\mathrm{c},i}(t+1)>0$. This contradicts the fact that $\sum_{j: (i,j)\in{\cal E}} Q_{\mathrm{s}, j}(t+1) > 0$.
    \item $Q_{\mathrm{c},i}(t)=0$\\
    In this case, there must be an arrival in queue $i$. However, under Algorithm~\ref{alg:matching}, the arrival must be matched at time slot $t$ because $\sum_{j: (i,j)\in{\cal E}} Q_{\mathrm{s}, j}(t+1) > 0$. This contradicts that fact that $Q_{\mathrm{c},i}(t+1)>0$.
\end{enumerate}

Consider the case where there exists a supply side queue $j$, $Q_{\mathrm{s},j}(t+1)>0$ and $\sum_{i: (i,j)\in{\cal E}} Q_{\mathrm{c}, i}(t+1) > 0$. There are two possibilities:
\begin{enumerate}
    \item $Q_{\mathrm{s},j}(t)>0$ and $\sum_{i: (i,j)\in{\cal E}} Q_{\mathrm{c}, i}(t) = 0$\\
    In this case, there must be arrivals in at least one of the queues $\{i: (i,j)\in{\cal E}\}$ at time $t$. However, under Algorithm~\ref{alg:matching}, these arrivals must be all matched at time slot $t$ because $Q_{\mathrm{s},j}(t)>0$ and $Q_{\mathrm{s},j}(t+1)>0$. This contradicts the fact that $\sum_{i: (i,j)\in{\cal E}} Q_{\mathrm{c}, i}(t+1) > 0$.
    \item $Q_{\mathrm{s},j}(t)=0$\\
    In this case, there must be an arrival in queue $j$. However, under Algorithm~\ref{alg:matching}, the arrival must be matched at time slot $t$ because $\sum_{i: (i,j)\in{\cal E}} Q_{\mathrm{c}, i}(t+1) > 0$. This contradicts that fact that $Q_{\mathrm{s},j}(t+1)>0$.
\end{enumerate}

Therefore, we prove by contradiction that for any demand-side queue $i$, if $Q_{\mathrm{c},i}(t+1)>0$, then $\sum_{j: (i,j)\in{\cal E}} Q_{\mathrm{s}, j}(t+1) = 0$; for any supply-side queue $j$, if $Q_{\mathrm{s},j}(t+1)>0$, then $\sum_{i: (i,j)\in{\cal E}} Q_{\mathrm{c}, i}(t+1) = 0$. By induction, Lemma~\ref{lemma:empty-queues} is proved.

\section{Proof of Lemma~\ref{lemma:shrunk-set}}
\label{app:proof-lemma-shrunk-set}

From the definition of ${\cal D}$ in~\eqref{equ:def-D}, we have
\begin{align*}
    & {\cal D} - \boldsymbol{x}_{\mathrm{ctr}}\nonumber\\
    = &  \left\{ \boldsymbol{x} - \boldsymbol{x}_{\mathrm{ctr}} \left| \boldsymbol{x} \in {\cal D}\right.\right\}\nonumber\\
    = & \Biggl\{  \tilde{\boldsymbol{x}} \Biggl|
    \text{ for all } i,j, ~
    \tilde{x}_{i,j} \ge -\frac{a_{\min}+1}{2N_{i,j}},
    \sum_{j'\in {\cal E}_{\mathrm{c},i}} \left[\tilde{x}_{i,j'} + \frac{a_{\min}+1}{2N_{i,j'}} \right] \in [a_{\min}, 1],\nonumber\\
    & \quad \sum_{i'\in {\cal E}_{\mathrm{s},j}} \left[ \tilde{x}_{i',j} + \frac{a_{\min}+1}{2N_{i',j}}\right] \in [a_{\min}, 1] 
    \Biggr. \Biggr\}\nonumber\\
    = & \Biggl\{  \tilde{\boldsymbol{x}} \Biggl| 
    \text{ for all } i,j, ~
    \tilde{x}_{i,j} \ge -\frac{a_{\min}+1}{2N_{i,j}}, \nonumber\\
    & \quad \sum_{j'\in {\cal E}_{\mathrm{c},i}} \tilde{x}_{i,j'}   \in \biggl[ - \biggl(\sum_{j'\in {\cal E}_{\mathrm{c},i}} \frac{a_{\min}+1}{2N_{i,j'}} - a_{\min}\biggr), 
    1  - \sum_{j'\in {\cal E}_{\mathrm{c},i}} \frac{a_{\min}+1}{2N_{i,j'}}\biggr],\nonumber\\
    & \quad \sum_{i'\in {\cal E}_{\mathrm{s},j}} \tilde{x}_{i',j} \in \biggl[  - \biggl(\sum_{i'\in {\cal E}_{\mathrm{s},j}} \frac{a_{\min}+1}{2N_{i',j}} -  a_{\min} \biggr), 
    1 - \sum_{i'\in {\cal E}_{\mathrm{s},j}} \frac{a_{\min}+1}{2N_{i',j}} \biggr]
    \Biggr. \Biggr\}.
\end{align*}
Shrunk by $(1-\frac{\delta}{r})$, the set  ${\cal D} - \boldsymbol{x}_{\mathrm{ctr}}$  becomes
\begin{align*}
    & (1-\frac{\delta}{r})({\cal D} - \boldsymbol{x}_{\mathrm{ctr}})\nonumber\\
    = & \Biggl\{  \tilde{\boldsymbol{x}} \Biggl| 
    \text{ for all } i,j, ~
    \tilde{x}_{i,j} \ge - (1-\frac{\delta}{r}) \frac{a_{\min}+1}{2N_{i,j}} ,\nonumber\\
    & \quad \sum_{j'\in {\cal E}_{\mathrm{c},i}} \tilde{x}_{i,j'}   \in \biggl[ - (1-\frac{\delta}{r})\biggl(\sum_{j'\in {\cal E}_{\mathrm{c},i}} \frac{a_{\min}+1}{2N_{i,j'}} - a_{\min}\biggr), 
    (1-\frac{\delta}{r}) \biggl(1  - \sum_{j'\in {\cal E}_{\mathrm{c},i}} \frac{a_{\min}+1}{2N_{i,j'}}\biggr) \biggr],\nonumber\\
    & \quad \sum_{i'\in {\cal E}_{\mathrm{s},j}} \tilde{x}_{i',j} \in \biggl[  - (1-\frac{\delta}{r})\biggl(\sum_{i'\in {\cal E}_{\mathrm{s},j}} \frac{a_{\min}+1}{2N_{i',j}} -  a_{\min} \biggr) , 
    (1-\frac{\delta}{r}) \biggl(1 - \sum_{i'\in {\cal E}_{\mathrm{s},j}} \frac{a_{\min}+1}{2N_{i',j}} \biggr)\biggr]
    \Biggr. \Biggr\}.
\end{align*}
Next, shifted by $\boldsymbol{x}_{\mathrm{ctr}}$, the set becomes
\begin{align*}
    & (1-\frac{\delta}{r})({\cal D} - \boldsymbol{x}_{\mathrm{ctr}}) + \boldsymbol{x}_{\mathrm{ctr}}\nonumber\\
    = & \left\{\tilde{\boldsymbol{x}} + \boldsymbol{x}_{\mathrm{ctr}} | \tilde{\boldsymbol{x}} \in (1-\frac{\delta}{r})({\cal D} - \boldsymbol{x}_{\mathrm{ctr}})\right\}\nonumber\\
    = & \Biggl\{  \boldsymbol{x}'\Biggl|
    \text{ for all } i,j, ~
    x_{i,j}' - \frac{a_{\min}+1}{2N_{i,j}}  \ge - (1-\frac{\delta}{r}) \frac{a_{\min}+1}{2N_{i,j}} \nonumber\\
    & \sum_{j'\in {\cal E}_{\mathrm{c},i}} \biggl( x'_{i,j'} - \frac{a_{\min}+1}{2N_{i,j'}} \biggr)  \in \biggl[ - (1-\frac{\delta}{r})\biggl(\sum_{j'\in {\cal E}_{\mathrm{c},i}} \frac{a_{\min}+1}{2N_{i,j'}} - a_{\min}\biggr), 
    (1-\frac{\delta}{r}) \biggl(1  - \sum_{j'\in {\cal E}_{\mathrm{c},i}} \frac{a_{\min}+1}{2N_{i,j'}}\biggr) \biggr],\nonumber\\
    & \sum_{i'\in {\cal E}_{\mathrm{s},j}} \biggl( x'_{i',j} - \frac{a_{\min}+1}{2N_{i',j}} \biggl) \in \biggl[  - (1-\frac{\delta}{r})\biggl(\sum_{i'\in {\cal E}_{\mathrm{s},j}} \frac{a_{\min}+1}{2N_{i',j}} -  a_{\min} \biggr) , 
    (1-\frac{\delta}{r}) \biggl(1 - \sum_{i'\in {\cal E}_{\mathrm{s},j}} \frac{a_{\min}+1}{2N_{i',j}} \biggr)\biggr]
    \Biggr. \Biggr\}\nonumber\\
    = & {\cal D}',
\end{align*}
where the last equality is by the definition of ${\cal D}'$ in~\eqref{equ:def-D'}.

Let $\boldsymbol{x}\in {\cal D}'$. We need to show that $\boldsymbol{x} + \delta \boldsymbol{u} \in {\cal D}$. Note that $\boldsymbol{x}\in {\cal D}' = (1-\frac{\delta}{r})({\cal D} - \boldsymbol{x}_{\mathrm{ctr}})+ \boldsymbol{x}_{\mathrm{ctr}}$. Hence, there exists $\boldsymbol{y}\in {\cal D}$ such that
\begin{align*}
    \boldsymbol{x} = (1-\frac{\delta}{r}) (\boldsymbol{y} - \boldsymbol{x}_{\mathrm{ctr}})+ \boldsymbol{x}_{\mathrm{ctr}}.
\end{align*}
Then we have
\begin{align}\label{equ:convex-combination}
    \boldsymbol{x} + \delta \boldsymbol{u} = & (1-\frac{\delta}{r}) (\boldsymbol{y} - \boldsymbol{x}_{\mathrm{ctr}})+ \boldsymbol{x}_{\mathrm{ctr}} + \delta \boldsymbol{u}\nonumber\\
    = & (1-\frac{\delta}{r}) (\boldsymbol{y} - \boldsymbol{x}_{\mathrm{ctr}})+ (1-\frac{\delta}{r})\boldsymbol{x}_{\mathrm{ctr}} + \frac{\delta}{r} \boldsymbol{x}_{\mathrm{ctr}}  + \delta \boldsymbol{u} \nonumber\\
    = & (1-\frac{\delta}{r}) \boldsymbol{y} + \frac{\delta}{r} (\boldsymbol{x}_{\mathrm{ctr}} + r \boldsymbol{u}).
\end{align}
We next show that $\boldsymbol{x}_{\mathrm{ctr}} + r \boldsymbol{u}\in {\cal D}$. Note that $\boldsymbol{u} \in \mathbb{B}$ and hence $u_{i,j} \in [-1,1]$ for all $i,j$. For any $i,j$, we have 
\begin{align*}
    x_{\mathrm{ctr}, i,j} + r u_{i,j} \ge x_{\mathrm{ctr}, i,j} - r.
\end{align*}
Recall the definition of $r$ in~\eqref{equ:def-r}. We know $r\le \frac{a_{\min}+1}{2N_{i,j}}$. Hence, we have
\begin{align}\label{equ:proof-lemma-set-D-1}
    x_{\mathrm{ctr}, i,j} + r u_{i,j} \ge x_{\mathrm{ctr}, i,j} - r \ge x_{\mathrm{ctr}, i,j} - \frac{a_{\min}+1}{2N_{i,j}} = 0.
\end{align}
We know from~\eqref{equ:def-r} that $r\le  \frac{1}{|{\cal E}_{\mathrm{c}, i}|} ( \sum_{j'\in {\cal E}_{\mathrm{c},i}} \frac{a_{\min}+1}{2N_{i,j'}} - a_{\min} )$. Hence, we have
\begin{align}\label{equ:proof-lemma-set-D-2}
    \sum_{j\in{\cal E}_{\mathrm{c},i}} \left(x_{\mathrm{ctr}, i,j} + r u_{i,j}\right)
    \ge \sum_{j\in{\cal E}_{\mathrm{c},i}} \left(x_{\mathrm{ctr}, i,j} - r\right) 
    \ge \sum_{j\in{\cal E}_{\mathrm{c},i}}  \left[ \frac{a_{\min}+1}{2N_{i,j}} - \frac{1}{|{\cal E}_{\mathrm{c}, i}|} \biggl( \sum_{j'\in {\cal E}_{\mathrm{c},i}} \frac{a_{\min}+1}{2N_{i,j'}} - a_{\min} \biggr)\right]  = a_{\min}.
\end{align}
We know from~\eqref{equ:def-r} that $r\le  \frac{1}{|{\cal E}_{\mathrm{s}, j}|} ( \sum_{i'\in {\cal E}_{\mathrm{s},j}} \frac{a_{\min}+1}{2N_{i',j}} - a_{\min} ) $. Hence, we have
\begin{align}\label{equ:proof-lemma-set-D-3}
    \sum_{i\in{\cal E}_{\mathrm{s},j}} \left(x_{\mathrm{ctr}, i,j} + r u_{i,j}\right)
    \ge \sum_{i\in{\cal E}_{\mathrm{s},j}} \left(x_{\mathrm{ctr}, i,j} - r\right) 
    \ge \sum_{i\in{\cal E}_{\mathrm{s},j}}  \left[ \frac{a_{\min}+1}{2N_{i,j}} - \frac{1}{|{\cal E}_{\mathrm{s}, j}|} \biggl( \sum_{i'\in {\cal E}_{\mathrm{s},j}} \frac{a_{\min}+1}{2N_{i',j}} - a_{\min} \biggr)\right]  = a_{\min}.
\end{align}
Notice that $x_{\mathrm{ctr}, i,j} + r u_{i,j} \le x_{\mathrm{ctr}, i,j} + r$ for all $i,j$. From~\eqref{equ:def-r} we have $r\le \frac{1}{|{\cal E}_{\mathrm{c}, i}|}( 1 -  \sum_{j'\in {\cal E}_{\mathrm{c},i}} \frac{a_{\min}+1}{2N_{i,j'}} )$. Hence we have
\begin{align}\label{equ:proof-lemma-set-D-4}
    \sum_{j\in{\cal E}_{\mathrm{c},i}} \left(x_{\mathrm{ctr}, i,j} + r u_{i,j}\right)
    \le \sum_{j\in{\cal E}_{\mathrm{c},i}} \left(x_{\mathrm{ctr}, i,j} + r\right) 
    \le \sum_{j\in{\cal E}_{\mathrm{c},i}}  \left[ \frac{a_{\min}+1}{2N_{i,j}} + \frac{1}{|{\cal E}_{\mathrm{c}, i}|}\biggl( 1 -  \sum_{j'\in {\cal E}_{\mathrm{c},i}} \frac{a_{\min}+1}{2N_{i,j'}} \biggr)\right]  = 1.
\end{align}
From~\eqref{equ:def-r} we have $r\le \frac{1}{|{\cal E}_{\mathrm{s}, j}|}( 1 -  \sum_{i'\in {\cal E}_{\mathrm{s},j}} \frac{a_{\min}+1}{2N_{i',j}} )$. Hence, we have
\begin{align}\label{equ:proof-lemma-set-D-5}
    \sum_{i\in{\cal E}_{\mathrm{s},j}} \left(x_{\mathrm{ctr}, i,j} + r u_{i,j}\right)
    \le \sum_{i\in{\cal E}_{\mathrm{s},j}} \left(x_{\mathrm{ctr}, i,j} + r\right) 
    \le \sum_{i\in{\cal E}_{\mathrm{s},j}}  \left[ \frac{a_{\min}+1}{2N_{i,j}} + \frac{1}{|{\cal E}_{\mathrm{s}, j}|}\biggl( 1 -  \sum_{i'\in {\cal E}_{\mathrm{s},j}} \frac{a_{\min}+1}{2N_{i',j}} \biggr)\right]  = 1.
\end{align}
From~\eqref{equ:proof-lemma-set-D-1}-\eqref{equ:proof-lemma-set-D-5} and the definition of the set ${\cal D}$ in~\eqref{equ:def-D}, we have
\begin{align*}
    \boldsymbol{x}_{\mathrm{ctr}} + r \boldsymbol{u} \in {\cal D}.
\end{align*}
Note that the constraints in ${\cal D}$ are affine and hence ${\cal D}$ is convex. Recall that from~\eqref{equ:convex-combination} we have $\boldsymbol{x} + \delta \boldsymbol{u} = (1-\frac{\delta}{r}) \boldsymbol{y} + \frac{\delta}{r} (\boldsymbol{x}_{\mathrm{ctr}} + r \boldsymbol{u})$. Since $\boldsymbol{y}\in {\cal D}$ and $\boldsymbol{x}_{\mathrm{ctr}} + r \boldsymbol{u} \in {\cal D}$, by convexity of ${\cal D}$, we have $\boldsymbol{x} + \delta \boldsymbol{u} \in {\cal D}$.

\section{Proof of Lemma~\ref{lemma:concentration}}
\label{app:proof-lemma-concentration}

Let 
\begin{align*}
    {\cal C}^{+}_{\mathrm{c},i}(k,m) \coloneqq \Biggl\{ \text{for all $\lambda_i(t^{+}_{\mathrm{c}, i}(k,m,n)) \in [0,1]$, } 
    \biggl| \frac{1}{N} \sum_{n=1}^{N} A_{\mathrm{c},i} (t^{+}_{\mathrm{c}, i}(k,m,n)) - \lambda_i(t^{+}_{\mathrm{c}, i}(k,m,n)) \biggr| < \epsilon
    \Biggr\}.
\end{align*}
Define similarly ${\cal C}^{-}_{\mathrm{c},i}(k,m)$, ${\cal C}^{+}_{\mathrm{s},j}(k,m)$, and ${\cal C}^{-}_{\mathrm{s},j}(k,m)$. Then we have
\begin{align}\label{equ:union-relation}
    & {\cal C}^{\mathrm{c}} = \nonumber\\
    &\bigcup_{k,m} \left(\left(\bigcup_{i=1}^{I} \left({\cal C}^{+}_{\mathrm{c},i}(k,m)\right)^{\mathrm{c}} \right)\cup \left(\bigcup_{i=1}^{I} \left({\cal C}^{-}_{\mathrm{c},i}(k,m)\right)^{\mathrm{c}}\right) 
    \cup \left(\bigcup_{j=1}^{J} \left({\cal C}^{+}_{\mathrm{s},j}(k,m)\right)^{\mathrm{c}} \right)\cup \left(\bigcup_{j=1}^{J} \left({\cal C}^{-}_{\mathrm{s},j}(k,m)\right)^{\mathrm{c}}\right)\right).
\end{align}
Note that all the arrivals follow independent Bernoulli distributions. Hence, we can equivalently generate the Bernoulli random variables by generating i.i.d. random variables $U^+_{\mathrm{c},i}(k,m,n)$, $U^-_{\mathrm{c},i}(k,m,n)$, $U^+_{\mathrm{s},j}(k,m,n)$, $U^-_{\mathrm{s},j}(k,m,n)$ uniformly distributed in $[0,1]$ at the beginning of the algorithm and
\begin{align*}
    A_{\mathrm{c},i} (t^{+}_{\mathrm{c}, i}(k,m,n)) = \mathbb{1} \left\{U^+_{\mathrm{c},i}(k,m,n) \le \lambda_i(t^{+}_{\mathrm{c}, i}(k,m,n))\right\},\nonumber\\
    A_{\mathrm{c},i} (t^{-}_{\mathrm{c}, i}(k,m,n)) = \mathbb{1} \left\{U^-_{\mathrm{c},i}(k,m,n) \le \lambda_i(t^{-}_{\mathrm{c}, i}(k,m,n))\right\},\nonumber\\
    A_{\mathrm{s},j} (t^{+}_{\mathrm{s}, j}(k,m,n)) = \mathbb{1} \left\{U^+_{\mathrm{s},j}(k,m,n) \le \mu_j(t^{+}_{\mathrm{s}, j}(k,m,n))\right\},\nonumber\\
    A_{\mathrm{s},j} (t^{-}_{\mathrm{s}, j}(k,m,n)) = \mathbb{1} \left\{U^-_{\mathrm{s},j}(k,m,n) \le \mu_j(t^{-}_{\mathrm{s}, j}(k,m,n))\right\}.
\end{align*}
Then the events ${\cal C}^{+}_{\mathrm{c},i}(k,m), {\cal C}^{-}_{\mathrm{c},i}(k,m)$, ${\cal C}^{+}_{\mathrm{s},j}(k,m)$, ${\cal C}^{-}_{\mathrm{s},j}(k,m)$ can be equivalently written as
\begin{align*}
    {\cal C}^{+}_{\mathrm{c},i}(k,m) = \Biggl\{ \text{for all $\lambda \in [0,1]$, } 
    \biggl| \frac{1}{N} \sum_{n=1}^{N} \mathbb{1} \left\{U^+_{\mathrm{c},i}(k,m,n) \le \lambda \right\} - \lambda \biggr| < \epsilon
    \Biggr\},\nonumber\\
    {\cal C}^{-}_{\mathrm{c},i}(k,m) = \Biggl\{ \text{for all $\lambda \in [0,1]$, } 
    \biggl| \frac{1}{N} \sum_{n=1}^{N} \mathbb{1} \left\{U^-_{\mathrm{c},i}(k,m,n) \le \lambda \right\} - \lambda \biggr| < \epsilon
    \Biggr\},\nonumber\\
    {\cal C}^{+}_{\mathrm{s},j}(k,m) = \Biggl\{ \text{for all $\mu \in [0,1]$, } 
    \biggl| \frac{1}{N} \sum_{n=1}^{N} \mathbb{1} \left\{U^+_{\mathrm{s},j}(k,m,n) \le \mu \right\} - \mu \biggr| < \epsilon
    \Biggr\},\nonumber\\
    {\cal C}^{-}_{\mathrm{s},j}(k,m) = \Biggl\{ \text{for all $\mu \in [0,1]$, } 
    \biggl| \frac{1}{N} \sum_{n=1}^{N} \mathbb{1} \left\{U^-_{\mathrm{s},j}(k,m,n) \le \mu \right\} - \mu \biggr| < \epsilon
    \Biggr\}.
\end{align*}

First consider the event ${\cal C}^{+}_{\mathrm{c},i}(k,m)$. We know that
\begin{align*}
    \left({\cal C}^{+}_{\mathrm{c},i}(k,m)\right) ^ {\mathrm{c}} = \Biggl\{ \text{there exists $\lambda \in [0,1]$ such that } 
    \biggl| \frac{1}{N} \sum_{n=1}^{N} \mathbb{1} \left\{U^+_{\mathrm{c},i}(k,m,n) \le \lambda \right\} - \lambda \biggr| \ge \epsilon
    \Biggr\}.
\end{align*}
We use an discretization idea that is similar to the covering argument for linear bandits~\cite{lattimore2020bandit}. Define a set ${\cal S}_{[0,1]}$, a discretized version of $[0,1]$, as follows:
\begin{align*}
    {\cal S}_{[0,1]} = \left\{0, \frac{\epsilon}{2}, \epsilon,  \frac{3\epsilon}{2},\ldots, \left\lfloor\frac{2}{\epsilon}\right\rfloor \frac{\epsilon}{2}, 1 \right \}.
\end{align*}
Then for any $\lambda\in[0,1]$, there exists two points $\lambda_{\mathrm{l}}, \lambda_{\mathrm{u}} \in {\cal S}_{[0,1]}$ such that $\lambda_{\mathrm{l}} \le \lambda \le \lambda_{\mathrm{u}}$, $|\lambda - \lambda_{\mathrm{l}}| \le \frac{\epsilon}{2}$, and $|\lambda - \lambda_{\mathrm{u}}| \le \frac{\epsilon}{2}$. Hence, we have
\begin{align}\label{equ:conc-upper-lower}
    \frac{1}{N} \sum_{n=1}^{N} \mathbb{1} \left\{U^+_{\mathrm{c},i}(k,m,n) \le \lambda_{\mathrm{l}} \right\} 
    \le \frac{1}{N} \sum_{n=1}^{N} \mathbb{1} \left\{U^+_{\mathrm{c},i}(k,m,n) \le \lambda \right\} \le \frac{1}{N} \sum_{n=1}^{N} \mathbb{1} \left\{U^+_{\mathrm{c},i}(k,m,n) \le \lambda_{\mathrm{u}} \right\}.
\end{align}
Then from~\eqref{equ:conc-upper-lower} we have
\begin{align*}
    \left({\cal C}^{+}_{\mathrm{c},i}(k,m)\right) ^ {\mathrm{c}} = & \Biggl\{ \text{there exists $\lambda \in [0,1]$ such that }
    \frac{1}{N} \sum_{n=1}^{N} \mathbb{1} \left\{U^+_{\mathrm{c},i}(k,m,n) \le \lambda \right\} \ge \lambda + \epsilon \nonumber\\
    & \quad \text{or } \frac{1}{N} \sum_{n=1}^{N} \mathbb{1} \left\{U^+_{\mathrm{c},i}(k,m,n) \le \lambda \right\} \le \lambda - \epsilon
    \Biggr\}\nonumber\\
    \subseteq & \Biggl\{ \text{there exists $\lambda \in [0,1]$ such that }
    \frac{1}{N} \sum_{n=1}^{N} \mathbb{1} \left\{U^+_{\mathrm{c},i}(k,m,n) \le \lambda_{\mathrm{u}} \right\} \ge \lambda + \epsilon \nonumber\\
    & \quad \text{or } \frac{1}{N} \sum_{n=1}^{N} \mathbb{1} \left\{U^+_{\mathrm{c},i}(k,m,n) \le \lambda_{\mathrm{l}} \right\} \le \lambda - \epsilon
    \Biggr\}.
\end{align*}
Since $\lambda_{\mathrm{l}} \le \lambda \le \lambda_{\mathrm{u}}$, $|\lambda - \lambda_{\mathrm{l}}| \le \frac{\epsilon}{2}$, and $|\lambda - \lambda_{\mathrm{u}}| \le \frac{\epsilon}{2}$, we further have
\begin{align*}
    \left({\cal C}^{+}_{\mathrm{c},i}(k,m)\right) ^ {\mathrm{c}} \subseteq & \Biggl\{ \text{there exists $\lambda \in [0,1]$ such that }
    \frac{1}{N} \sum_{n=1}^{N} \mathbb{1} \left\{U^+_{\mathrm{c},i}(k,m,n) \le \lambda_{\mathrm{u}} \right\} \ge \lambda_{\mathrm{u}} + \frac{\epsilon}{2} \nonumber\\
    & \quad \text{or } \frac{1}{N} \sum_{n=1}^{N} \mathbb{1} \left\{U^+_{\mathrm{c},i}(k,m,n) \le \lambda_{\mathrm{l}} \right\} \le \lambda_{\mathrm{l}} - \frac{\epsilon}{2}
    \Biggr\}\nonumber\\
    = & \Biggl\{ \text{there exists $\lambda_{\mathrm{u}} \in {\cal S}_{[0,1]}$ such that }
    \frac{1}{N} \sum_{n=1}^{N} \mathbb{1} \left\{U^+_{\mathrm{c},i}(k,m,n) \le \lambda_{\mathrm{u}} \right\} \ge \lambda_{\mathrm{u}} + \frac{\epsilon}{2} \Biggr\} \nonumber\\
    & \cup \Biggl\{ \text{there exists $\lambda_{\mathrm{l}} \in {\cal S}_{[0,1]}$ such that } \frac{1}{N} \sum_{n=1}^{N} \mathbb{1} \left\{U^+_{\mathrm{c},i}(k,m,n) \le \lambda_{\mathrm{l}} \right\} \le \lambda_{\mathrm{l}} - \frac{\epsilon}{2}
    \Biggr\}\nonumber\\
    = & \Biggl\{ \text{there exists $\lambda \in {\cal S}_{[0,1]}$ such that }
    \biggl| \frac{1}{N} \sum_{n=1}^{N} \mathbb{1} \left\{U^+_{\mathrm{c},i}(k,m,n) \le \lambda \right\} - \lambda \biggr| \ge \frac{\epsilon}{2} \Biggr\}.
\end{align*}
Hence, we have
\begin{align}\label{equ:conc-union}
    \prob \left( \left({\cal C}^{+}_{\mathrm{c},i}(k,m)\right) ^ {\mathrm{c}} \right) \le & \prob \left(\text{there exists $\lambda \in {\cal S}_{[0,1]}$ such that }
    \biggl| \frac{1}{N} \sum_{n=1}^{N} \mathbb{1} \left\{U^+_{\mathrm{c},i}(k,m,n) \le \lambda \right\} - \lambda \biggr| \ge \frac{\epsilon}{2}  \right) \nonumber\\
    \le & \sum_{\lambda\in {\cal S}_{[0,1]}} \prob\left( \biggl| \frac{1}{N} \sum_{n=1}^{N} \mathbb{1} \left\{U^+_{\mathrm{c},i}(k,m,n) \le \lambda \right\} - \lambda \biggr| \ge \frac{\epsilon}{2} \right),
\end{align}
where the last inequality is by union bound. By Hoeffding's inequality, we have
\begin{align}\label{equ:conc-hoeffding}
    \prob \left( \biggl| \frac{1}{N} \sum_{n=1}^{N} \mathbb{1} \left\{U^+_{\mathrm{c},i}(k,m,n) \le \lambda \right\} - \lambda \biggr| \ge \frac{\epsilon}{2} \right) 
    \le 2 \exp \left( -\frac{N\epsilon^2}{2} \right).
\end{align}
Hence, from~\eqref{equ:conc-union} and \eqref{equ:conc-hoeffding}, we have
\begin{align*}
    \prob \left( \left({\cal C}^{+}_{\mathrm{c},i}(k,m)\right) ^ {\mathrm{c}} \right) 
    \le 2 \left| {\cal S}_{[0,1]} \right| \exp \left( -\frac{N\epsilon^2}{2} \right) \le \left(\frac{4}{\epsilon} + 4\right) \exp \left( - \frac{N \epsilon^2 }{2} \right),
\end{align*}
where the last inequality is due to the fact that $\left| {\cal S}_{[0,1]} \right| \le \frac{2}{\epsilon} + 2$. Since $N=\left\lceil\frac{\beta\ln (1/\epsilon)}{\epsilon^2}\right\rceil$, we have
\begin{align}
    \prob \left( \left({\cal C}^{+}_{\mathrm{c},i}(k,m)\right) ^ {\mathrm{c}} \right)  \le 4 \epsilon^{\frac{\beta}{2}-1} + 4 \epsilon^{\frac{\beta}{2}} \le 8 \epsilon^{\frac{\beta}{2}-1}.\label{equ:conc-sub-1}
\end{align}
Similarly, we have
\begin{align}
    \prob \left( \left({\cal C}^{-}_{\mathrm{c},i}(k,m)\right) ^ {\mathrm{c}} \right)  \le &  8 \epsilon^{\frac{\beta}{2}-1},\label{equ:conc-sub-2}\\
    \prob \left( \left({\cal C}^{+}_{\mathrm{s},j}(k,m)\right) ^ {\mathrm{c}} \right)  \le &  8 \epsilon^{\frac{\beta}{2}-1},\label{equ:conc-sub-3}\\
    \prob \left( \left({\cal C}^{-}_{\mathrm{s},j}(k,m)\right) ^ {\mathrm{c}} \right)  \le &  8 \epsilon^{\frac{\beta}{2}-1},\label{equ:conc-sub-4}
\end{align}
for all $k,m,i,j$. Therefore, by \eqref{equ:union-relation}, \eqref{equ:conc-sub-1}-\eqref{equ:conc-sub-4}, and the union bound, we have
\begin{align}\label{equ:bound-conc}
    \prob ({\cal C}^{\mathrm{c}}) 
    \le & \left\lceil \frac{T}{2MN}\right\rceil   \times M  \times 2 (I + J) \times 8 \epsilon^{\frac{\beta}{2}-1} \le \frac{8T(I+J) \epsilon^{\frac{\beta}{2}-1}  }{N} 
    + 16 M (I + J)  \epsilon^{\frac{\beta}{2}-1}\nonumber\\
    \le & \frac{8}{\beta} T(I+J)\epsilon^{\frac{\beta}{2}+1} + 32 (I+J) \epsilon^{\frac{\beta}{2}-1} \log_2 (1/\epsilon) = \Theta\left(T\epsilon^{\frac{\beta}{2}+1} + \epsilon^{\frac{\beta}{2}-1} \log (1/\epsilon) \right),
\end{align}
where the last inequality holds since $N=\left\lceil\frac{\beta\ln (1/\epsilon)}{\epsilon^2}\right\rceil \ge \frac{\beta}{\epsilon^2}$ by $\epsilon< 1/e$, and $M=\lceil \log_2 (1/\epsilon) \rceil$. Lemma~\ref{lemma:concentration} is proved.

\section{Proof of Lemma~\ref{lemma:bound-time-large-queue}}
\label{app:proof-lemma-bound-time-large-queue}

We use Lyapunov function method to prove Lemma~\ref{lemma:bound-time-large-queue}. Define a Lyapunov function $V_{\mathrm{c}}(t)\coloneqq \sum_i Q_{\mathrm{c},i}^2(t)$. We will bound the Lyapunov drift $V_{\mathrm{c}}(t+1) - V_{\mathrm{c}}(t)$ given $\boldsymbol{Q}_{\mathrm{c}}(t)$, where $\boldsymbol{Q}_{\mathrm{c}}(t)$ denotes a vector containing $Q_{\mathrm{c},i}(t)$ for all $i$. Using the queue dynamics \eqref{equ:dynamics}, we have
 \begin{align*}
     & \expt [V_{\mathrm{c}}(t+1) - V_{\mathrm{c}}(t) | \boldsymbol{Q}_{\mathrm{c}}(t) ]\nonumber\\
     = & \expt \biggl[ \sum_i ( Q_{\mathrm{c},i}^2(t+1) - Q_{\mathrm{c},i}^2(t) ) \biggl|\biggr. \boldsymbol{Q}_{\mathrm{c}}(t)\biggr]\nonumber\\
     = & \expt \Biggl[ \sum_i \biggl[ \biggl(  Q_{\mathrm{c},i}(t) + A_{\mathrm{c},i}(t) - \sum_{j\in {\cal E}_{\mathrm{c},i}} X_{i,j}(t)  \biggr)^2 - Q_{\mathrm{c},i}^2(t) \biggr] \Biggl|\Biggr. \boldsymbol{Q}_{\mathrm{c}}(t)\Biggr]\nonumber\\
     = & \expt \Biggl[ \sum_i \biggl[ \biggl( A_{\mathrm{c},i}(t) - \sum_{j\in {\cal E}_{\mathrm{c},i}} X_{i,j}(t)  \biggr)^2 
     + 2 Q_{\mathrm{c},i}(t) \biggl( A_{\mathrm{c},i}(t) - \sum_{j\in {\cal E}_{\mathrm{c},i}} X_{i,j}(t) \biggr) \biggr] \Biggl|\Biggr. \boldsymbol{Q}_{\mathrm{c}}(t)\Biggr].
 \end{align*}
 Since $A_{\mathrm{c},i}(t) \in \{0,1\}$ and $X_{i,j}(t)\in \{0,1\}$, we have $(A_{\mathrm{c},i}(t) - \sum_{j\in {\cal E}_{\mathrm{c},i}} X_{i,j}(t))^2 \le \max\{ 1, |{\cal E}_{\mathrm{c},i}|^2 \} =  |{\cal E}_{\mathrm{c},i}|^2 $. Hence, we have
\begin{align}\label{equ:drift-1}
     \expt [V_{\mathrm{c}}(t+1) - V_{\mathrm{c}}(t) | \boldsymbol{Q}_{\mathrm{c}}(t) ] \le 
     \sum_i |{\cal E}_{\mathrm{c},i}|^2 + 2 \sum_i Q_{\mathrm{c},i}(t) \biggl( \expt [ A_{\mathrm{c},i}(t) | \boldsymbol{Q}_{\mathrm{c}}(t) ] -  \expt \biggl[\sum_{j\in {\cal E}_{\mathrm{c},i}} X_{i,j}(t) \biggl|\biggr. \boldsymbol{Q}_{\mathrm{c}}(t) \biggr] \biggr).
\end{align}
Recall that $\lambda'_i(t),\mu'_j(t)$ are the arrival rates as if there is no queue length control.
For $k\ge 2$, i.e., $t\ge 2MN+1$, the arrival is zero when the queue length is greater than or equal to the threshold $q^{\mathrm{th}}$. Hence, for $t\ge 2MN + 1$, the term $\expt [ A_{\mathrm{c},i}(t) | \boldsymbol{Q}_{\mathrm{c}}(t) ]$ in \eqref{equ:drift-1} can be written as
\begin{align}\label{equ:drift-1-1}
    \expt [ A_{\mathrm{c},i}(t) | \boldsymbol{Q}_{\mathrm{c}}(t) ] = & \mathbb{1} \left\{ Q_{\mathrm{c},i}(t) < q^{\mathrm{th}} \right\} \expt [ A_{\mathrm{c},i}(t) | \boldsymbol{Q}_{\mathrm{c}}(t) ]\nonumber\\
    = & \mathbb{1} \left\{ Q_{\mathrm{c},i}(t) < q^{\mathrm{th}} \right\}
    \expt \bigl[  
        \expt [ A_{\mathrm{c},i}(t) 
        | \lambda'_i(t), \boldsymbol{Q}_{\mathrm{c}}(t) ]
    \bigl|\bigr. \boldsymbol{Q}_{\mathrm{c}}(t) \bigr]\nonumber\\
    = & \mathbb{1} \left\{ Q_{\mathrm{c},i}(t) < q^{\mathrm{th}} \right\}
    \expt \bigl[  \lambda'_i(t)  \bigl|\bigr. \boldsymbol{Q}_{\mathrm{c}}(t) \bigr],
\end{align}
where the second equality is by the law of iterated expectation and the last equality holds since the arrival $A_{\mathrm{c},i}(t)$ is independent of the queue length $\boldsymbol{Q}_{\mathrm{c}}(t)$ given $\lambda'_i(t)$ and that $Q_{\mathrm{c},i}(t) < q^{\mathrm{th}}$.
Next, we will rewrite the term $\sum_i Q_{\mathrm{c},i}(t) \expt [\sum_{j\in {\cal E}_{\mathrm{c},i}} X_{i,j}(t) | \boldsymbol{Q}_{\mathrm{c}}(t) ]$ in \eqref{equ:drift-1}. We have
\begin{align}\label{equ:drift-1-2}
    \sum_i Q_{\mathrm{c},i}(t) \expt \biggl[\sum_{j\in {\cal E}_{\mathrm{c},i}} X_{i,j}(t) \biggl|\biggr. \boldsymbol{Q}_{\mathrm{c}}(t) \biggr]
    = & \expt \biggl[ \sum_i \sum_{j\in {\cal E}_{\mathrm{c},i}} Q_{\mathrm{c},i}(t) X_{i,j}(t)
    \biggl|\biggr. \boldsymbol{Q}_{\mathrm{c}}(t) \biggr]\nonumber\\
    = & \expt \biggl[ \sum_j \sum_{i\in {\cal E}_{\mathrm{s},j}} Q_{\mathrm{c},i}(t) X_{i,j}(t)
    \biggl|\biggr. \boldsymbol{Q}_{\mathrm{c}}(t) \biggr] \nonumber\\
    = & \expt \biggl[ \sum_j \sum_{i\in {\cal E}_{\mathrm{s},j}}
    \mathbb{1} \left\{ Q_{\mathrm{c},i}(t) > 0 \right\}
    Q_{\mathrm{c},i}(t) X_{i,j}(t)
    \biggl|\biggr. \boldsymbol{Q}_{\mathrm{c}}(t) \biggr],
\end{align}
where the second equality is by exchanging the summations. Note that by Lemma~\ref{lemma:empty-queues} we know that if $Q_{\mathrm{c},i}(t)>0$, then $Q_{\mathrm{s},j}(t)=0$ for all $j\in {\cal E}_{\mathrm{c},i}$. Hence, in this case, $X_{i,j}(t)=1$ (there is one match) only if $A_{\mathrm{s},j}(t)=1$. Therefore, combining this observation with \eqref{equ:drift-1-2}, we have
\begin{align}\label{equ:drift-1-3}
    \sum_i Q_{\mathrm{c},i}(t) \expt \biggl[\sum_{j\in {\cal E}_{\mathrm{c},i}} X_{i,j}(t) \biggl|\biggr. \boldsymbol{Q}_{\mathrm{c}}(t) \biggr]
    = \expt \biggl[ \sum_j \sum_{i\in {\cal E}_{\mathrm{s},j}}
    \mathbb{1} \left\{ Q_{\mathrm{c},i}(t) > 0 \right\} 
    \mathbb{1} \left\{ A_{\mathrm{s},j}(t) = 1 \right\} 
    Q_{\mathrm{c},i}(t) X_{i,j}(t)
    \biggl|\biggr. \boldsymbol{Q}_{\mathrm{c}}(t) \biggr].
\end{align}
Note that $Q_{\mathrm{s},j}(t)=0$ and $A_{\mathrm{s},j}(t) = 1$ imply that 
$\sum_{i\in {\cal E}_{\mathrm{s},j}} X_{i,j}(t)\le 1$. Also note that the summation $\sum_{i\in {\cal E}_{\mathrm{s},j}} \mathbb{1} \left\{ Q_{\mathrm{c},i}(t) > 0 \right\} \mathbb{1} \left\{ A_{\mathrm{s},j}(t) = 1 \right\} Q_{\mathrm{c},i}(t) X_{i,j}(t)$ is nonzero only if $\sum_{i\in {\cal E}_{\mathrm{s},j}} Q_{\mathrm{c},i}(t) > 0 $. Hence, we have $\sum_{i\in {\cal E}_{\mathrm{s},j}} X_{i,j}(t) = 1$.
Recall that $i^*_j(t)$ is the customer queue picked by Algorithm~\ref{alg:matching} for server queue $j$ at time $t$.
Hence, \eqref{equ:drift-1-3} can be further written as
\begin{align}\label{equ:drift-1-4}
    \sum_i Q_{\mathrm{c},i}(t) \expt \biggl[\sum_{j\in {\cal E}_{\mathrm{c},i}} X_{i,j}(t) \biggl|\biggr. \boldsymbol{Q}_{\mathrm{c}}(t) \biggr]
    = & \expt \biggl[ \sum_j
    \mathbb{1} \left\{ A_{\mathrm{s},j}(t) = 1 \right\} 
    \mathbb{1} \biggl\{ \sum_{i\in {\cal E}_{\mathrm{s},j}} Q_{\mathrm{c},i}(t)>0 \biggr\}
    Q_{\mathrm{c},i^*_j(t)}(t) 
    \biggl|\biggr. \boldsymbol{Q}_{\mathrm{c}}(t) \biggr]\nonumber\\
    = &  \sum_j
    \mathbb{1} \biggl\{ \sum_{i\in {\cal E}_{\mathrm{s},j}} Q_{\mathrm{c},i}(t)>0 \biggr\}
    \expt \bigl[ \mathbb{1} \left\{ A_{\mathrm{s},j}(t) = 1 \right\}  Q_{\mathrm{c},i^*_j(t)}(t) 
    \bigl|\bigr. \boldsymbol{Q}_{\mathrm{c}}(t) \bigr].
\end{align}
By the law of iterated expectation, the term $\expt [ \mathbb{1} \left\{ A_{\mathrm{s},j}(t) = 1 \right\}  Q_{\mathrm{c},i^*_j(t)}(t) | \boldsymbol{Q}_{\mathrm{c}}(t) ]$ in \eqref{equ:drift-1-4} can be written as
\begin{align}\label{equ:drift-1-5}
    \expt \bigl[ \mathbb{1} \left\{ A_{\mathrm{s},j}(t) = 1 \right\}  Q_{\mathrm{c},i^*_j(t)}(t) 
    \bigl|\bigr. \boldsymbol{Q}_{\mathrm{c}}(t) \bigr] =  & \expt \bigl[ A_{\mathrm{s},j}(t)  Q_{\mathrm{c},i^*_j(t)}(t) 
    \bigl|\bigr. \boldsymbol{Q}_{\mathrm{c}}(t) \bigr]\nonumber\\
    = & \expt \left[  \expt \bigl[ A_{\mathrm{s},j}(t)  Q_{\mathrm{c},i^*_j(t)}(t) 
     \bigl|\bigr. \boldsymbol{Q}_{\mathrm{c}}(t), i^*_j(t), \mu'_j(t)  \bigr] \Bigl|\bigr. \boldsymbol{Q}_{\mathrm{c}}(t) \right] \nonumber\\
    = & \expt \left[   \expt \bigl[ A_{\mathrm{s},j}(t)   
     \bigl|\bigr. \boldsymbol{Q}_{\mathrm{c}}(t), i^*_j(t), \mu'_j(t)  \bigr] Q_{\mathrm{c},i^*_j(t)}(t) \Bigl|\bigr. \boldsymbol{Q}_{\mathrm{c}}(t) \right].
\end{align}
Note that when $\sum_{i\in {\cal E}_{\mathrm{s},j}} Q_{\mathrm{c},i}(t)>0$, we have $Q_{\mathrm{s},j}(t)=0<q^{\mathrm{th}}$ by Lemma~\ref{lemma:empty-queues}. Hence, for $k\ge 2$, i.e., $t\ge 2MN$, we know that the arrival rate for server queue $j$ at time $t$ is $\mu'_j(t)$ when $\sum_{i\in {\cal E}_{\mathrm{s},j}} Q_{\mathrm{c},i}(t)>0$. Therefore, when $\sum_{i\in {\cal E}_{\mathrm{s},j}} Q_{\mathrm{c},i}(t)>0$, we have
\begin{align}\label{equ:drift-1-6}
    \expt \bigl[ A_{\mathrm{s},j}(t)   
     \bigl|\bigr. \boldsymbol{Q}_{\mathrm{c}}(t), i^*_j(t), \mu'_j(t)  \bigr] = \mu'_j(t),
\end{align}
since $A_{\mathrm{s},j}(t)$ is independent of $\boldsymbol{Q}_{\mathrm{c}}(t), i^*_j(t)$ given $\mu'_j(t)$. Combining \eqref{equ:drift-1-4}, \eqref{equ:drift-1-5}, and \eqref{equ:drift-1-6}, we have
\begin{align}\label{equ:drift-1-7}
    \sum_i Q_{\mathrm{c},i}(t) \expt \biggl[\sum_{j\in {\cal E}_{\mathrm{c},i}} X_{i,j}(t) \biggl|\biggr. \boldsymbol{Q}_{\mathrm{c}}(t) \biggr]
    = & \sum_j
    \mathbb{1} \biggl\{ \sum_{i\in {\cal E}_{\mathrm{s},j}} Q_{\mathrm{c},i}(t)>0 \biggr\}
    \expt \left[   \mu'_j(t) Q_{\mathrm{c},i^*_j(t)}(t) \Bigl|\bigr. \boldsymbol{Q}_{\mathrm{c}}(t) \right]\nonumber\\
    = & \sum_j
    \expt \left[   \mu'_j(t) Q_{\mathrm{c},i^*_j(t)}(t) \Bigl|\bigr. \boldsymbol{Q}_{\mathrm{c}}(t) \right],
\end{align}
where the last equality holds since $Q_{\mathrm{c},i^*_j(t)}(t)=0$ if $\sum_{i\in {\cal E}_{\mathrm{s},j}} Q_{\mathrm{c},i}(t)=0$. From \eqref{equ:drift-1}, \eqref{equ:drift-1-1}, and \eqref{equ:drift-1-7}, we have
\begin{align}\label{equ:drift-2}
    & \expt [V_{\mathrm{c}}(t+1) - V_{\mathrm{c}}(t) | \boldsymbol{Q}_{\mathrm{c}}(t) ] \nonumber\\
    \le & 
     \sum_i |{\cal E}_{\mathrm{c},i}|^2 + 2 \sum_i Q_{\mathrm{c},i}(t) \mathbb{1} \left\{ Q_{\mathrm{c},i}(t) < q^{\mathrm{th}} \right\}
    \expt \bigl[  \lambda'_i(t)  \bigl|\bigr. \boldsymbol{Q}_{\mathrm{c}}(t) \bigr]  - 2 \sum_j
    \expt \left[   \mu'_j(t) Q_{\mathrm{c},i^*_j(t)}(t) \Bigl|\bigr. \boldsymbol{Q}_{\mathrm{c}}(t) \right].
\end{align}
We then add the indicator of the ``good'' event ${\cal C}$ into \eqref{equ:drift-2}, obtaining
\begin{align}\label{equ:drift-2-1}
    \expt \bigl[  \lambda'_i(t)  \bigl|\bigr. \boldsymbol{Q}_{\mathrm{c}}(t) \bigr]
    = & \expt \bigl[  \lambda'_i(t) \mathbb{1}_{{\cal C}} \bigl|\bigr. \boldsymbol{Q}_{\mathrm{c}}(t) \bigr]
    + \expt \bigl[  \lambda'_i(t) \mathbb{1}_{{\cal C}^{\mathrm{c}}} \bigl|\bigr. \boldsymbol{Q}_{\mathrm{c}}(t) \bigr]\nonumber\\
    \le & \expt \bigl[  \lambda'_i(t) \mathbb{1}_{{\cal C}} \bigl|\bigr. \boldsymbol{Q}_{\mathrm{c}}(t) \bigr]
    + \expt \bigl[ \mathbb{1}_{{\cal C}^{\mathrm{c}}} \bigl|\bigr. \boldsymbol{Q}_{\mathrm{c}}(t) \bigr],
\end{align}
where the inequality is due to the fact that $\lambda'_i(t) \le 1$. Similarly, we have
\begin{align}\label{equ:drift-2-2}
    \expt \left[   \mu'_j(t) Q_{\mathrm{c},i^*_j(t)}(t) \Bigl|\bigr. \boldsymbol{Q}_{\mathrm{c}}(t) \right] = & \expt \left[   \mu'_j(t) Q_{\mathrm{c},i^*_j(t)}(t) \mathbb{1}_{{\cal C}}
    \Bigl|\bigr. \boldsymbol{Q}_{\mathrm{c}}(t)  \right]
    + \expt \left[   \mu'_j(t) Q_{\mathrm{c},i^*_j(t)}(t) \mathbb{1}_{{\cal C}^{\mathrm{c}}} 
    \Bigl|\bigr. \boldsymbol{Q}_{\mathrm{c}}(t) \right]\nonumber\\
    \ge & \expt \left[   \mu'_j(t) Q_{\mathrm{c},i^*_j(t)}(t) \mathbb{1}_{{\cal C}}
    \Bigl|\bigr. \boldsymbol{Q}_{\mathrm{c}}(t)  \right].
\end{align}
Substituting \eqref{equ:drift-2-1} and \eqref{equ:drift-2-2} into \eqref{equ:drift-2}, we have
\begin{align}\label{equ:drift-3}
    & \expt [V_{\mathrm{c}}(t+1) - V_{\mathrm{c}}(t) | \boldsymbol{Q}_{\mathrm{c}}(t) ] \nonumber\\
    \le & 
     \sum_i |{\cal E}_{\mathrm{c},i}|^2 + 2 \sum_i Q_{\mathrm{c},i}(t) \mathbb{1} \left\{ Q_{\mathrm{c},i}(t) < q^{\mathrm{th}} \right\}
     \expt \bigl[  \lambda'_i(t) \mathbb{1}_{{\cal C}} \bigl|\bigr. \boldsymbol{Q}_{\mathrm{c}}(t) \bigr]\nonumber\\
    & + 2 \sum_i Q_{\mathrm{c},i}(t) \mathbb{1} \left\{ Q_{\mathrm{c},i}(t) < q^{\mathrm{th}} \right\} \expt \bigl[ \mathbb{1}_{{\cal C}^{\mathrm{c}}} \bigl|\bigr. \boldsymbol{Q}_{\mathrm{c}}(t) \bigr]  - 2 \sum_j
    \expt \left[   \mu'_j(t) Q_{\mathrm{c},i^*_j(t)}(t) \mathbb{1}_{{\cal C}}
    \Bigl|\bigr. \boldsymbol{Q}_{\mathrm{c}}(t)  \right]\nonumber\\
    \le & \sum_i |{\cal E}_{\mathrm{c},i}|^2 + 2 \sum_i Q_{\mathrm{c},i}(t) \mathbb{1} \left\{ Q_{\mathrm{c},i}(t) < q^{\mathrm{th}} \right\}
     \expt \bigl[  \lambda'_i(t) \mathbb{1}_{{\cal C}} \bigl|\bigr. \boldsymbol{Q}_{\mathrm{c}}(t) \bigr]\nonumber\\
     & + 2 I q^{\mathrm{th}} \expt \bigl[ \mathbb{1}_{{\cal C}^{\mathrm{c}}} \bigl|\bigr. \boldsymbol{Q}_{\mathrm{c}}(t) \bigr]  - 2 \sum_j
    \expt \left[   \mu'_j(t) Q_{\mathrm{c},i^*_j(t)}(t) \mathbb{1}_{{\cal C}}
    \Bigl|\bigr. \boldsymbol{Q}_{\mathrm{c}}(t)  \right],
\end{align}
where the last inequality is due to the fact that $ Q_{\mathrm{c},i}(t) \mathbb{1} \left\{ Q_{\mathrm{c},i}(t) < q^{\mathrm{th}} \right\} \le q^{\mathrm{th}}$.

Next, we will relate $\lambda'_i(t)$ and $\mu'_j(t)$ to the arrival rates $\sum_{j\in {\cal E}_{\mathrm{c},i}} x_{i,j}(k(t))$ and $\sum_{i\in {\cal E}_{\mathrm{s},j}} x_{i,j}(k(t))$, where we recall that $k(t)$ denotes the outer iteration at time $t$.
We have the following lemma under the pricing algorithm~\ref{alg:pricing} and bisection algorithm~\ref{alg:bisection}.
\begin{lemma}\label{lemma:bisection-error}
    Let Assumption~\ref{assum:1} and Assumption~\ref{assum:4} hold. Suppose $2e_{\mathrm{c},i}\le p_{\mathrm{c},i,\max} - p_{\mathrm{c},i,\min}$ and $2e_{\mathrm{s},j}\le p_{\mathrm{s},j,\max} - p_{\mathrm{s},j,\min}$. Then under the event ${\cal C}$, for $t\ge 2MN+1$, we have
    \begin{align*}
    \biggl|\lambda'_i(t) - \sum_{j\in {\cal E}_{\mathrm{c},i}} x_{i,j}(k(t))\biggr|
    \le & 2 L_{F^{-1}_i} e_{\mathrm{c},i} + 2\epsilon \left[1 + L_{F^{-1}_i} \left( p_{\mathrm{c},i,\max} - p_{\mathrm{c},i,\min}\right)\right]  + \delta \sqrt{|{\cal E}_{\mathrm{c},i}|}, \mbox{ for all } i;\nonumber\\
    \biggl|\mu'_j(t) - \sum_{i\in {\cal E}_{\mathrm{s},j}} x_{i,j}(k(t))\biggr|
    \le & 2 L_{G^{-1}_j} e_{\mathrm{s},j} + 2\epsilon \left[1 + L_{G^{-1}_j} \left( p_{\mathrm{s},j,\max} - p_{\mathrm{s},j,\min}\right)\right]  + \delta \sqrt{|{\cal E}_{\mathrm{s},j}|}, \mbox{ for all } j.
    \end{align*}
\end{lemma}
\noindent Lemma~\ref{lemma:bisection-error} shows that the difference between the actual arrival rate of any queue without queue length control and the target arrival rate generated in the optimization process is bounded and small. 
Proof of Lemma~\ref{lemma:bisection-error} can be found in Appendix~\ref{app:proof-lemma-bisection-error}.

For sufficiently large $T$, we have $2e_{\mathrm{c},i}\le p_{\mathrm{c},i,\max} - p_{\mathrm{c},i,\min}$ and $2e_{\mathrm{s},j}\le p_{\mathrm{s},j,\max} - p_{\mathrm{s},j,\min}$ by the definitions of $e_{\mathrm{c},i}$ and $e_{\mathrm{s},j}$ in Appendix~\ref{app:exact-expression} and $\epsilon < \delta$. Then using Lemma~\ref{lemma:bisection-error} and the fact that $\mathbb{1}_{\cal C}\le 1$, \eqref{equ:drift-3} can be bounded by
\begin{align*}
    & \expt [V_{\mathrm{c}}(t+1) - V_{\mathrm{c}}(t) | \boldsymbol{Q}_{\mathrm{c}}(t) ] \nonumber\\
    \le & \sum_i |{\cal E}_{\mathrm{c},i}|^2 + 2 \sum_i Q_{\mathrm{c},i}(t) \mathbb{1} \left\{ Q_{\mathrm{c},i}(t) < q^{\mathrm{th}} \right\}
     \expt \biggl[  \sum_{j\in {\cal E}_{\mathrm{c},i}} x_{i,j}(k(t)) \mathbb{1}_{{\cal C}} \biggl|\biggr. \boldsymbol{Q}_{\mathrm{c}}(t) \biggr]\nonumber\\
    & + 2 \sum_i Q_{\mathrm{c},i}(t) \mathbb{1} \left\{ Q_{\mathrm{c},i}(t) < q^{\mathrm{th}} \right\} \left( 2 L_{F^{-1}_i} e_{\mathrm{c},i} + 2\epsilon \left[1 + L_{F^{-1}_i} \left( p_{\mathrm{c},i,\max} - p_{\mathrm{c},i,\min}\right)\right]  + \delta \sqrt{|{\cal E}_{\mathrm{c},i}|} \right)\nonumber\\
    & + 2 I q^{\mathrm{th}} \expt \bigl[ \mathbb{1}_{{\cal C}^{\mathrm{c}}} \bigl|\bigr. \boldsymbol{Q}_{\mathrm{c}}(t) \bigr]  - 2 \sum_j
    \expt \left[   \sum_{i\in {\cal E}_{\mathrm{s},j}} x_{i,j}(k(t))  Q_{\mathrm{c},i^*_j(t)}(t) \mathbb{1}_{{\cal C}}
    \Bigl|\bigr. \boldsymbol{Q}_{\mathrm{c}}(t)  \right]\nonumber\\
    & + 2 \sum_j
    \left(2 L_{G^{-1}_j} e_{\mathrm{s},j} + 2\epsilon \left[1 + L_{G^{-1}_j} \left( p_{\mathrm{s},j,\max} - p_{\mathrm{s},j,\min}\right)\right]  + \delta \sqrt{|{\cal E}_{\mathrm{s},j}|}\right)
    \expt \left[  Q_{\mathrm{c},i^*_j(t)}(t) \mathbb{1}_{{\cal C}}
    \Bigl|\bigr. \boldsymbol{Q}_{\mathrm{c}}(t)  \right],
\end{align*}
for $t\ge 2MN+1$.
Since $Q_{\mathrm{c},i}(t) \mathbb{1} \left\{ Q_{\mathrm{c},i}(t) < q^{\mathrm{th}} \right\} \le q^{\mathrm{th}}$ and $Q_{\mathrm{c},i^*_j(t)}(t) \le \sum_i Q_{\mathrm{c},i}(t)$, we further obtain for $t\ge 2MN+1$,
\begin{align}
    & \expt [V_{\mathrm{c}}(t+1) - V_{\mathrm{c}}(t) | \boldsymbol{Q}_{\mathrm{c}}(t) ] \nonumber\\
    \le &  2 \sum_i Q_{\mathrm{c},i}(t) \mathbb{1} \left\{ Q_{\mathrm{c},i}(t) < q^{\mathrm{th}} \right\}
     \expt \biggl[  \sum_{j\in {\cal E}_{\mathrm{c},i}} x_{i,j}(k(t)) \mathbb{1}_{{\cal C}} \biggl|\biggr. \boldsymbol{Q}_{\mathrm{c}}(t) \biggr]\label{equ:drift-4-1}\\
    & - 2 \sum_j
    \expt \left[   \sum_{i\in {\cal E}_{\mathrm{s},j}} x_{i,j}(k(t))  Q_{\mathrm{c},i^*_j(t)}(t) \mathbb{1}_{{\cal C}}
    \Bigl|\bigr. \boldsymbol{Q}_{\mathrm{c}}(t)  \right]\label{equ:drift-4-2}\\
    & + \sum_i |{\cal E}_{\mathrm{c},i}|^2 + 2 q^{\mathrm{th}} \sum_i \left( 2 L_{F^{-1}_i} e_{\mathrm{c},i} + 2\epsilon \left[1 + L_{F^{-1}_i} \left( p_{\mathrm{c},i,\max} - p_{\mathrm{c},i,\min}\right)\right]  + \delta \sqrt{|{\cal E}_{\mathrm{c},i}|} \right)\label{equ:drift-4-3}\\
    & + 2 I q^{\mathrm{th}} \expt \bigl[ \mathbb{1}_{{\cal C}^{\mathrm{c}}} \bigl|\bigr. \boldsymbol{Q}_{\mathrm{c}}(t) \bigr] \label{equ:drift-4-4}\\
    & + 2 \sum_i Q_{\mathrm{c},i}(t) \expt \left[ \mathbb{1}_{{\cal C}}
    \Bigl|\bigr. \boldsymbol{Q}_{\mathrm{c}}(t)  \right] \sum_j
    \left(2 L_{G^{-1}_j} e_{\mathrm{s},j} + 2\epsilon \left[1 + L_{G^{-1}_j} \left( p_{\mathrm{s},j,\max} - p_{\mathrm{s},j,\min}\right)\right]  + \delta \sqrt{|{\cal E}_{\mathrm{s},j}|}\right)\label{equ:drift-4-4-end},
\end{align}
where \eqref{equ:drift-4-1} and \eqref{equ:drift-4-2} can be combined as
\begin{align}\label{equ:drift-4-5}
    \eqref{equ:drift-4-1} + \eqref{equ:drift-4-2} 
    = & 2\sum_{(i,j)\in {\cal E}} 
    \expt \biggl[ x_{i,j}(k(t)) \mathbb{1}_{{\cal C}} 
    \Bigl[ Q_{\mathrm{c},i}(t) \mathbb{1} \left\{ Q_{\mathrm{c},i}(t) < q^{\mathrm{th}} \right\} 
    -  Q_{\mathrm{c},i^*_j(t)}(t) \Bigr] 
    \biggl|\biggr. \boldsymbol{Q}_{\mathrm{c}}(t) \biggr]\nonumber\\
    = & 2\sum_{(i,j)\in {\cal E}} 
    \expt \biggl[ x_{i,j}(k(t)) \mathbb{1}_{{\cal C}} 
    \Bigl[ Q_{\mathrm{c},i}(t) - Q_{\mathrm{c},i}(t) \mathbb{1} \left\{ Q_{\mathrm{c},i}(t) \ge q^{\mathrm{th}} \right\} 
    -  Q_{\mathrm{c},i^*_j(t)}(t) \Bigr] 
    \biggl|\biggr. \boldsymbol{Q}_{\mathrm{c}}(t) \biggr].
\end{align}
Recall the matching algorithm~\ref{alg:matching}. $i^*_j(t)$ is the customer queue chosen to match the arrival in the server queue $j$ at time $t$ and $i^*_j(t)=\argmax_{i: (i,j)\in {\cal E}} \tilde{Q}_{\mathrm{c}, i}$, where $\tilde{Q}_{\mathrm{c}, i}$ is an intermediate queue length during time slot $t$. Note that by the matching algorithm~\ref{alg:matching}, we know that during one time slot there are at most $J$ departures from any customer queue $i$. Hence, $Q_{\mathrm{c},i}(t) - \tilde{Q}_{\mathrm{c}, i} \le J$ for any $i$. Hence, we have
\begin{align}\label{equ:drift-4-6}
    Q_{\mathrm{c}, i}(t) - Q_{\mathrm{c},i^*_j(t)}(t) 
    \le &
    \max_{i: (i,j)\in {\cal E}} Q_{\mathrm{c}, i}(t) - Q_{\mathrm{c},i^*_j(t)}(t) 
    \le \max_{i: (i,j)\in {\cal E}} \tilde{Q}_{\mathrm{c}, i} + J - Q_{\mathrm{c},i^*_j(t)}(t) = J.
\end{align}
From \eqref{equ:drift-4-5} and \eqref{equ:drift-4-6}, we have
\begin{align}\label{equ:drift-4-7}
    \eqref{equ:drift-4-1} + \eqref{equ:drift-4-2} 
    \le & 2 \sum_{(i,j)\in {\cal E}} 
    \expt \biggl[ x_{i,j}(k(t)) \mathbb{1}_{{\cal C}} 
    \Bigl[ J - Q_{\mathrm{c},i}(t) \mathbb{1} \left\{ Q_{\mathrm{c},i}(t) \ge q^{\mathrm{th}} \right\} 
    \Bigr] 
    \biggl|\biggr. \boldsymbol{Q}_{\mathrm{c}}(t) \biggr].
\end{align}
From the pricing algorithm~\ref{alg:pricing}, we know that $\boldsymbol{x}(k)\in {\cal D}'$ for all $k$ and we know that ${\cal D}' \subseteq {\cal D}$ by Lemma~\ref{lemma:shrunk-set}. Hence, we have $\sum_{j\in {\cal E}_{\mathrm{c},i}} x_{i,j}(k(t))\in [a_{\min}, 1]$ and $\sum_{i\in {\cal E}_{\mathrm{s},j}} x_{i,j}(k(t))\in [a_{\min}, 1]$. Since $\sum_{(i,j)\in {\cal E}} x_{i,j}(k(t)) = \sum_i \sum_{j\in {\cal E}_{\mathrm{c},i}} x_{i,j}(k(t)) = \sum_j \sum_{i\in {\cal E}_{\mathrm{s},j}} x_{i,j}(k(t))$, we have $\sum_{(i,j)\in {\cal E}} x_{i,j}(k(t)) \le \max\{I, J\}$. From this inequality and \eqref{equ:drift-4-7}, we have
\begin{align}\label{equ:drift-4-8}
    \eqref{equ:drift-4-1} + \eqref{equ:drift-4-2} \le & 2 \max\{I, J\} J 
    - 2 \sum_{(i,j)\in {\cal E}} 
    \expt \biggl[ x_{i,j}(k(t)) \mathbb{1}_{{\cal C}} 
    Q_{\mathrm{c},i}(t) \mathbb{1} \left\{ Q_{\mathrm{c},i}(t) \ge q^{\mathrm{th}} \right\} 
    \biggl|\biggr. \boldsymbol{Q}_{\mathrm{c}}(t) \biggr]\nonumber\\
    = & 2 \max\{I, J\} J 
    - 2 \sum_i \expt \biggl[ \sum_{j\in {\cal E}_{\mathrm{c},i}} x_{i,j}(k(t))
    \mathbb{1}_{{\cal C}}
    Q_{\mathrm{c},i}(t)
    \mathbb{1} \left\{ Q_{\mathrm{c},i}(t) \ge q^{\mathrm{th}} \right\}
    \biggl|\biggr. \boldsymbol{Q}_{\mathrm{c}}(t) \biggr]\nonumber\\
    \le & 2 \max\{I, J\} J 
    - 2 a_{\min} \sum_i \expt \biggl[ 
    \mathbb{1}_{{\cal C}}
    Q_{\mathrm{c},i}(t)
    \mathbb{1} \left\{ Q_{\mathrm{c},i}(t) \ge q^{\mathrm{th}} \right\}
    \biggl|\biggr. \boldsymbol{Q}_{\mathrm{c}}(t) \biggr]
\end{align}
where in the equality we divide the summation $\sum_{(i,j)\in {\cal E}}$ into two summations, and in the last inequality we use $\sum_{j\in {\cal E}_{\mathrm{c},i}} x_{i,j}(k(t)) \ge a_{\min}$.
Combining \eqref{equ:drift-4-1}-\eqref{equ:drift-4-4-end} and \eqref{equ:drift-4-8}, we have for $t\ge 2MN+1$,
\begin{align*}
    & \expt [V_{\mathrm{c}}(t+1) - V_{\mathrm{c}}(t) | \boldsymbol{Q}_{\mathrm{c}}(t) ] \nonumber\\
    \le & 2 \max\{I, J\} J 
    - 2 a_{\min} \sum_i \expt \biggl[ 
    \mathbb{1}_{{\cal C}}
    Q_{\mathrm{c},i}(t)
    \mathbb{1} \left\{ Q_{\mathrm{c},i}(t) \ge q^{\mathrm{th}} \right\}
    \biggl|\biggr. \boldsymbol{Q}_{\mathrm{c}}(t) \biggr]\nonumber\\
    & + \sum_i |{\cal E}_{\mathrm{c},i}|^2 + 2 q^{\mathrm{th}} \sum_i \left( 2 L_{F^{-1}_i} e_{\mathrm{c},i} + 2\epsilon \left[1 + L_{F^{-1}_i} \left( p_{\mathrm{c},i,\max} - p_{\mathrm{c},i,\min}\right)\right]  + \delta \sqrt{|{\cal E}_{\mathrm{c},i}|} \right)\nonumber\\
    & + 2 I q^{\mathrm{th}} \expt \bigl[ \mathbb{1}_{{\cal C}^{\mathrm{c}}} \bigl|\bigr. \boldsymbol{Q}_{\mathrm{c}}(t) \bigr] \nonumber\\
    & + 2 \sum_i Q_{\mathrm{c},i}(t) \expt \left[ \mathbb{1}_{{\cal C}}
    \Bigl|\bigr. \boldsymbol{Q}_{\mathrm{c}}(t)  \right] \sum_j
    \left(2 L_{G^{-1}_j} e_{\mathrm{s},j} + 2\epsilon \left[1 + L_{G^{-1}_j} \left( p_{\mathrm{s},j,\max} - p_{\mathrm{s},j,\min}\right)\right]  + \delta \sqrt{|{\cal E}_{\mathrm{s},j}|}\right).
\end{align*}
Taking expectations on both sides, we have
\begin{align*}
    & \expt [V_{\mathrm{c}}(t+1) - V_{\mathrm{c}}(t) ] \nonumber\\
    \le & 2 \max\{I, J\} J 
    - 2 a_{\min} \sum_i \expt \Bigl[ 
    \mathbb{1}_{{\cal C}}
    Q_{\mathrm{c},i}(t)
    \mathbb{1} \Bigl\{ Q_{\mathrm{c},i}(t) \ge q^{\mathrm{th}} \Bigr\}
    \Bigr]\nonumber\\
    & + \sum_i |{\cal E}_{\mathrm{c},i}|^2 + 2 q^{\mathrm{th}} \sum_i \left( 2 L_{F^{-1}_i} e_{\mathrm{c},i} + 2\epsilon \left[1 + L_{F^{-1}_i} \left( p_{\mathrm{c},i,\max} - p_{\mathrm{c},i,\min}\right)\right]  + \delta \sqrt{|{\cal E}_{\mathrm{c},i}|} \right)\nonumber\\
    & + 2 I q^{\mathrm{th}} \expt \bigl[ \mathbb{1}_{{\cal C}^{\mathrm{c}}} \bigr] 
    + 2 \sum_i \expt \bigl[ \mathbb{1}_{{\cal C}} Q_{\mathrm{c},i}(t) \bigr] 
    \sum_j
    \left(2 L_{G^{-1}_j} e_{\mathrm{s},j} + 2\epsilon \left[1 + L_{G^{-1}_j} \left( p_{\mathrm{s},j,\max} - p_{\mathrm{s},j,\min}\right)\right]  + \delta \sqrt{|{\cal E}_{\mathrm{s},j}|}\right).
\end{align*}
Note that for the last term in the above, we have
\begin{align*}
    Q_{\mathrm{c},i}(t)
    = & Q_{\mathrm{c},i}(t) \mathbb{1} \left\{ Q_{\mathrm{c},i}(t) \ge q^{\mathrm{th}} \right\} + Q_{\mathrm{c},i}(t) \mathbb{1} \left\{ Q_{\mathrm{c},i}(t) < q^{\mathrm{th}} \right\}\nonumber\\
    \le & Q_{\mathrm{c},i}(t) \mathbb{1} \left\{ Q_{\mathrm{c},i}(t) \ge q^{\mathrm{th}} \right\} + q^{\mathrm{th}}.
\end{align*}
Hence, we have
\begin{align}
    & \expt [V_{\mathrm{c}}(t+1) - V_{\mathrm{c}}(t) ] \nonumber\\
    \le & 2 \max\{I, J\} J 
    - 2 a_{\min} \sum_i \expt \Bigl[ 
    \mathbb{1}_{{\cal C}}
    Q_{\mathrm{c},i}(t)
    \mathbb{1} \Bigl\{ Q_{\mathrm{c},i}(t) \ge q^{\mathrm{th}} \Bigr\}
    \Bigr]\label{equ:drift-negative}\\
    & + \sum_i |{\cal E}_{\mathrm{c},i}|^2 + 2 q^{\mathrm{th}} \sum_i \left( 2 L_{F^{-1}_i} e_{\mathrm{c},i} + 2\epsilon \left[1 + L_{F^{-1}_i} \left( p_{\mathrm{c},i,\max} - p_{\mathrm{c},i,\min}\right)\right]  + \delta \sqrt{|{\cal E}_{\mathrm{c},i}|} \right)
    + 2 I q^{\mathrm{th}} \expt \bigl[ \mathbb{1}_{{\cal C}^{\mathrm{c}}} \bigr]\nonumber\\
    & 
    + 2 \sum_i \expt \bigl[ \mathbb{1}_{{\cal C}} Q_{\mathrm{c},i}(t) \mathbb{1} \bigl\{ Q_{\mathrm{c},i}(t) \ge q^{\mathrm{th}} \bigr\} \bigr] 
    \sum_j
    \left(2 L_{G^{-1}_j} e_{\mathrm{s},j} + 2\epsilon \left[1 + L_{G^{-1}_j} \left( p_{\mathrm{s},j,\max} - p_{\mathrm{s},j,\min}\right)\right]  + \delta \sqrt{|{\cal E}_{\mathrm{s},j}|}\right)\label{equ:drift-small-positive}\\
    & + 2 I q^{\mathrm{th}} 
    \sum_j
    \left(2 L_{G^{-1}_j} e_{\mathrm{s},j} + 2\epsilon \left[1 + L_{G^{-1}_j} \left( p_{\mathrm{s},j,\max} - p_{\mathrm{s},j,\min}\right)\right]  + \delta \sqrt{|{\cal E}_{\mathrm{s},j}|}\right).\nonumber
\end{align}
For sufficiently large $T$ and $\epsilon < \delta$, by the definition of $e_{\mathrm{s},j}$, $e_{\mathrm{s},j}$, $\epsilon$, and $\delta$ will be sufficiently small such that
\[
    \sum_j
    \left(2 L_{G^{-1}_j} e_{\mathrm{s},j} + 2\epsilon \left[1 + L_{G^{-1}_j} \left( p_{\mathrm{s},j,\max} - p_{\mathrm{s},j,\min}\right)\right]  + \delta \sqrt{|{\cal E}_{\mathrm{s},j}|}\right) \le \frac{a_{\min}}{2}.
\]
Hence, the positive term \eqref{equ:drift-small-positive} can be absorbed into the negative drift in \eqref{equ:drift-negative} and thus
\begin{align}\label{equ:drift-intermediate}
    & \expt [V_{\mathrm{c}}(t+1) - V_{\mathrm{c}}(t) ] \nonumber\\
    \le & 2 \max\{I, J\} J 
    - a_{\min} \sum_i \expt \Bigl[ 
    \mathbb{1}_{{\cal C}}
    Q_{\mathrm{c},i}(t)
    \mathbb{1} \Bigl\{ Q_{\mathrm{c},i}(t) \ge q^{\mathrm{th}} \Bigr\}
    \Bigr]\nonumber\\
    & + \sum_i |{\cal E}_{\mathrm{c},i}|^2 + 2 q^{\mathrm{th}} \sum_i \left( 2 L_{F^{-1}_i} e_{\mathrm{c},i} + 2\epsilon \left[1 + L_{F^{-1}_i} \left( p_{\mathrm{c},i,\max} - p_{\mathrm{c},i,\min}\right)\right]  + \delta \sqrt{|{\cal E}_{\mathrm{c},i}|} \right)
    + 2 I q^{\mathrm{th}} \expt \bigl[ \mathbb{1}_{{\cal C}^{\mathrm{c}}} \bigr]\nonumber\\
    & + 2 I q^{\mathrm{th}} \sum_j
    \left(2 L_{G^{-1}_j} e_{\mathrm{s},j} + 2\epsilon \left[1 + L_{G^{-1}_j} \left( p_{\mathrm{s},j,\max} - p_{\mathrm{s},j,\min}\right)\right]  + \delta \sqrt{|{\cal E}_{\mathrm{s},j}|}\right).
\end{align}
For the above negative drift term $- a_{\min} \sum_i \expt [ 
\mathbb{1}_{{\cal C}} Q_{\mathrm{c},i}(t) \mathbb{1} \{ Q_{\mathrm{c},i}(t) \ge q^{\mathrm{th}} \}]$, from the indicator, we know $Q_{\mathrm{c},i}(t) \ge q^{\mathrm{th}}$. Hence, we have
\begin{align*}
    & \expt [V_{\mathrm{c}}(t+1) - V_{\mathrm{c}}(t) ] \nonumber\\
    \le & 2 \max\{I, J\} J 
    - a_{\min} q^{\mathrm{th}} \sum_i \expt \Bigl[ 
    \mathbb{1}_{{\cal C}}
    \mathbb{1} \Bigl\{ Q_{\mathrm{c},i}(t) \ge q^{\mathrm{th}} \Bigr\}
    \Bigr]\nonumber\\
    & + \sum_i |{\cal E}_{\mathrm{c},i}|^2 + 2 q^{\mathrm{th}} \sum_i \left( 2 L_{F^{-1}_i} e_{\mathrm{c},i} + 2\epsilon \left[1 + L_{F^{-1}_i} \left( p_{\mathrm{c},i,\max} - p_{\mathrm{c},i,\min}\right)\right]  + \delta \sqrt{|{\cal E}_{\mathrm{c},i}|} \right)
    + 2 I q^{\mathrm{th}} \expt \bigl[ \mathbb{1}_{{\cal C}^{\mathrm{c}}} \bigr]\nonumber\\
    & + 2 I q^{\mathrm{th}} \sum_j
    \left(2 L_{G^{-1}_j} e_{\mathrm{s},j} + 2\epsilon \left[1 + L_{G^{-1}_j} \left( p_{\mathrm{s},j,\max} - p_{\mathrm{s},j,\min}\right)\right]  + \delta \sqrt{|{\cal E}_{\mathrm{s},j}|}\right).
\end{align*}
Then we have
\begin{align*}
    &  \sum_i \expt \Bigl[ 
    \mathbb{1}_{{\cal C}}
    \mathbb{1} \Bigl\{ Q_{\mathrm{c},i}(t) \ge q^{\mathrm{th}} \Bigr\}
    \Bigr]\nonumber\\
    \le & \frac{\expt [V_{\mathrm{c}}(t)  - V_{\mathrm{c}}(t+1)]}{q^{\mathrm{th}} a_{\min}}\nonumber\\
    & +  \frac{2\max\{I, J\} J + \sum_i |{\cal E}_{\mathrm{c},i}|^2}{q^{\mathrm{th}} a_{\min}}  
    + \frac{2 \sum_i \left( 2 L_{F^{-1}_i} e_{\mathrm{c},i} + 2\epsilon \left[1 + L_{F^{-1}_i} \left( p_{\mathrm{c},i,\max} - p_{\mathrm{c},i,\min}\right)\right]  + \delta \sqrt{|{\cal E}_{\mathrm{c},i}|} \right)}
    {a_{\min}}\nonumber\\
    & + \frac{ 2 I \expt \bigl[ \mathbb{1}_{{\cal C}^{\mathrm{c}}} \bigr] }{a_{\min}}
    + \frac{2 I \sum_j
    \left(2 L_{G^{-1}_j} e_{\mathrm{s},j} + 2\epsilon \left[1 + L_{G^{-1}_j} \left( p_{\mathrm{s},j,\max} - p_{\mathrm{s},j,\min}\right)\right]  + \delta \sqrt{|{\cal E}_{\mathrm{s},j}|}\right)}{a_{\min}},
\end{align*}
which holds for $t\ge 2MN+1$.
Taking summation over $t$ from $2MN+1$ to $T$, we have
\begin{align*}
    &  \sum_{t=2MN+1}^{T} \sum_i \expt \Bigl[ 
    \mathbb{1}_{{\cal C}}
    \mathbb{1} \Bigl\{ Q_{\mathrm{c},i}(t) \ge q^{\mathrm{th}} \Bigr\}
    \Bigr]\nonumber\\
    \le & \frac{\expt [V_{\mathrm{c}}(2MN+1)]
    +T\bigl(
    2\max\{I, J\} J + \sum_i |{\cal E}_{\mathrm{c},i}|^2\bigr)}
    {q^{\mathrm{th}} a_{\min}}\nonumber\\
    & + \frac{2 T \sum_i \left( 2 L_{F^{-1}_i} e_{\mathrm{c},i} + 2\epsilon \left[1 + L_{F^{-1}_i} \left( p_{\mathrm{c},i,\max} - p_{\mathrm{c},i,\min}\right)\right]  + \delta \sqrt{|{\cal E}_{\mathrm{c},i}|} \right)}
    {a_{\min}}\nonumber\\
    & + \frac{ 2 T I \expt \bigl[ \mathbb{1}_{{\cal C}^{\mathrm{c}}} \bigr] }{a_{\min}} 
    + \frac{2 T I \sum_j
    \left(2 L_{G^{-1}_j} e_{\mathrm{s},j} + 2\epsilon \left[1 + L_{G^{-1}_j} \left( p_{\mathrm{s},j,\max} - p_{\mathrm{s},j,\min}\right)\right]  + \delta \sqrt{|{\cal E}_{\mathrm{s},j}|}\right)}{a_{\min}}.
\end{align*}
Note that $\expt [ \mathbb{1}_{{\cal C}^{\mathrm{c}}} ] = \prob ( {\cal C}^{\mathrm{c}}) \le \Theta\left(T\epsilon^{\frac{\beta}{2}+1} + \epsilon^{\frac{\beta}{2}-1} \log (1/\epsilon) \right) $ by Lemma~\ref{lemma:concentration}. Hence, we have
\begin{align}\label{equ:drift-5}
    &  \sum_{t=2MN+1}^{T} \sum_i \expt \Bigl[ 
    \mathbb{1}_{{\cal C}}
    \mathbb{1} \Bigl\{ Q_{\mathrm{c},i}(t) \ge q^{\mathrm{th}} \Bigr\}
    \Bigr]\nonumber\\
    & \le \frac{\expt [V_{\mathrm{c}}(2MN+1)]
    +T\bigl(
    2\max\{I, J\} J + \sum_i |{\cal E}_{\mathrm{c},i}|^2\bigr)}
    {q^{\mathrm{th}} a_{\min}} \nonumber\\
    & + \frac{2 T \sum_i \left( 2 L_{F^{-1}_i} e_{\mathrm{c},i} + 2\epsilon \left[1 + L_{F^{-1}_i} \left( p_{\mathrm{c},i,\max} - p_{\mathrm{c},i,\min}\right)\right]  + \delta \sqrt{|{\cal E}_{\mathrm{c},i}|} \right)
    + \Theta\left(T^2\epsilon^{\frac{\beta}{2}+1} + T\epsilon^{\frac{\beta}{2}-1} \log (1/\epsilon) \right)}
    {a_{\min}}\nonumber\\
    & + \frac{2 T I \sum_j
    \left(2 L_{G^{-1}_j} e_{\mathrm{s},j} + 2\epsilon \left[1 + L_{G^{-1}_j} \left( p_{\mathrm{s},j,\max} - p_{\mathrm{s},j,\min}\right)\right]  + \delta \sqrt{|{\cal E}_{\mathrm{s},j}|}\right)}{a_{\min}}
\end{align}
Note that we do not reject arrivals to control the queue lengths during the first outer iteration ($k=1$) and control it with threshold $q^{\mathrm{th}}$ when $k\ge 2$. Therefore, we have $V_{\mathrm{c}}(2MN+1) = \sum_i Q^2_{\mathrm{c},i}(2MN+1) \le 4 I M^2 N^2$. Substituting this inequality into \eqref{equ:drift-5}, we obtain
\begin{align*}
    &  \sum_{t=2MN+1}^{T} \sum_i \expt \Bigl[ 
    \mathbb{1}_{{\cal C}}
    \mathbb{1} \Bigl\{ Q_{\mathrm{c},i}(t) \ge q^{\mathrm{th}} \Bigr\}
    \Bigr]\nonumber\\
    & \le \frac{4 I M^2 N^2
    +T\bigl(
    2\max\{I, J\} J + \sum_i |{\cal E}_{\mathrm{c},i}|^2\bigr)}
    {q^{\mathrm{th}} a_{\min}}
    \nonumber\\
    & + \frac{2 T \sum_i \left( 2 L_{F^{-1}_i} e_{\mathrm{c},i} + 2\epsilon \left[1 + L_{F^{-1}_i} \left( p_{\mathrm{c},i,\max} - p_{\mathrm{c},i,\min}\right)\right]  + \delta \sqrt{|{\cal E}_{\mathrm{c},i}|} \right)
    + \Theta\left(T^2\epsilon^{\frac{\beta}{2}+1} + T\epsilon^{\frac{\beta}{2}-1} \log (1/\epsilon) \right)}
    {a_{\min}}\nonumber\\
    & + \frac{2 T I \sum_j
    \left(2 L_{G^{-1}_j} e_{\mathrm{s},j} + 2\epsilon \left[1 + L_{G^{-1}_j} \left( p_{\mathrm{s},j,\max} - p_{\mathrm{s},j,\min}\right)\right]  + \delta \sqrt{|{\cal E}_{\mathrm{s},j}|}\right)}{a_{\min}}
\end{align*}
Since $M=\left \lceil \log_2 (1/\epsilon) \right \rceil$ and $N= \left \lceil \frac{\beta\ln (1/\epsilon)}{\epsilon^2} \right \rceil$, we have
\begin{align*}
    &  \sum_{t=2MN+1}^{T} \sum_i \expt \Bigl[ 
    \mathbb{1}_{{\cal C}}
    \mathbb{1} \Bigl\{ Q_{\mathrm{c},i}(t) \ge q^{\mathrm{th}} \Bigr\}
    \Bigr]\nonumber\\
    \le & \Theta \left(\frac{\log^4(1/\epsilon)}{\epsilon^4 q^{\mathrm{th}}}
    + \frac{T}{q^{\mathrm{th}}}
    + T^2\epsilon^{\frac{\beta}{2}+1} + T\epsilon^{\frac{\beta}{2}-1} \log (1/\epsilon) 
    + T \left(\frac{\eta \epsilon}{\delta} + \eta + \delta + \epsilon\right) \right).
\end{align*}
Lemma~\ref{lemma:bound-time-large-queue} is proved for the customer side.

For the server side, we similarly define a Lyapunov function $V_{\mathrm{s}}(t)\coloneqq \sum_j Q_{\mathrm{s},j}^2(t)$. Using the same Lyapunov drift method, we can obtain
\begin{align*}
    &  \sum_{t=2MN+1}^{T} \sum_j \expt \Bigl[ 
    \mathbb{1}_{{\cal C}}
    \mathbb{1} \Bigl\{ Q_{\mathrm{s},j}(t) \ge q^{\mathrm{th}} \Bigr\}
    \Bigr]\nonumber\\
    \le & 
    \Theta \left(\frac{\log^4(1/\epsilon)}{\epsilon^4 q^{\mathrm{th}}}
    + \frac{T}{q^{\mathrm{th}}}
    + T^2\epsilon^{\frac{\beta}{2}+1} + T\epsilon^{\frac{\beta}{2}-1} \log (1/\epsilon) 
    + T \left(\frac{\eta \epsilon}{\delta} + \eta + \delta + \epsilon\right) \right).
\end{align*}
Lemma~\ref{lemma:bound-time-large-queue} is proved for the server side.

\section{Proof of Lemma \ref{lemma:optimization}}
\label{app:proof-lemma-optimization}

Recall that $f$ is concave and $\boldsymbol{x}^*\in \argmax_{\boldsymbol{x}\in {\cal D}} f(\boldsymbol{x})$. Let $\boldsymbol{x}^{*\prime}$ be a vector in ${\cal D}'$ such that
\begin{align*}
    f(\boldsymbol{x}^{*\prime}) = \max_{\boldsymbol{x}\in {\cal D}' } f(\boldsymbol{x}).
\end{align*}
Define a function $\tilde{f}:{\cal D}'\rightarrow \mathbb{R}$ such that
\begin{align*}
    \tilde{f} (\boldsymbol{x}) \coloneqq \expt_{\boldsymbol{v}} [ f(\boldsymbol{x} + \delta \boldsymbol{v}) ],
\end{align*}
where $\boldsymbol{v}$ is a random vector uniformly distributed in the unit ball $\mathbb{B}$.
Note that $\tilde{f}$ is also concave. To see this, for any $w\in (0,1)$ and any $\boldsymbol{x}, \boldsymbol{y}\in {\cal D}'$,
\begin{align*}
    \tilde{f} (w \boldsymbol{x} + (1-w) \boldsymbol{y}) 
    = & \expt_{\boldsymbol{v}} [f (w \boldsymbol{x} + (1-w) \boldsymbol{y} + \delta \boldsymbol{v})]\nonumber\\
    = & \expt_{\boldsymbol{v}} [f (w (\boldsymbol{x} + \delta \boldsymbol{v}) + (1-w) (\boldsymbol{y} + \delta \boldsymbol{v}))]\nonumber\\
    \ge & \expt_{\boldsymbol{v}} [ w f (\boldsymbol{x} + \delta \boldsymbol{v}) + (1-w) f (\boldsymbol{y} + \delta \boldsymbol{v})]\nonumber\\
    = & w \expt_{\boldsymbol{v}} [f (\boldsymbol{x} + \delta \boldsymbol{v})]
    + (1-w)\expt_{\boldsymbol{v}} [f (\boldsymbol{y} + \delta \boldsymbol{v})]\nonumber\\
    = & w \tilde{f} (\boldsymbol{x}) + (1-w) \tilde{f} (\boldsymbol{y}),
\end{align*}
where the inequality is by the concavity of the function $f$.
Let $\tilde{\boldsymbol{x}}^*$ be a vector in ${\cal D}'$ such that
\begin{align*}
    \tilde{f} (\tilde{\boldsymbol{x}}^*) = \max_{\boldsymbol{x}\in {\cal D}'} \tilde{f}(\boldsymbol{x}).
\end{align*}
Then we have
\begin{align}\label{equ:transform-subopti-1}
    f(\boldsymbol{x}^*) - f(\boldsymbol{x}(k))
    = & \tilde{f}(\tilde{\boldsymbol{x}}^*) - \tilde{f} (\boldsymbol{x}(k))  
    + \left[ f(\boldsymbol{x}^*)  - f(\boldsymbol{x}^{*\prime}) \right] 
    + \left[ f(\boldsymbol{x}^{*\prime}) - \tilde{f}(\boldsymbol{x}^{*\prime})\right]\nonumber\\
    & + \left[\tilde{f}(\boldsymbol{x}^{*\prime}) - \tilde{f}(\tilde{\boldsymbol{x}}^*)\right]
    + \left[\tilde{f} (\boldsymbol{x}(k)) - f(\boldsymbol{x}(k))\right].
\end{align}
Since $\tilde{\boldsymbol{x}}^*$ maximizes the function $\tilde{f}$ in ${\cal D}'$ and $\boldsymbol{x}^{*\prime}\in {\cal D}'$, we have $\tilde{f}(\boldsymbol{x}^{*\prime})\le \tilde{f}(\tilde{\boldsymbol{x}}^*)$. Hence, \eqref{equ:transform-subopti-1} can be bounded by
\begin{align}\label{equ:transform-subopti-2}
    f(\boldsymbol{x}^*) - f(\boldsymbol{x}(k))
    \le \tilde{f}(\tilde{\boldsymbol{x}}^*) - \tilde{f} (\boldsymbol{x}(k))  
    + \left[ f(\boldsymbol{x}^*)  - f(\boldsymbol{x}^{*\prime}) \right] 
    + \left[ f(\boldsymbol{x}^{*\prime}) - \tilde{f}(\boldsymbol{x}^{*\prime})\right] 
    + \left[\tilde{f} (\boldsymbol{x}(k)) - f(\boldsymbol{x}(k))\right].
\end{align}
Note that for any $\boldsymbol{x}\in {\cal D}'$, we have
\begin{align}\label{equ:diff-f-tilde-f}
    |\tilde{f}(\boldsymbol{x}) - f(\boldsymbol{x})|  = | \expt_{\boldsymbol{v}} [ f(\boldsymbol{x} + \delta \boldsymbol{v}) -  f(\boldsymbol{x}) ] |.
\end{align}
We will show that the function $f$ is Lipschitz. Recall the definition of $f$. For any $\boldsymbol{x}, \boldsymbol{y}\in {\cal D}$, we have
\begin{align}\label{equ:lipschitz-f}
    f(\boldsymbol{x}) - f(\boldsymbol{y})
    \le & \Biggl| \sum_i \biggl(\sum_{j\in {\cal E}_{\mathrm{c},i}} x_{i,j}\biggr) F_i\biggl(\sum_{j\in {\cal E}_{\mathrm{c},i}} x_{i,j}\biggr) 
    - \sum_j \biggl( \sum_{i\in {\cal E}_{\mathrm{s},j}} x_{i,j} \biggr) G_j\biggl( \sum_{i\in {\cal E}_{\mathrm{s},j}} x_{i,j} \biggr)\nonumber\\
     & - \Biggl[ \sum_i \biggl(\sum_{j\in {\cal E}_{\mathrm{c},i}} y_{i,j}\biggr) F_i\biggl(\sum_{j\in {\cal E}_{\mathrm{c},i}} y_{i,j}\biggr) 
    - \sum_j \biggl( \sum_{i\in {\cal E}_{\mathrm{s},j}} y_{i,j} \biggr) G_j\biggl( \sum_{i\in {\cal E}_{\mathrm{s},j}} y_{i,j} \biggr) \Biggr]\Biggr|\nonumber\\
    \le & \sum_i (L_{F_i} + p_{\mathrm{c},i,\max})\Biggl|\sum_{j\in {\cal E}_{\mathrm{c},i}} x_{i,j} - \sum_{j\in {\cal E}_{\mathrm{c},i}} y_{i,j}\Biggr| 
    + \sum_j (L_{G_j} + p_{\mathrm{s},j,\max})\Biggl|\sum_{i\in {\cal E}_{\mathrm{s},j}} x_{i,j} - \sum_{i\in {\cal E}_{\mathrm{s},j}} y_{i,j}\Biggr|\nonumber\\
    \le & \sum_i (L_{F_i} + p_{\mathrm{c},i,\max})\sum_{j\in {\cal E}_{\mathrm{c},i}} \Bigl|x_{i,j} -  y_{i,j}\Bigr| 
    + \sum_j (L_{G_j} + p_{\mathrm{s},j,\max})\sum_{i\in {\cal E}_{\mathrm{s},j}} \Bigl|x_{i,j} -  y_{i,j}\Bigr|\nonumber\\
    \le & \Bigl[ \max_i (L_{F_i} + p_{\mathrm{c},i,\max})  + \max_j (L_{G_j} + p_{\mathrm{s},j,\max})\Bigr] \sum_{(i,j)\in {\cal E}} \Bigl|x_{i,j} -  y_{i,j}\Bigr|\nonumber\\
    \le & \Bigl[ \max_i (L_{F_i} + p_{\mathrm{c},i,\max})  + \max_j (L_{G_j} + p_{\mathrm{s},j,\max})\Bigr] \sqrt{|{\cal E}|} \| \boldsymbol{x} - \boldsymbol{y} \|_2,
\end{align}
where the second inequality is by \eqref{equ:diff-bound-f-1} and \eqref{equ:diff-bound-f-2}, the third inequality is by the triangle inequality, and the last inequality is by Cauchy-Schwarz inequality.
From\eqref{equ:diff-f-tilde-f}, \eqref{equ:lipschitz-f}, and the fact that $\|\boldsymbol{v}\|_2 \le 1$, we have
\begin{align}\label{equ:diff-f-tilde-f-final}
    |\tilde{f}(\boldsymbol{x}) - f(\boldsymbol{x})| 
    \le  \delta \Bigl[ \max_i (L_{F_i} + p_{\mathrm{c},i,\max})  + \max_j (L_{G_j} + p_{\mathrm{s},j,\max})\Bigr] \sqrt{|{\cal E}|},
\end{align}
for any $\boldsymbol{x}\in {\cal D}'$.
From \eqref{equ:transform-subopti-2} and \eqref{equ:diff-f-tilde-f-final}, we have
\begin{align}\label{equ:transform-subopti-3}
    f(\boldsymbol{x}^*) - f(\boldsymbol{x}(k))
    \le & \tilde{f}(\tilde{\boldsymbol{x}}^*) - \tilde{f} (\boldsymbol{x}(k))  
    + \left[ f(\boldsymbol{x}^*)  - f(\boldsymbol{x}^{*\prime}) \right] \nonumber\\
    & + 2 \delta \Bigl[ \max_i (L_{F_i} + p_{\mathrm{c},i,\max})  + \max_j (L_{G_j} + p_{\mathrm{s},j,\max})\Bigr] \sqrt{|{\cal E}|}.
\end{align}
Next, we will bound the term $f(\boldsymbol{x}^*)  - f(\boldsymbol{x}^{*\prime})$. 
From \eqref{equ:lipschitz-f}, for any $\boldsymbol{x}'\in {\cal D}'$, we have
\begin{align}\label{equ:diff-any-d-prime}
    f(\boldsymbol{x}^*) 
    \le & f(\boldsymbol{x}') + \Bigl[ \max_i (L_{F_i} + p_{\mathrm{c},i,\max})  + \max_j (L_{G_j} + p_{\mathrm{s},j,\max})\Bigr] \sqrt{|{\cal E}|} \|\boldsymbol{x}^* - \boldsymbol{x}'\|_2\nonumber\\
    \le & f(\boldsymbol{x}^{*\prime}) + \Bigl[ \max_i (L_{F_i} + p_{\mathrm{c},i,\max})  + \max_j (L_{G_j} + p_{\mathrm{s},j,\max})\Bigr] \sqrt{|{\cal E}|} \|\boldsymbol{x}^* - \boldsymbol{x}'\|_2,
\end{align}
where the last inequality holds since $f(\boldsymbol{x}^{*\prime}) = \max_{\boldsymbol{x}\in {\cal D}' } f(\boldsymbol{x})$.
By Lemma~\ref{lemma:shrunk-set}, we can set
\begin{align*}
    \boldsymbol{x}' = \left(1-\frac{\delta}{r}\right)(\boldsymbol{x}^* - \boldsymbol{x}_{\mathrm{ctr}}) +  \boldsymbol{x}_{\mathrm{ctr}}\in {\cal D}'.
\end{align*}
Then we have
\begin{align}\label{equ:diff-x-star-x-prime}
    \|\boldsymbol{x}^* - \boldsymbol{x}'\|_2 = \biggl\|\boldsymbol{x}^* - \biggl[\biggl(1-\frac{\delta}{r}\biggr)(\boldsymbol{x}^* - \boldsymbol{x}_{\mathrm{ctr}}) +  \boldsymbol{x}_{\mathrm{ctr}}\biggr]\biggr\|_2
    = \frac{\delta}{r} \|\boldsymbol{x}^* - \boldsymbol{x}_{\mathrm{ctr}}\|_2.
\end{align}
From \eqref{equ:diff-any-d-prime} and \eqref{equ:diff-x-star-x-prime}, we have
\begin{align}\label{equ:transform-subopti-4}
    f(\boldsymbol{x}^*)  - f(\boldsymbol{x}^{*\prime}) \le & \Bigl[ \max_i (L_{F_i} + p_{\mathrm{c},i,\max})  + \max_j (L_{G_j} + p_{\mathrm{s},j,\max})\Bigr] \sqrt{|{\cal E}|} \frac{\delta}{r} \|\boldsymbol{x}^* - \boldsymbol{x}_{\mathrm{ctr}}\|_2\nonumber\\
    \le & \frac{\delta}{r} \Bigl[ \max_i (L_{F_i} + p_{\mathrm{c},i,\max})  + \max_j (L_{G_j} + p_{\mathrm{s},j,\max})\Bigr] |{\cal E}|,
\end{align}
where the last inequality holds since $x^*_{i,j}-x_{\mathrm{ctr},i,j}\in [-1,1]$. From \eqref{equ:transform-subopti-3} and \eqref{equ:transform-subopti-4}, we have
\begin{align*}
    f(\boldsymbol{x}^*) - f(\boldsymbol{x}(k))
    \le & \tilde{f}(\tilde{\boldsymbol{x}}^*) - \tilde{f} (\boldsymbol{x}(k))  
    + \delta \Bigl( \frac{|{\cal E}|}{r} + 2 \sqrt{|{\cal E}|} \Bigr) \Bigl[ \max_i (L_{F_i} + p_{\mathrm{c},i,\max})  + \max_j (L_{G_j} + p_{\mathrm{s},j,\max})\Bigr].
\end{align*}
Hence, we have
\begin{align}\label{equ:transform-subopti-5}
    \sum_{k=1}^{K} \expt \Bigl[f(\boldsymbol{x}^*) - f(\boldsymbol{x}(k))
    \Bigl|\Bigr. {\cal C} \Bigr] 
    \le & \sum_{k=1}^{K} \expt \Bigl[
    \tilde{f}(\tilde{\boldsymbol{x}}^*) - \tilde{f} (\boldsymbol{x}(k))  
    \Bigl|\Bigr. {\cal C} \Bigr]\nonumber\\
    & + K \delta \Bigl( \frac{|{\cal E}|}{r} + 2 \sqrt{|{\cal E}|} \Bigr) \Bigl[ \max_i (L_{F_i} + p_{\mathrm{c},i,\max})  + \max_j (L_{G_j} + p_{\mathrm{s},j,\max})\Bigr].
\end{align}

Next, we will bound the term $\sum_{k=1}^{K} \expt [\tilde{f}(\tilde{\boldsymbol{x}}^*) - \tilde{f} (\boldsymbol{x}(k)) | {\cal C} ]$. We will first bound the following term:
\begin{align}\label{equ:prove-subopti-1}
    & \expt \Bigl[ \bigl\| \boldsymbol{x}(k) + \eta \hat{\boldsymbol{g}}(k) - \tilde{\boldsymbol{x}}^* \bigr\|_2^2 
    \bigl|\bigr. \boldsymbol{x}(k), {\cal C} \Bigr]\nonumber\\
    = & \expt \Bigl[ 
    \bigl\| \boldsymbol{x}(k) - \tilde{\boldsymbol{x}}^*  \bigr\|_2^2 
    + \eta^2 \bigl\| \hat{\boldsymbol{g}}(k) \bigr\|_2^2 
    + 2\eta \bigl( \boldsymbol{x}(k) - \tilde{\boldsymbol{x}}^* \bigr)^\top \hat{\boldsymbol{g}}(k)
    \bigl|\bigr. \boldsymbol{x}(k), {\cal C} \Bigr]\nonumber\\
    = & \bigl\| \boldsymbol{x}(k) - \tilde{\boldsymbol{x}}^*  \bigr\|_2^2
    + \eta^2 \expt \Bigl[   \bigl\| \hat{\boldsymbol{g}}(k) \bigr\|_2^2 \bigl|\bigr. \boldsymbol{x}(k), {\cal C} \Bigr]
    + 2\eta \bigl( \boldsymbol{x}(k) - \tilde{\boldsymbol{x}}^* \bigr)^\top
    \expt \bigl[  \hat{\boldsymbol{g}}(k) \bigl|\bigr. \boldsymbol{x}(k), {\cal C} \bigr].
\end{align}
Recall from Algorithm~\ref{alg:pricing} that
\begin{align*}
    \hat{\boldsymbol{g}}(k) = & \frac{|{\cal E}|}{2\delta} \Biggl[ \left( \sum_{i=1}^{I} \lambda_i^+(k) p_{\mathrm{c},i}^+(k, M) - \sum_{j=1}^{J} \mu_j^+(k) p_{\mathrm{s},j}^+(k, M) \right)\nonumber\\
    & - \left( \sum_{i=1}^{I} \lambda_i^-(k) p_{\mathrm{c},i}^-(k, M) - \sum_{j=1}^{J} \mu_j^-(k) p_{\mathrm{s},j}^-(k, M)
    \right) \Biggr] \boldsymbol{u}(k).
\end{align*}
Note that we have shown in the proof of Lemma~\ref{lemma:bisection-error} in Appendix~\ref{app:proof-lemma-bisection-error} that under the event ${\cal C}$,
\begin{align}\label{equ:prove-subopti-2}
    \| \hat{\boldsymbol{g}}(k) \|_2 
    \le & |{\cal E}| \left( \sum_{i=1}^{I} |{\cal E}_{\mathrm{c},i}| ( L_{F_i} + p_{\mathrm{c}, i, \max} ) + \sum_{j=1}^{J}  |{\cal E}_{\mathrm{s},j}| ( L_{G_j} + p_{\mathrm{s}, j, \max} ) \right)\nonumber\\
    & + \frac{2|{\cal E}|\epsilon}{\delta} 
     \Biggl[ \sum_{i=1}^{I} L_{F_i}
     \left(1 + L_{F^{-1}_i} \left(p_{\mathrm{c},i,\max} - p_{\mathrm{c},i,\min}\right)\right)
     + \sum_{j=1}^{J} L_{G_j}
    \left( 1 + L_{G^{-1}_j}
    \left(p_{\mathrm{s},j,\max} - p_{\mathrm{s},j,\min}\right)\right)\Biggr]\nonumber\\
    = & \Theta \left( 1 + \frac{\epsilon}{\delta} \right).
\end{align}
Hence, using \eqref{equ:prove-subopti-1} and \eqref{equ:prove-subopti-2}, we obtain
\begin{align}\label{equ:prove-subopti-3}
    & \expt \Bigl[ \bigl\| \boldsymbol{x}(k) + \eta \hat{\boldsymbol{g}}(k) - \tilde{\boldsymbol{x}}^* \bigr\|_2^2 
    \bigl|\bigr. \boldsymbol{x}(k), {\cal C} \Bigr]\nonumber\\
    \le & \bigl\| \boldsymbol{x}(k) - \tilde{\boldsymbol{x}}^*  \bigr\|_2^2 
    + \eta^2  \Theta \left( 1 + \frac{\epsilon^2}{\delta^2} + \frac{\epsilon}{\delta} \right) 
    + 2\eta \bigl( \boldsymbol{x}(k) - \tilde{\boldsymbol{x}}^* \bigr)^\top
    \expt \bigl[ \hat{\boldsymbol{g}}(k) \bigl|\bigr. \boldsymbol{x}(k), {\cal C} \bigr].
\end{align}
Define for any $\boldsymbol{x}\in{\cal D}', \boldsymbol{u}\in \{\boldsymbol{u}\in \mathbb{R}^{|{\cal E}|}: \|\boldsymbol{u}\|_2 = 1\}$,
\begin{align*}
    \boldsymbol{g}(\boldsymbol{x}, \boldsymbol{u}) \coloneqq \frac{|{\cal E}|}{2\delta} [f(\boldsymbol{x}+\delta \boldsymbol{u}) - f(\boldsymbol{x} - \delta \boldsymbol{u})] \boldsymbol{u}.
\end{align*}
Then from \eqref{equ:prove-subopti-3} we have
\begin{align}\label{equ:prove-subopti-4}
    & \expt \Bigl[ \bigl\| \boldsymbol{x}(k) + \eta \hat{\boldsymbol{g}}(k) - \tilde{\boldsymbol{x}}^* \bigr\|_2^2 
    \bigl|\bigr. \boldsymbol{x}(k), {\cal C} \Bigr]\nonumber\\
    \le & \bigl\| \boldsymbol{x}(k) - \tilde{\boldsymbol{x}}^*  \bigr\|_2^2 
    + \eta^2  \Theta \left( 1 + \frac{\epsilon^2}{\delta^2} + \frac{\epsilon}{\delta} \right) 
    + 2\eta \bigl( \boldsymbol{x}(k) - \tilde{\boldsymbol{x}}^* \bigr)^\top
    \expt \bigl[ \boldsymbol{g}(\boldsymbol{x}(k), \boldsymbol{u}(k)) \bigl|\bigr. \boldsymbol{x}(k), {\cal C} \bigr] \nonumber\\
    & +   2\eta \bigl( \boldsymbol{x}(k) - \tilde{\boldsymbol{x}}^* \bigr)^\top
    \expt \bigl[ \hat{\boldsymbol{g}}(k) - \boldsymbol{g}(\boldsymbol{x}(k), \boldsymbol{u}(k)) \bigl|\bigr. \boldsymbol{x}(k), {\cal C} \bigr].
\end{align}
By Cauchy-Schwarz inequality, we have
\begin{align}\label{equ:prove-subopti-4-1}
    & \bigl( \boldsymbol{x}(k) - \tilde{\boldsymbol{x}}^* \bigr)^\top
    \expt \bigl[ \hat{\boldsymbol{g}}(k) - \boldsymbol{g}(\boldsymbol{x}(k), \boldsymbol{u}(k)) \bigl|\bigr. \boldsymbol{x}(k), {\cal C} \bigr]\nonumber\\
    \le &  \bigl\| \boldsymbol{x}(k) - \tilde{\boldsymbol{x}}^* \bigr\|_2
    \bigl\| \expt \bigl[ \hat{\boldsymbol{g}}(k) - \boldsymbol{g}(\boldsymbol{x}(k), \boldsymbol{u}(k)) \bigl|\bigr. \boldsymbol{x}(k), {\cal C} \bigr] \bigr\|_2\nonumber\\
    \le & \bigl\| \boldsymbol{x}(k) - \tilde{\boldsymbol{x}}^* \bigr\|_2
     \expt \bigl[  \bigl\| \hat{\boldsymbol{g}}(k) - \boldsymbol{g}(\boldsymbol{x}(k), \boldsymbol{u}(k)) \bigr\|_2 \bigl|\bigr. \boldsymbol{x}(k), {\cal C} \bigr] \nonumber\\
     \le & \sqrt{|{\cal E}|}
     \expt \bigl[  \bigl\| \hat{\boldsymbol{g}}(k) - \boldsymbol{g}(\boldsymbol{x}(k), \boldsymbol{u}(k)) \bigr\|_2 \bigl|\bigr. \boldsymbol{x}(k), {\cal C} \bigr],
\end{align}
where the second inequality is by Jensen's inequality and the last inequality holds since $x_{i,j}(k) - \tilde{x}_{i,j}^*\in[-1,1]$.
Next, we will bound the difference $\| \hat{\boldsymbol{g}}(k) - \boldsymbol{g}(\boldsymbol{x}(k), \boldsymbol{u}(k)) \|_2$ under the event ${\cal C}$. Under the event ${\cal C}$, We have
\begin{align}\label{equ:diff-already-proved-1}
    & \left| \left( \sum_{i=1}^{I} \lambda_i^+(k) p_{\mathrm{c},i}^+(k, M) - \sum_{j=1}^{J} \mu_j^+(k) p_{\mathrm{s},j}^+(k, M) \right)
    - f(\boldsymbol{x}(k)+\delta \boldsymbol{u}(k)) \right| \nonumber\\
    = & \left| \left( \sum_{i=1}^{I} \lambda_i^+(k) p_{\mathrm{c},i}^+(k, M) - \sum_{j=1}^{J} \mu_j^+(k) p_{\mathrm{s},j}^+(k, M) \right)
    - f(\boldsymbol{x}^+(k)) \right|\nonumber\\
    = & \Biggl| \left( \sum_{i=1}^{I} \lambda_i^+(k) p_{\mathrm{c},i}^+(k, M)
     - \sum_{j=1}^{J} \mu_j^+(k) p_{\mathrm{s},j}^+(k, M) \right) 
     - \left( \sum_{i=1}^{I} \lambda_i^+(k) F_i(\lambda_i^+(k))
     - \sum_{j=1}^{J} \mu_j^+(k) G_j(\mu_j^+(k))
    \biggr) \right|.
\end{align}
Similarly, we have
\begin{align}\label{equ:diff-already-proved-2}
    & \left| \left( \sum_{i=1}^{I} \lambda_i^-(k) p_{\mathrm{c},i}^-(k, M) - \sum_{j=1}^{J} \mu_j^-(k) p_{\mathrm{s},j}^-(k, M) \right)
    - f(\boldsymbol{x}(k)-\delta \boldsymbol{u}(k)) \right| \nonumber\\
    = & \Biggl| \biggl( \sum_{i=1}^{I} \lambda_i^-(k) p_{\mathrm{c},i}^-(k, M)
     - \sum_{j=1}^{J} \mu_j^-(k) p_{\mathrm{s},j}^-(k, M) \biggr) 
     - \biggl( \sum_{i=1}^{I} \lambda_i^-(k) F_i(\lambda_i^-(k))
     - \sum_{j=1}^{J} \mu_j^-(k) G_j(\mu_j^-(k))
    \biggr) \Biggr|.
\end{align}
In the proof of Lemma~\ref{lemma:bisection-error} in Appendix~\ref{app:proof-lemma-bisection-error}, we have shown that
\begin{align*}
    \eqref{equ:diff-already-proved-1} \le & 2\epsilon \Biggl[ \sum_{i=1}^{I} L_{F_i}
     \left(1 + L_{F^{-1}_i} \left(p_{\mathrm{c},i,\max} - p_{\mathrm{c},i,\min}\right)\right) 
     + \sum_{j=1}^{J} L_{G_j}
    \left( 1 + L_{G^{-1}_j}
    \left(p_{\mathrm{s},j,\max} - p_{\mathrm{s},j,\min}\right)\right)\Biggr],\nonumber\\
    \eqref{equ:diff-already-proved-2} \le & 2\epsilon \Biggl[ \sum_{i=1}^{I} L_{F_i}
     \left(1 + L_{F^{-1}_i} \left(p_{\mathrm{c},i,\max} - p_{\mathrm{c},i,\min}\right)\right)
     + \sum_{j=1}^{J} L_{G_j}
    \left( 1 + L_{G^{-1}_j}
    \left(p_{\mathrm{s},j,\max} - p_{\mathrm{s},j,\min}\right)\right)\Biggr],
\end{align*}
i.e., \eqref{equ:error-gradient-1} and \eqref{equ:error-gradient-2}. Using the above inequalities, we obtain
\begin{align}\label{equ:prove-subopti-4-2}
    & \| \hat{\boldsymbol{g}}(k) - \boldsymbol{g}(\boldsymbol{x}(k), \boldsymbol{u}(k)) \|_2 \le \frac{|{\cal E}|}{2\delta} (\eqref{equ:diff-already-proved-1} + \eqref{equ:diff-already-proved-2}) \nonumber\\
    \le & \frac{2 \epsilon |{\cal E}|}{\delta}
    \Biggl[ \sum_{i=1}^{I} L_{F_i}
     \left(1 + L_{F^{-1}_i} \left(p_{\mathrm{c},i,\max} - p_{\mathrm{c},i,\min}\right)\right) 
     + \sum_{j=1}^{J} L_{G_j}
    \left( 1 + L_{G^{-1}_j}
    \left(p_{\mathrm{s},j,\max} - p_{\mathrm{s},j,\min}\right)\right)\Biggr]\nonumber\\
    = & \Theta\left(\frac{\epsilon}{\delta}\right).
\end{align}
Substituting \eqref{equ:prove-subopti-4-2} into \eqref{equ:prove-subopti-4-1}, we have
\begin{align}\label{equ:prove-subopti-4-3}
    \bigl( \boldsymbol{x}(k) - \tilde{\boldsymbol{x}}^* \bigr)^\top
    \expt \bigl[ \hat{\boldsymbol{g}}(k) - \boldsymbol{g}(\boldsymbol{x}(k), \boldsymbol{u}(k)) \bigl|\bigr. \boldsymbol{x}(k), {\cal C} \bigr]
    \le \Theta\left(\frac{\epsilon}{\delta}\right).
\end{align}
Substituting \eqref{equ:prove-subopti-4-3} into \eqref{equ:prove-subopti-4}, we have
\begin{align}\label{equ:prove-subopti-5}
    & \expt \Bigl[ \bigl\| \boldsymbol{x}(k) + \eta \hat{\boldsymbol{g}}(k) - \tilde{\boldsymbol{x}}^* \bigr\|_2^2 
    \bigl|\bigr. \boldsymbol{x}(k), {\cal C} \Bigr]\nonumber\\
    \le & \bigl\| \boldsymbol{x}(k) - \tilde{\boldsymbol{x}}^*  \bigr\|_2^2 
    + \eta^2  \Theta \left( 1 + \frac{\epsilon^2}{\delta^2} + \frac{\epsilon}{\delta} \right) 
    + 2\eta \bigl( \boldsymbol{x}(k) - \tilde{\boldsymbol{x}}^* \bigr)^\top
    \expt \bigl[ \boldsymbol{g}(\boldsymbol{x}(k), \boldsymbol{u}(k)) \bigl|\bigr. \boldsymbol{x}(k), {\cal C} \bigr]
    +   2\eta \Theta\left(\frac{\epsilon}{\delta}\right)\nonumber\\
    \le & \bigl\| \boldsymbol{x}(k) - \tilde{\boldsymbol{x}}^*  \bigr\|_2^2 
    + 2\eta \bigl( \boldsymbol{x}(k) - \tilde{\boldsymbol{x}}^* \bigr)^\top
    \expt \bigl[ \boldsymbol{g}(\boldsymbol{x}(k), \boldsymbol{u}(k)) \bigl|\bigr. \boldsymbol{x}(k), {\cal C} \bigr]
    + \Theta \left( \eta^2 + \frac{\eta^2\epsilon^2}{\delta^2}  + \frac{\eta\epsilon}{\delta} \right),
\end{align}
where the last inequality is by $\eta < 1$.
Note that the expectation in the term $\expt [ \boldsymbol{g}(\boldsymbol{x}(k), \boldsymbol{u}(k)) | \boldsymbol{x}(k), {\cal C} ]$ in \eqref{equ:prove-subopti-5} is taken over the random variable $\boldsymbol{u}(k)$. By the definition of the event ${\cal C}$ in Section~\ref{sec:roadmap-add-concentration}, we know that the randomness of ${\cal C}$ includes only the random arrivals at the time slots when there is no queue length control, i.e., $A_{\mathrm{c},i} (t^{+}_{\mathrm{c}, i}(k,m,n)), A_{\mathrm{c},i} (t^{-}_{\mathrm{c}, i}(k,m,n)), A_{\mathrm{s},j} (t^{+}_{\mathrm{s}, j}(k,m,n)), A_{\mathrm{s},j} (t^{-}_{\mathrm{s}, j}(k,m,n))$. These random arrivals can be viewed as generated before time $t=1$ and are independent of $\boldsymbol{u}(k)$. Hence, the event ${\cal C}$ is independent of $\boldsymbol{u}(k)$. Also note that $\boldsymbol{u}(k)$ is independent of $\boldsymbol{x}(k)$. Hence, the event ${\cal C}$ is independent of $\boldsymbol{u}(k)$ given $\boldsymbol{x}(k)$. Hence, we have
\begin{align}\label{equ:independent-u-c}
    \expt [ \boldsymbol{g}(\boldsymbol{x}(k), \boldsymbol{u}(k)) | \boldsymbol{x}(k), {\cal C} ] = \expt [ \boldsymbol{g}(\boldsymbol{x}(k), \boldsymbol{u}(k)) | \boldsymbol{x}(k) ].
\end{align}
It is shown in \cite{agarwal2010optimal} that 
\begin{align}\label{equ:unbiased-g}
    \expt [ \boldsymbol{g}(\boldsymbol{x}(k), \boldsymbol{u}(k)) | \boldsymbol{x}(k) ] = \nabla \tilde{f}(\boldsymbol{x}(k)).
\end{align}
Therefore, from \eqref{equ:prove-subopti-5}, \eqref{equ:independent-u-c}, and \eqref{equ:unbiased-g}, we have
\begin{align}\label{equ:prove-subopti-6}
    & \expt \Bigl[ \bigl\| \boldsymbol{x}(k) + \eta \hat{\boldsymbol{g}}(k) - \tilde{\boldsymbol{x}}^* \bigr\|_2^2 
    \bigl|\bigr. \boldsymbol{x}(k), {\cal C} \Bigr]\nonumber\\
    \le & \bigl\| \boldsymbol{x}(k) - \tilde{\boldsymbol{x}}^*  \bigr\|_2^2 
    + 2\eta \bigl( \boldsymbol{x}(k) - \tilde{\boldsymbol{x}}^* \bigr)^\top
    \nabla \tilde{f}(\boldsymbol{x}(k))
    + \Theta \left( \eta^2 + \frac{\eta^2\epsilon^2}{\delta^2} + \frac{\eta\epsilon}{\delta} \right).
\end{align}
By the concavity of the function $\tilde{f}$, we have
\begin{align}\label{equ:concav-f-tilde}
    \bigl( \boldsymbol{x}(k) - \tilde{\boldsymbol{x}}^* \bigr)^\top
    \nabla \tilde{f}(\boldsymbol{x}(k)) \le \tilde{f} (\boldsymbol{x}(k)) - \tilde{f}(\tilde{\boldsymbol{x}}^*).
\end{align}
From \eqref{equ:prove-subopti-6} and \eqref{equ:concav-f-tilde} we have
\begin{align}\label{equ:prove-subopti-7}
    & \expt \Bigl[ \bigl\| \boldsymbol{x}(k) + \eta \hat{\boldsymbol{g}}(k) - \tilde{\boldsymbol{x}}^* \bigr\|_2^2 
    \bigl|\bigr. \boldsymbol{x}(k), {\cal C} \Bigr]\nonumber\\
    \le & \bigl\| \boldsymbol{x}(k) - \tilde{\boldsymbol{x}}^*  \bigr\|_2^2 
    + 2\eta \bigl[ \tilde{f} (\boldsymbol{x}(k)) - \tilde{f}(\tilde{\boldsymbol{x}}^*) \bigr]
    + \Theta \left( \eta^2 + \frac{\eta^2\epsilon^2}{\delta^2} + \frac{\eta\epsilon}{\delta} \right).
\end{align}
Since ${\cal D}'$ is convex and $\tilde{\boldsymbol{x}}^*\in {\cal D}'$, by the nonexpansiveness of projection~\cite{bauschke10convex}, we have
\begin{align}\label{equ:diff-projection}
    \| \boldsymbol{x}(k+1) - \tilde{\boldsymbol{x}}^* \|^2_2 = \| \Pi_{{\cal D}'}(\boldsymbol{x}(k)+\eta \hat{\boldsymbol{g}}(k) ) - \Pi_{{\cal D}'} (\tilde{\boldsymbol{x}}^*) \|^2_2 \le \| \boldsymbol{x}(k)+\eta \hat{\boldsymbol{g}}(k) - \tilde{\boldsymbol{x}}^* \|_2.
\end{align}
From \eqref{equ:prove-subopti-7} and \eqref{equ:diff-projection}, we have
\begin{align*}
    \expt \Bigl[ \bigl\| \boldsymbol{x}(k+1) - \tilde{\boldsymbol{x}}^* \bigr\|_2^2 
    \bigl|\bigr. \boldsymbol{x}(k), {\cal C} \Bigr]
    \le  \bigl\| \boldsymbol{x}(k) - \tilde{\boldsymbol{x}}^*  \bigr\|_2^2 
    + 2\eta \bigl[ \tilde{f} (\boldsymbol{x}(k)) - \tilde{f}(\tilde{\boldsymbol{x}}^*) \bigr]
    + \Theta \left( \eta^2 + \frac{\eta^2\epsilon^2}{\delta^2} + \frac{\eta\epsilon}{\delta} \right).
\end{align*}
Taking expectation on both sides over the randomness of $\boldsymbol{x}(k)$, we obtain
\begin{align*}
    & \expt \Bigl[ \bigl\| \boldsymbol{x}(k+1) - \tilde{\boldsymbol{x}}^* \bigr\|_2^2 
    \bigl|\bigr. {\cal C} \Bigr]\nonumber\\
    \le & \expt \Bigl[ \bigl\| \boldsymbol{x}(k) - \tilde{\boldsymbol{x}}^*  \bigr\|_2^2 \bigl|\bigr. {\cal C} \Bigr] 
    + 2\eta \expt \bigl[ \tilde{f} (\boldsymbol{x}(k)) - \tilde{f}(\tilde{\boldsymbol{x}}^*) \bigl|\bigr. {\cal C} \bigr]
    + \Theta \left( \eta^2 + \frac{\eta^2\epsilon^2}{\delta^2} + \frac{\eta\epsilon}{\delta} \right).
\end{align*}
Rearranging terms, we obtain
\begin{align*}
    \tilde{f}(\tilde{\boldsymbol{x}}^*) - \expt \bigl[ \tilde{f} (\boldsymbol{x}(k)) \bigl|\bigr. {\cal C} \bigr] \le \frac{1}{2\eta} \expt \Bigl[ \bigl\| \boldsymbol{x}(k) - \tilde{\boldsymbol{x}}^* \bigr\|_2^2 
    - \bigl\| \boldsymbol{x}(k+1) - \tilde{\boldsymbol{x}}^* \bigr\|_2^2 
    \Bigl|\Bigr. {\cal C} \Bigr] + \Theta\left( \eta + \frac{\eta \epsilon^2}{\delta^2} + \frac{\epsilon}{\delta}\right).
\end{align*}
Taking summation from $k=1$ to $K$, we obtain
\begin{align}\label{equ:prove-subopti-final}
    \sum_{k=1}^{K} \tilde{f}(\tilde{\boldsymbol{x}}^*) - \expt \bigl[ \tilde{f} (\boldsymbol{x}(k)) \bigl|\bigr. {\cal C} \bigr] \le & \frac{1}{2\eta} \bigl\| \boldsymbol{x}(1) - \tilde{\boldsymbol{x}}^* \bigr\|_2^2
    + \Theta \left( K\eta + \frac{K\eta \epsilon^2}{\delta^2} + \frac{K\epsilon}{\delta}\right)\nonumber\\
    \le & \Theta \left( \frac{1}{\eta} +  K\eta + \frac{K\eta \epsilon^2}{\delta^2} + \frac{K\epsilon}{\delta}\right),
\end{align}
where the last inequality holds since $x_{i,j}(1) - \tilde{x}_{i,j}^*\in [-1,1]$ for all $i,j$.

Combining \eqref{equ:transform-subopti-5} and \eqref{equ:prove-subopti-final}, we obtain
\begin{align}
    \sum_{k=1}^{K} \expt \Bigl[f(\boldsymbol{x}^*) - f(\boldsymbol{x}(k))
    \Bigl|\Bigr. {\cal C} \Bigr] 
    \le & \Theta \left( \frac{1}{\eta} +  K\eta + \frac{K\eta \epsilon^2}{\delta^2} + \frac{K\epsilon}{\delta}\right)\nonumber\\
    & + K \delta \Bigl( \frac{|{\cal E}|}{r} + 2 \sqrt{|{\cal E}|} \Bigr) \Bigl[ \max_i (L_{F_i} + p_{\mathrm{c},i,\max})  + \max_j (L_{G_j} + p_{\mathrm{s},j,\max})\Bigr]\nonumber\\
    = & \Theta \left( \frac{1}{\eta} +  K\eta + \frac{K\eta \epsilon^2}{\delta^2} + \frac{K\epsilon}{\delta} + K\delta\right).
\end{align}
Lemma~\ref{lemma:optimization} is proved.

\section{Proof of Lemma~\ref{lemma:improved-bound-time-large-queue}}
\label{app:proof-lemma-improved-bound-time-large-queue}

Similar to the proof of Lemma~\ref{lemma:bound-time-large-queue}, we use Lyapunov drift method to prove Lemma~\ref{lemma:improved-bound-time-large-queue}. The main difference is that the summation starts from $t=1$ rather than from $t=2MN+1$ under the \textit{balanced pricing algorithm}. Since we have near-balanced arrival rates for all time $t$, we can obtain a negative drift.

Recall the Lyapunov function $V_{\mathrm{c}}(t)\coloneqq \sum_i Q_{\mathrm{c},i}^2(t)$. Recall that $\lambda'_i(t)\coloneqq F^{-1}_i(p^{+/-}_{\mathrm{c},i}(k(t),m(t)))$ and $\mu'_j(t)\coloneqq G^{-1}_j(p^{+/-}_{\mathrm{s},j}(k(t),m(t)))$. 
Following the same argument as the proof of Lemma~\ref{lemma:bound-time-large-queue} in Appendix~\ref{app:proof-lemma-bound-time-large-queue}, we obtain
\begin{align}\label{equ:improved-drift-3}
    & \expt [V_{\mathrm{c}}(t+1) - V_{\mathrm{c}}(t) | \boldsymbol{Q}_{\mathrm{c}}(t) ] \nonumber\\
    \le & \sum_i |{\cal E}_{\mathrm{c},i}|^2 + 2 \sum_i Q_{\mathrm{c},i}(t) \mathbb{1} \left\{ Q_{\mathrm{c},i}(t) < q^{\mathrm{th}} \right\}
     \expt \bigl[  \lambda'_i(t) \mathbb{1}_{{\cal C}} \bigl|\bigr. \boldsymbol{Q}_{\mathrm{c}}(t) \bigr]\nonumber\\
     & + 2 I q^{\mathrm{th}} \expt \bigl[ \mathbb{1}_{{\cal C}^{\mathrm{c}}} \bigl|\bigr. \boldsymbol{Q}_{\mathrm{c}}(t) \bigr]  - 2 \sum_j
    \expt \left[   \mu'_j(t) Q_{\mathrm{c},i^*_j(t)}(t) \mathbb{1}_{{\cal C}}
    \Bigl|\bigr. \boldsymbol{Q}_{\mathrm{c}}(t)  \right]
\end{align}

Next, we will relate $\lambda'_i(t)$ and $\mu'_j(t)$ to the arrival rates $\sum_{j\in {\cal E}_{\mathrm{c},i}} x_{i,j}(k(t))$ and $\sum_{i\in {\cal E}_{\mathrm{s},j}} x_{i,j}(k(t))$, where we recall that $k(t)$ denotes the outer iteration at time $t$. In fact, we can obtain a lemma similar to Lemma~\ref{lemma:bisection-error} under the \textit{balanced pricing algorithm}.
\begin{lemma}\label{lemma:improved-bisection-error}
    Let Assumption~\ref{assum:1}, Assumption~\ref{assum:4}, and Assumption~\ref{assum:5} hold. Suppose $2e_{\mathrm{c},i}\le p_{\mathrm{c},i,\max} - p_{\mathrm{c},i,\min}$ and $2e_{\mathrm{s},j}\le p_{\mathrm{s},j,\max} - p_{\mathrm{s},j,\min}$. Then under the balanced pricing algorithm and the event ${\cal C}$, for all $t\ge 1$, we have
    \begin{align*}
    \biggl|\lambda'_i(t) - \sum_{j\in {\cal E}_{\mathrm{c},i}} x_{i,j}(k(t))\biggr|
    \le & 2 L_{F^{-1}_i} e_{\mathrm{c},i} + 2\epsilon \left[1 + L_{F^{-1}_i} \left( p_{\mathrm{c},i,\max} - p_{\mathrm{c},i,\min}\right)\right]  + \delta \sqrt{|{\cal E}_{\mathrm{c},i}|}, \mbox{ for all } i;\nonumber\\
    \biggl|\mu'_j(t) - \sum_{i\in {\cal E}_{\mathrm{s},j}} x_{i,j}(k(t))\biggr|
    \le & 2 L_{G^{-1}_j} e_{\mathrm{s},j} + 2\epsilon \left[1 + L_{G^{-1}_j} \left( p_{\mathrm{s},j,\max} - p_{\mathrm{s},j,\min}\right)\right]  + \delta \sqrt{|{\cal E}_{\mathrm{s},j}|}, \mbox{ for all } j.
    \end{align*}
\end{lemma}
\noindent The difference between Lemma~\ref{lemma:improved-bisection-error} and Lemma~\ref{lemma:bisection-error} is that Lemma~\ref{lemma:improved-bisection-error} states that the arrival rates are nearly balanced for \textit{all time $t$}.
Proof of Lemma~\ref{lemma:improved-bisection-error} can be found in Appendix~\ref{app:proof-lemma-improved-bisection-error}. By Lemma~\ref{lemma:improved-bisection-error} and \eqref{equ:improved-drift-3}, following the same proof as that in the proof of Lemma~\ref{lemma:bound-time-large-queue} in Appendix~\ref{app:proof-lemma-bound-time-large-queue}, we can obtain
\begin{align}\label{equ:improved-drift-4}
    &  \sum_i \expt \Bigl[ 
    \mathbb{1}_{{\cal C}}
    \mathbb{1} \Bigl\{ Q_{\mathrm{c},i}(t) \ge q^{\mathrm{th}} \Bigr\}
    \Bigr]\nonumber\\
    \le & \frac{\expt [V_{\mathrm{c}}(t)  - V_{\mathrm{c}}(t+1)]}{q^{\mathrm{th}} a_{\min}}
    +  \frac{2\max\{I, J\} J + \sum_i |{\cal E}_{\mathrm{c},i}|^2}{q^{\mathrm{th}} a_{\min}}  
    + \frac{\Theta(\frac{\eta \epsilon}{\delta} + \eta + \delta + \epsilon)}
    {a_{\min}}
    + \frac{ 2I \expt \bigl[ \mathbb{1}_{{\cal C}^{\mathrm{c}}} \bigr] }{a_{\min}}\nonumber\\
    \le & \frac{\expt [V_{\mathrm{c}}(t)  - V_{\mathrm{c}}(t+1)]}{ q^{\mathrm{th}} a_{\min}}
    +  \frac{2\max\{I, J\} J + \sum_i |{\cal E}_{\mathrm{c},i}|^2}{q^{\mathrm{th}} a_{\min}}  
    + \frac{\Theta(\frac{\eta \epsilon}{\delta} + \eta + \delta + \epsilon)}
    {a_{\min}}
     + \frac{ \Theta\left(T\epsilon^{\frac{\beta}{2}+1} + \epsilon^{\frac{\beta}{2}-1} \log (1/\epsilon) \right) }{a_{\min}},
\end{align}
where the last inequality is by Lemma~\ref{lemma:concentration}. Note that the inequality \eqref{equ:improved-drift-4} holds for all $t\ge 1$. Taking summation over $t$ from $t=1$ to $T$ and using the assumption that all queues are empty at $t=1$, we have
\begin{align}\label{equ:improved-drift-final-1}
    \sum_{t=1}^{T} \sum_i \expt \Bigl[ 
    \mathbb{1}_{{\cal C}}
    \mathbb{1} \Bigl\{ Q_{\mathrm{c},i}(t) \ge q^{\mathrm{th}} \Bigr\}
    \Bigr]
    \le \Theta \left(
    \frac{T}{q^{\mathrm{th}}}
    + T^2\epsilon^{\frac{\beta}{2}+1} + T\epsilon^{\frac{\beta}{2}-1} \log (1/\epsilon)
    + T \left(\frac{\eta \epsilon}{\delta} + \eta + \delta + \epsilon\right)
    \right).
\end{align}

For the server side, we similarly define a Lyapunov function $V_{\mathrm{s}}(t)\coloneqq \sum_j Q_{\mathrm{s},j}^2(t)$. Using the same Lyapunov drift method, we can obtain
\begin{align}\label{equ:improved-drift-final-2}
    \sum_{t=1}^{T} \sum_j \expt \Bigl[ 
    \mathbb{1}_{{\cal C}}
    \mathbb{1} \Bigl\{ Q_{\mathrm{s},j}(t) \ge q^{\mathrm{th}} \Bigr\}
    \Bigr]
    \le \Theta \left(
    \frac{T}{q^{\mathrm{th}}}
    + T^2\epsilon^{\frac{\beta}{2}+1} + T\epsilon^{\frac{\beta}{2}-1} \log (1/\epsilon)
    + T \left(\frac{\eta \epsilon}{\delta} + \eta + \delta + \epsilon\right)
    \right).
\end{align}
From \eqref{equ:improved-drift-final-1} and \eqref{equ:improved-drift-final-2}, Lemma~\ref{lemma:improved-bound-time-large-queue} is proved.

\section{Proof of Lemma~\ref{lemma:bisection-error}}
\label{app:proof-lemma-bisection-error}

Under the event ${\cal C}$, we will first prove the following:
for all outer iteration $k$,
\begin{align}
    \left| \lambda_i^{+/-}(k) - F^{-1}_i(p^{+/-}_{\mathrm{c},i}(k,M)) \right| 
    \le & 2\epsilon \left[1 + L_{F^{-1}_i} \left( p_{\mathrm{c},i,\max} - p_{\mathrm{c},i,\min}\right)\right], \mbox{ for all } i,\label{equ:lemma-price-error-1}\\
    \left| \mu_j^{+/-}(k) - G^{-1}_j(p^{+/-}_{\mathrm{s},j}(k,M)) \right| 
    \le & 2 \epsilon \left[ 1 + L_{G^{-1}_j}
    \left( p_{\mathrm{s},j,\max} - p_{\mathrm{s},j,\min}\right)\right], \mbox{ for all } j. \label{equ:lemma-price-error-2}
\end{align}
We will prove \eqref{equ:lemma-price-error-1} and \eqref{equ:lemma-price-error-2} by induction. 

First consider the first outer iteration $k=1$. For the first outer iteration, we are running the bisection algorithm with the initial price intervals $[p_{\mathrm{c},i,\min}, p_{\mathrm{c},i,\max}]$ and $[p_{\mathrm{s},j,\min}, p_{\mathrm{s},j,\max}]$, i.e., $\underline{p}^{+/-}_{\mathrm{c},i}(k,1)=p_{\mathrm{c},i,\min}$, $\bar{p}^{+/-}_{\mathrm{c},i}(k,1)=p_{\mathrm{c},i,\max}$, $\underline{p}^{+/-}_{\mathrm{s},j}(k,1)=p_{\mathrm{s},j,\min}$, $\bar{p}^{+/-}_{\mathrm{s},j}(k,1)=p_{\mathrm{s},j,\max}$ for all $i,j$. Hence, for $k=1$, we have
\begin{align*}
    \lambda_i^{+/-}(k) \in & [0,1] \subseteq [F^{-1}_i(\bar{p}^{+/-}_{\mathrm{c},i}(k,1))-\epsilon, F^{-1}_i(\underline{p}^{+/-}_{\mathrm{c},i}(k,1)) + \epsilon], \text{for all $i$,}\nonumber\\
    \mu_j^{+/-}(k) \in & [0,1] \subseteq [G^{-1}_j(\underline{p}^{+/-}_{\mathrm{s},j}(k,1))-\epsilon, G^{-1}_j(\bar{p}^{+/-}_{\mathrm{s},j}(k,1)) + \epsilon], \text{for all $j$}.
\end{align*}
We will show by induction that for $k=1$ and all $m=1,\ldots,M$,
\begin{align}
    \lambda_i^{+/-}(k) \in & [F^{-1}_i(\bar{p}^{+/-}_{\mathrm{c},i}(k,m))-\epsilon, F^{-1}_i(\underline{p}^{+/-}_{\mathrm{c},i}(k,m)) + \epsilon], \text{for all $i$,}\label{equ:arr-in-interval-1}\\
    \mu_j^{+/-}(k) \in & [G^{-1}_j(\underline{p}^{+/-}_{\mathrm{s},j}(k,m))-\epsilon, G^{-1}_j(\bar{p}^{+/-}_{\mathrm{s},j}(k,m)) + \epsilon], \text{for all $j$}.\label{equ:arr-in-interval-2}
\end{align}
Suppose \eqref{equ:arr-in-interval-1} and \eqref{equ:arr-in-interval-2} hold for $m$. We need to show that they hold for $m+1$. For $\lambda_i^{+/-}(k)$, there are two possible cases:
\begin{itemize}[leftmargin=*]
    \item If $\hat{\lambda}_i^{+/-}(k,m) > \lambda_i^{+/-}(k)$, then by Algorithm~\ref{alg:bisection}, we have $\underline{p}^{+/-}_{\mathrm{c}, i}(k,m+1) = p^{+/-}_{\mathrm{c}, i}  (k,m)$ and $\bar{p}^{+/-}_{\mathrm{c}, i}(k,m+1) = \bar{p}^{+/-}_{\mathrm{c}, i}(k,m)$. Hence, we have
    \begin{align*}
        F_i^{-1}(\underline{p}^{+/-}_{\mathrm{c}, i}(k,m+1)) = F_i^{-1}(p^{+/-}_{\mathrm{c}, i}  (k,m)) > \hat{\lambda}_i^{+/-}(k,m) - \epsilon > \lambda_i^{+/-}(k) - \epsilon,
    \end{align*}
    where the first inequality holds since the event ${\cal C}$ holds. And by induction hypothesis, we have
    \begin{align*}
        F_i^{-1}(\bar{p}^{+/-}_{\mathrm{c}, i}(k,m+1)) = F_i^{-1}(\bar{p}^{+/-}_{\mathrm{c}, i}(k,m)) \le \lambda_i^{+/-}(k) + \epsilon.
    \end{align*}
    
    \item If $\hat{\lambda}_i^{+/-}(k,m) \le \lambda_i^{+/-}(k)$, then by Algorithm~\ref{alg:bisection}, we have $\bar{p}^{+/-}_{\mathrm{c},i}(k,m+1) = p^{+/-}_{\mathrm{c}, i}(k,m)$ and $\underline{p}^{+/-}_{\mathrm{c},i}(k,m+1) = \underline{p}^{+/-}_{\mathrm{c},i}(k,m)$. Hence, we have
    \begin{align*}
        F_i^{-1}(\bar{p}^{+/-}_{\mathrm{c},i}(k,m+1)) = F_i^{-1}(p^{+/-}_{\mathrm{c}, i}  (k,m)) < \hat{\lambda}_i^{+/-}(k,m) + \epsilon  \le \lambda_i^{+/-}(k) + \epsilon,
    \end{align*}
    where the first inequality holds since the event ${\cal C}$ holds. And by induction hypothesis, we have
    \begin{align*}
        F_i^{-1}(\underline{p}^{+/-}_{\mathrm{c},i}(k,m+1)) = F_i^{-1}(\underline{p}^{+/-}_{\mathrm{c},i}(k,m)) \ge \lambda_i^{+/-}(k) - \epsilon.
    \end{align*}
\end{itemize}
Therefore, we have $\lambda_i^{+/-}(k) \in[F^{-1}_i(\bar{p}^{+/-}_{\mathrm{c},i}(k,m+1))-\epsilon, F^{-1}_i(\underline{p}^{+/-}_{\mathrm{c},i}(k,m+1)) + \epsilon]$ for all $i$.

For $\mu_j^{+/-}(k)$, similarly, there are two possible cases:
\begin{itemize}[leftmargin=*]
    \item If $\hat{\mu}_j^{+/-}(k,m) > \mu_j^{+/-}(k)$, then by Algorithm~\ref{alg:bisection}, we have $\bar{p}^{+/-}_{\mathrm{s}, j}(k,m+1) = p^{+/-}_{\mathrm{s}, j}(k,m)$ and $\underline{p}^{+/-}_{\mathrm{s}, j}(k,m+1) = \underline{p}^{+/-}_{\mathrm{s}, j}(k,m)$. Hence, we have
    \begin{align*}
        G_j^{-1}(\bar{p}^{+/-}_{\mathrm{s}, j}(k,m+1)) = G_j^{-1}(p^{+/-}_{\mathrm{s}, j}(k,m)) > \hat{\mu}_j^{+/-}(k,m) - \epsilon > \mu_j^{+/-}(k) - \epsilon,
    \end{align*}
    where the first inequality holds since the event ${\cal C}$ holds. And by induction hypothesis, we have
    \begin{align*}
        G_j^{-1}(\underline{p}^{+/-}_{\mathrm{s}, j}(k,m+1)) = G_j^{-1}(\underline{p}^{+/-}_{\mathrm{s}, j}(k,m)) \le \mu_j^{+/-}(k) + \epsilon.
    \end{align*}
    
    \item If $\hat{\mu}_j^{+/-}(k,m) \le \mu_j^{+/-}(k)$, then by Algorithm~\ref{alg:bisection}, we have $\underline{p}^{+/-}_{\mathrm{s}, j}(k,m+1) = p^{+/-}_{\mathrm{s}, j}(k,m)$ and $\bar{p}^{+/-}_{\mathrm{s}, j}(k,m+1) = \bar{p}^{+/-}_{\mathrm{s}, j}(k,m)$. Hence, we have
    \begin{align*}
        G_j^{-1}(\underline{p}^{+/-}_{\mathrm{s}, j}(k,m+1)) = G_j^{-1}(p^{+/-}_{\mathrm{s}, j}(k,m)) < \hat{\mu}_j^{+/-}(k,m) + \epsilon  \le \mu_j^{+/-}(k) + \epsilon,
    \end{align*}
    where the first inequality holds since the event ${\cal C}$ holds. And by induction hypothesis, we have
    \begin{align*}
        G_j^{-1}(\bar{p}^{+/-}_{\mathrm{s}, j}(k,m+1)) = G_j^{-1}(\bar{p}^{+/-}_{\mathrm{s}, j}(k,m)) \ge \mu_j^{+/-}(k) - \epsilon.
    \end{align*}
\end{itemize}
Therefore, we have $\mu_j^{+/-}(k) \in[G^{-1}_j(\underline{p}^{+/-}_{\mathrm{s},j}(k,m+1))-\epsilon, G^{-1}_j(\bar{p}^{+/-}_{\mathrm{s},j}(k,m+1)) + \epsilon]$ for all $i$. 

Therefore, we have proved that \eqref{equ:arr-in-interval-1} and \eqref{equ:arr-in-interval-2} hold for $k=1$ and all $m$.

Since in each bisection iteration, the algorithm divides the search interval by $2$, we have, for $k=1$ and $m=M$,
\begin{align*}
    \bar{p}^{+/-}_{\mathrm{c},i}(1,M) - \underline{p}^{+/-}_{\mathrm{c},i}(1,M)
    \le & \frac{p_{\mathrm{c},i,\max} - p_{\mathrm{c},i,\min}}{2^{M-1}},\nonumber\\
    \bar{p}^{+/-}_{\mathrm{s},j}(1,M) - \underline{p}^{+/-}_{\mathrm{s},j}(1,M)
    \le & \frac{p_{\mathrm{s},j,\max} - p_{\mathrm{s},j,\min}}{2^{M-1}}.
\end{align*}
By Assumption~\ref{assum:1}, $F_i^{-1}$ and $G_j^{-1}$ are Lipschitz. Hence, we have 
\begin{align}
    F^{-1}_i(\underline{p}^{+/-}_{\mathrm{c},i}(1,M)) - F^{-1}_i(\bar{p}^{+/-}_{\mathrm{c},i}(1,M)) 
    \le & L_{F^{-1}_i} \frac{p_{\mathrm{c},i,\max} - p_{\mathrm{c},i,\min}}{2^{M-1}},\label{equ:diff-arrs-1}\\
    G^{-1}_j(\bar{p}^{+/-}_{\mathrm{s},j}(1,M)) - G^{-1}_j(\underline{p}^{+/-}_{\mathrm{s},j}(1,M)) 
    \le & L_{G^{-1}_j}
    \frac{p_{\mathrm{s},j,\max} - p_{\mathrm{s},j,\min}}{2^{M-1}}.\label{equ:diff-arrs-2}
\end{align}
Note that $p^{+/-}_{\mathrm{c},i}(1,M) \in [\underline{p}^{+/-}_{\mathrm{c},i}(1,M), \bar{p}^{+/-}_{\mathrm{c},i}(1,M)]$ and $p^{+/-}_{\mathrm{s},j}(1,M) \in [\underline{p}^{+/-}_{\mathrm{s},j}(1,M), \bar{p}^{+/-}_{\mathrm{s},j}(1,M)]$. Hence, we have
\begin{align}
    F^{-1}_i(p^{+/-}_{\mathrm{c},i}(1,M)) 
    \in &  [F^{-1}_i(\bar{p}^{+/-}_{\mathrm{c},i}(1,M))-\epsilon, F^{-1}_i(\underline{p}^{+/-}_{\mathrm{c},i}(1,M))+\epsilon],\label{equ:arr-in-interval-3}\\
    G^{-1}_j (p^{+/-}_{\mathrm{s},j}(1,M)) 
    \in & [G^{-1}_j(\underline{p}^{+/-}_{\mathrm{s},j}(1,M))-\epsilon, G^{-1}_j(\bar{p}^{+/-}_{\mathrm{s},j}(1,M))+\epsilon].\label{equ:arr-in-interval-4}
\end{align}
From \eqref{equ:arr-in-interval-1}, \eqref{equ:arr-in-interval-2}, \eqref{equ:arr-in-interval-3}, and \eqref{equ:arr-in-interval-4}, we have
\begin{align*}
    \left| \lambda_i^{+/-}(1) - F^{-1}_i(p^{+/-}_{\mathrm{c},i}(1,M)) \right| 
    \le & F^{-1}_i(\underline{p}^{+/-}_{\mathrm{c},i}(1,M)) - F^{-1}_i(\bar{p}^{+/-}_{\mathrm{c},i}(1,M)) + 2\epsilon,\\
    \left| \mu_j^{+/-}(1) - G^{-1}_j(p^{+/-}_{\mathrm{s},j}(1,M)) \right| 
    \le & G^{-1}_j(\bar{p}^{+/-}_{\mathrm{s},j}(1,M)) - G^{-1}_j(\underline{p}^{+/-}_{\mathrm{s},j}(1,M)) + 2\epsilon.
\end{align*}
Substituting \eqref{equ:diff-arrs-1} and \eqref{equ:diff-arrs-2} into the above inequalities, we have
\begin{align*}
    \left| \lambda_i^{+/-}(1) - F^{-1}_i(p^{+/-}_{\mathrm{c},i}(1,M)) \right| 
    \le & L_{F^{-1}_i} \frac{p_{\mathrm{c},i,\max} - p_{\mathrm{c},i,\min}}{2^{M-1}} + 2\epsilon,\\
    \left| \mu_j^{+/-}(1) - G^{-1}_j(p^{+/-}_{\mathrm{s},j}(1,M)) \right| 
    \le & L_{G^{-1}_j}
    \frac{p_{\mathrm{s},j,\max} - p_{\mathrm{s},j,\min}}{2^{M-1}} + 2\epsilon.
\end{align*}
Recall that $M=\log_2 \frac{1}{\epsilon}$. Then we have
\begin{align*}
    \left| \lambda_i^{+/-}(1) - F^{-1}_i(p^{+/-}_{\mathrm{c},i}(1,M)) \right| 
    \le & 2\epsilon \left[1 + L_{F^{-1}_i} \left(p_{\mathrm{c},i,\max} - p_{\mathrm{c},i,\min}\right)\right],\\
    \left| \mu_j^{+/-}(1) - G^{-1}_j(p^{+/-}_{\mathrm{s},j}(1,M)) \right| 
    \le & 2 \epsilon \left[ 1 + L_{G^{-1}_j}
    \left(p_{\mathrm{s},j,\max} - p_{\mathrm{s},j,\min}\right)\right],
\end{align*}
i.e., \eqref{equ:lemma-price-error-1} and \eqref{equ:lemma-price-error-2} hold for $k=1$.

Next, we will use induction on $k$ to show that \eqref{equ:arr-in-interval-1} and \eqref{equ:arr-in-interval-2} hold for all $k$ and all $m$ and 
\begin{align}
    \left| \lambda_i^{+/-}(k) - F^{-1}_i(p^{+/-}_{\mathrm{c},i}(k,M)) \right| 
    \le & 2\epsilon \left[1 + L_{F^{-1}_i} \left(p_{\mathrm{c},i,\max} - p_{\mathrm{c},i,\min}\right)\right],\label{equ:price-error-1}\\
    \left| \mu_j^{+/-}(k) - G^{-1}_j(p^{+/-}_{\mathrm{s},j}(k,M)) \right| 
    \le & 2 \epsilon \left[ 1 + L_{G^{-1}_j}
    \left(p_{\mathrm{s},j,\max} - p_{\mathrm{s},j,\min}\right)\right].\label{equ:price-error-2}
\end{align}
for all $k$. Note that we have already shown that \eqref{equ:arr-in-interval-1}, \eqref{equ:arr-in-interval-2}, \eqref{equ:price-error-1}, and \eqref{equ:price-error-2} hold for $k=1$. 
Suppose \eqref{equ:arr-in-interval-1}-\eqref{equ:arr-in-interval-2} hold for $k$ and all $m$, and \eqref{equ:price-error-1}-\eqref{equ:price-error-2} hold for $k$. We need to show that \eqref{equ:arr-in-interval-1}-\eqref{equ:arr-in-interval-2} hold for $k+1$ and all $m$, and \eqref{equ:price-error-1}-\eqref{equ:price-error-2} hold for $k+1$. We will first show that \eqref{equ:arr-in-interval-1}-\eqref{equ:arr-in-interval-2} hold for $k+1$ and $m=1$, i.e.,
\begin{align}
    \lambda_i^{+/-}(k+1) \in & [F^{-1}_i(\bar{p}^{+/-}_{\mathrm{c},i}(k+1,1))-\epsilon, F^{-1}_i(\underline{p}^{+/-}_{\mathrm{c},i}(k+1,1)) + \epsilon], \text{for all $i$,}\label{equ:next-arr-in-interval-1}\\
    \mu_j^{+/-}(k+1) \in & [G^{-1}_j(\underline{p}^{+/-}_{\mathrm{s},j}(k+1,1))-\epsilon, G^{-1}_j(\bar{p}^{+/-}_{\mathrm{s},j}(k+1,1)) + \epsilon], \text{for all $j$}.\label{equ:next-arr-in-interval-2}
\end{align}
By Algorithm~\ref{alg:pricing}, we have $\underline{p}^{+/-}_{\mathrm{c},i}(k+1,1) = p^{+/-}_{\mathrm{c},i}(k,M) - e_{\mathrm{c},i}$, $\bar{p}^{+/-}_{\mathrm{c},i}(k+1,1) = p^{+/-}_{\mathrm{c},i}(k,M) + e_{\mathrm{c},i}$, $\underline{p}^{+/-}_{\mathrm{s},j}(k+1,1) = p^{+/-}_{\mathrm{s},j}(k,M) - e_{\mathrm{s},j}$, and $\bar{p}^{+/-}_{\mathrm{s},j}(k+1,1) = p^{+/-}_{\mathrm{s},j}(k,M) + e_{\mathrm{s},j}$. Therefore, proving \eqref{equ:next-arr-in-interval-1} and \eqref{equ:next-arr-in-interval-2} is equivalent to proving the following:
\begin{align}
    \lambda_i^{+/-}(k+1) \in & [F^{-1}_i(p^{+/-}_{\mathrm{c},i}(k,M) + e_{\mathrm{c},i}) - \epsilon, F^{-1}_i(p^{+/-}_{\mathrm{c},i}(k,M) - e_{\mathrm{c},i}) + \epsilon], \text{for all $i$,}\label{equ:next-arr-in-interval-3}\\
    \mu_j^{+/-}(k+1) \in & [G^{-1}_j(p^{+/-}_{\mathrm{s},j}(k,M) - e_{\mathrm{s},j})-\epsilon, G^{-1}_j(p^{+/-}_{\mathrm{s},j}(k,M) + e_{\mathrm{s},j}) + \epsilon], \text{for all $j$}\label{equ:next-arr-in-interval-4}.
\end{align}
Note that by Algorithm~\ref{alg:pricing}, we have
\begin{align*}
    \| \boldsymbol{x}(k+1) - \boldsymbol{x}(k) \|_2 = &  \| \Pi_{{\cal D}'}(\boldsymbol{x}(k)+\eta \hat{\boldsymbol{g}}(k) ) - \boldsymbol{x}(k) \|_2\nonumber\\
    = &  \| \Pi_{{\cal D}'}(\boldsymbol{x}(k)+\eta \hat{\boldsymbol{g}}(k) ) - \Pi_{{\cal D}'} (\boldsymbol{x}(k)) \|_2,
\end{align*}
where the second equality holds since $\boldsymbol{x}(k) \in {\cal D}'$. Note that the constraints in ${\cal D}'$ are affine and hence ${\cal D}'$ is convex. Since ${\cal D}'$ is convex, by the nonexpansiveness of projection~\cite{bauschke10convex}, we have
\begin{align}\label{equ:diff-x-bound}
    \| \boldsymbol{x}(k+1) - \boldsymbol{x}(k) \|_2 = \| \Pi_{{\cal D}'}(\boldsymbol{x}(k)+\eta \hat{\boldsymbol{g}}(k) ) - \Pi_{{\cal D}'} (\boldsymbol{x}(k)) \|_2 \le \| \boldsymbol{x}(k)+\eta \hat{\boldsymbol{g}}(k) - \boldsymbol{x}(k)  \|_2 = \eta \| \hat{\boldsymbol{g}}(k) \|_2,
\end{align}
where
\begin{align}
    \| \hat{\boldsymbol{g}}(k) \|_2 = & \Biggl\| \frac{|{\cal E}|}{2\delta} \biggl[ \biggl( \sum_{i=1}^{I} \lambda_i^+(k) p_{\mathrm{c},i}^+(k, M)
     - \sum_{j=1}^{J} \mu_j^+(k) p_{\mathrm{s},j}^+(k, M) \biggr) \nonumber\\
     & \qquad - \biggl( \sum_{i=1}^{I} \lambda_i^-(k) p_{\mathrm{c},i}^-(k, M) - \sum_{j=1}^{J} \mu_j^-(k) p_{\mathrm{s},j}^-(k, M)
    \biggr) \biggr] \boldsymbol{u}(k) \Biggr\|_2\nonumber\\
    = & \frac{|{\cal E}|}{2\delta} \Biggl| \biggl( \sum_{i=1}^{I} \lambda_i^+(k) p_{\mathrm{c},i}^+(k, M)
     - \sum_{j=1}^{J} \mu_j^+(k) p_{\mathrm{s},j}^+(k, M) \biggr) \nonumber\\
     & \qquad - \biggl( \sum_{i=1}^{I} \lambda_i^-(k) p_{\mathrm{c},i}^-(k, M) - \sum_{j=1}^{J} \mu_j^-(k) p_{\mathrm{s},j}^-(k, M)
    \biggr) \Biggr|.\label{equ:norm-g-hat-k}
\end{align}
Note that we have
\begin{align}\label{equ:error-gradient-1}
    & \Biggl| \biggl( \sum_{i=1}^{I} \lambda_i^+(k) p_{\mathrm{c},i}^+(k, M)
     - \sum_{j=1}^{J} \mu_j^+(k) p_{\mathrm{s},j}^+(k, M) \biggr) 
     - \biggl( \sum_{i=1}^{I} \lambda_i^+(k) F_i(\lambda_i^+(k))
     - \sum_{j=1}^{J} \mu_j^+(k) G_j(\mu_j^+(k))
    \biggr) \Biggr|\nonumber\\
    \le & \sum_{i=1}^{I} \lambda_i^+(k) 
    \left| p_{\mathrm{c},i}^+(k, M) - F_i(\lambda_i^+(k)) \right|
    + \sum_{j=1}^{J} \mu_j^+(k)
    \left| p_{\mathrm{s},j}^+(k, M) - G_j(\mu_j^+(k)) \right|\nonumber\\
    \le & \sum_{i=1}^{I}
    \left| p_{\mathrm{c},i}^+(k, M) - F_i(\lambda_i^+(k)) \right|
    + \sum_{j=1}^{J}
    \left| p_{\mathrm{s},j}^+(k, M) - G_j(\mu_j^+(k)) \right|\nonumber\\
    \le & \sum_{i=1}^{I} L_{F_i}
    \left| F_i^{-1} (p_{\mathrm{c},i}^+(k, M)) - \lambda_i^+(k) \right|
    + \sum_{j=1}^{J} L_{G_j}
    \left| G_j^{-1} (p_{\mathrm{s},j}^+(k, M)) - \mu_j^+(k) \right|\nonumber\\
    \le & 2\epsilon \Biggl[ \sum_{i=1}^{I} L_{F_i}
     \left(1 + L_{F^{-1}_i} \left(p_{\mathrm{c},i,\max} - p_{\mathrm{c},i,\min}\right)\right) 
     + \sum_{j=1}^{J} L_{G_j}
    \left( 1 + L_{G^{-1}_j}
    \left(p_{\mathrm{s},j,\max} - p_{\mathrm{s},j,\min}\right)\right)\Biggr],
\end{align}
where the first inequality is by triangle inequality, the second inequality is due to the fact that $\lambda_i^+(k), \mu_j^+(k) \le 1$ by Lemma~\ref{lemma:shrunk-set}, the third inequality uses the Lipschitz property of $F_i$ and $G_j$ in Assumption~\ref{assum:1}, and the last inequality is by induction hypothesis. Similarly, we have
\begin{align}\label{equ:error-gradient-2}
    & \Biggl| \biggl( \sum_{i=1}^{I} \lambda_i^-(k) p_{\mathrm{c},i}^-(k, M)
     - \sum_{j=1}^{J} \mu_j^-(k) p_{\mathrm{s},j}^-(k, M) \biggr) 
     - \biggl( \sum_{i=1}^{I} \lambda_i^-(k) F_i(\lambda_i^-(k))
     - \sum_{j=1}^{J} \mu_j^-(k) G_j(\mu_j^-(k))
    \biggr) \Biggr|\nonumber\\
    \le & 2\epsilon \Biggl[ \sum_{i=1}^{I} L_{F_i}
     \left(1 + L_{F^{-1}_i} \left(p_{\mathrm{c},i,\max} - p_{\mathrm{c},i,\min}\right)\right)
     + \sum_{j=1}^{J} L_{G_j}
    \left( 1 + L_{G^{-1}_j}
    \left(p_{\mathrm{s},j,\max} - p_{\mathrm{s},j,\min}\right)\right)\Biggr].
\end{align}
Hence, by \eqref{equ:norm-g-hat-k}, \eqref{equ:error-gradient-1}, \eqref{equ:error-gradient-2}, and triangle inequality, we have
\begin{align}\label{equ:norm-g-hat-k-bound-temp}
    & \| \hat{\boldsymbol{g}}(k) \|_2 \nonumber\\
    \le &
    \frac{|{\cal E}|}{2\delta} \Biggl| \biggl( \sum_{i=1}^{I} \lambda_i^+(k) F_i(\lambda_i^+(k))
     - \sum_{j=1}^{J} \mu_j^+(k) G_j(\mu_j^+(k)) \biggr)  - \biggl( \sum_{i=1}^{I} \lambda_i^-(k) F_i(\lambda_i^-(k)) - \sum_{j=1}^{J} \mu_j^-(k) G_j(\mu_j^-(k))
    \biggr) \Biggr|\nonumber\\
    & + \frac{2|{\cal E}|\epsilon}{\delta} 
      \Biggl[ \sum_{i=1}^{I} L_{F_i}
     \left(1 + L_{F^{-1}_i} \left(p_{\mathrm{c},i,\max} - p_{\mathrm{c},i,\min}\right)\right)
     + \sum_{j=1}^{J} L_{G_j}
    \left( 1 + L_{G^{-1}_j}
    \left(p_{\mathrm{s},j,\max} - p_{\mathrm{s},j,\min}\right)\right)\Biggr].
\end{align}
Note that by the Lipschitz property of $F_i$ in Assumption~\ref{assum:1}, we have
\begin{align*}
    & |\lambda_i^+(k) F_i(\lambda_i^+(k)) - \lambda_i^-(k) F_i(\lambda_i^-(k))|\nonumber\\
    = & |\lambda_i^+(k) F_i(\lambda_i^+(k)) - \lambda_i^+(k) F_i(\lambda_i^-(k)) + \lambda_i^+(k) F_i(\lambda_i^-(k)) - \lambda_i^-(k) F_i(\lambda_i^-(k))|\nonumber\\
    \le & |\lambda_i^+(k) F_i(\lambda_i^+(k)) - \lambda_i^+(k) F_i(\lambda_i^-(k))| + |\lambda_i^+(k) F_i(\lambda_i^-(k)) - \lambda_i^-(k) F_i(\lambda_i^-(k))|\nonumber\\
    \le & \lambda_i^+(k) L_{F_i} | \lambda_i^+(k) - \lambda_i^-(k) | 
    +  F_i(\lambda_i^-(k)) | \lambda_i^+(k) - \lambda_i^-(k) |\nonumber\\
    = & \left[\lambda_i^+(k) L_{F_i} + F_i(\lambda_i^-(k)) \right]
    | \lambda_i^+(k) - \lambda_i^-(k) |.
\end{align*}
Note that $|\lambda_i^+(k) - \lambda_i^-(k)| = |\sum_{j\in {\cal E}_{\mathrm{c},i}} (x^{+}_{i,j}(k) -  x^{-}_{i,j}(k)) | \le 2 |{\cal E}_{\mathrm{c},i}| \delta $ by definition. Hence, we have
\begin{align}\label{equ:lipschitz-f-1}
    |\lambda_i^+(k) F_i(\lambda_i^+(k)) - \lambda_i^-(k) F_i(\lambda_i^-(k))| \le  2 \delta  |{\cal E}_{\mathrm{c},i}| \left[\lambda_i^+(k) L_{F_i} + F_i(\lambda_i^-(k)) \right] \le 2 \delta  |{\cal E}_{\mathrm{c},i}| ( L_{F_i} + p_{\mathrm{c}, i, \max} ),
\end{align}
where the last inequality holds since $\lambda_i^+(k) \le 1$ and $F_i(\lambda_i^-(k)) \le p_{\mathrm{c}, i, \max}$. Similarly, we have
\begin{align}\label{equ:lipschitz-f-2}
    |\mu_j^+(k) G_j(\mu_j^+(k)) - \mu_j^-(k) G_j(\mu_j^-(k))| \le  2 \delta  |{\cal E}_{\mathrm{s},j}| ( L_{G_j} + p_{\mathrm{s}, j, \max} ).
\end{align}
From \eqref{equ:norm-g-hat-k-bound-temp}, \eqref{equ:lipschitz-f-1}, and \eqref{equ:lipschitz-f-2}, we have
\begin{align*}
    \| \hat{\boldsymbol{g}}(k) \|_2 \le & \frac{|{\cal E}|}{2\delta} \left( \sum_{i=1}^{I} 2 \delta  |{\cal E}_{\mathrm{c},i}| ( L_{F_i} + p_{\mathrm{c}, i, \max} ) + \sum_{j=1}^{J}  2 \delta  |{\cal E}_{\mathrm{s},j}| ( L_{G_j} + p_{\mathrm{s}, j, \max} ) \right)\nonumber\\
    & + \frac{2|{\cal E}|\epsilon}{\delta} 
      \Biggl[ \sum_{i=1}^{I} L_{F_i}
     \left(1 + L_{F^{-1}_i} \left(p_{\mathrm{c},i,\max} - p_{\mathrm{c},i,\min}\right)\right)
     + \sum_{j=1}^{J} L_{G_j}
    \left( 1 + L_{G^{-1}_j}
    \left(p_{\mathrm{s},j,\max} - p_{\mathrm{s},j,\min}\right)\right)\Biggr]\nonumber\\
    = & |{\cal E}| \left( \sum_{i=1}^{I} |{\cal E}_{\mathrm{c},i}| ( L_{F_i} + p_{\mathrm{c}, i, \max} ) + \sum_{j=1}^{J}  |{\cal E}_{\mathrm{s},j}| ( L_{G_j} + p_{\mathrm{s}, j, \max} ) \right)\nonumber\\
    & + \frac{2|{\cal E}|\epsilon}{\delta} 
     \Biggl[ \sum_{i=1}^{I} L_{F_i}
     \left(1 + L_{F^{-1}_i} \left(p_{\mathrm{c},i,\max} - p_{\mathrm{c},i,\min}\right)\right)
     + \sum_{j=1}^{J} L_{G_j}
    \left( 1 + L_{G^{-1}_j}
    \left(p_{\mathrm{s},j,\max} - p_{\mathrm{s},j,\min}\right)\right)\Biggr].
\end{align*}
Substituting the above bound into \eqref{equ:diff-x-bound}, we have
\begin{align}\label{equ:diff-x-bound-2}
    & \| \boldsymbol{x}(k+1) - \boldsymbol{x}(k) \|_2 \nonumber\\
    \le &
    \eta |{\cal E}| \left( \sum_{i=1}^{I} |{\cal E}_{\mathrm{c},i}| ( L_{F_i} + p_{\mathrm{c}, i, \max} ) + \sum_{j=1}^{J}  |{\cal E}_{\mathrm{s},j}| ( L_{G_j} + p_{\mathrm{s}, j, \max} ) \right)\nonumber\\
    & + \frac{2|{\cal E}|\eta \epsilon}{\delta} 
     \Biggl[ \sum_{i=1}^{I} L_{F_i}
     \left(1 + L_{F^{-1}_i} \left(p_{\mathrm{c},i,\max} - p_{\mathrm{c},i,\min}\right)\right)
     + \sum_{j=1}^{J} L_{G_j}
    \left( 1 + L_{G^{-1}_j}
    \left(p_{\mathrm{s},j,\max} - p_{\mathrm{s},j,\min}\right)\right)\Biggr].
\end{align}
Then we have
\begin{align}\label{equ:diff-lambda-bound}
    |\lambda_i^{+}(k+1) - \lambda_i^{+}(k)| = & \biggl| \sum_{j\in {\cal E}_{\mathrm{c},i}} x_{i,j}^{+} (k+1) -  \sum_{j\in {\cal E}_{\mathrm{c},i}} x_{i,j}^{+} (k) \biggr|\nonumber\\
    \le & \sum_{j\in {\cal E}_{\mathrm{c},i}} \left|x_{i,j}^{+} (k+1) - x_{i,j}^{+}(k)\right|  \nonumber\\
    \le &  \| \boldsymbol{x}^{+}(k+1) - \boldsymbol{x}^{+}(k) \|_1 \nonumber\\
    \le &  \sqrt{|{\cal E}|} \| \boldsymbol{x}^{+}(k+1) - \boldsymbol{x}^{+}(k) \|_2\nonumber\\
    = &  \sqrt{|{\cal E}|} \| \boldsymbol{x}(k+1) - \boldsymbol{x}(k) + \delta (\boldsymbol{u}(k+1)-\boldsymbol{u}(k)) \|_2\nonumber\\
    \le & \sqrt{|{\cal E}|} \| \boldsymbol{x}(k+1) - \boldsymbol{x}(k)\|_2 + 2 \delta \sqrt{|{\cal E}|},
\end{align}
where the first inequality is by triangle inequality, the third inequality is by Cauchy-Schwarz inequality, and the last inequality uses triangle inequality and the fact that $\|\boldsymbol{u}(k)\|_2 = 1$ for all $k$. From \eqref{equ:diff-x-bound-2} and \eqref{equ:diff-lambda-bound}, we have
\begin{align*}
    & |\lambda_i^{+}(k+1) - \lambda_i^{+}(k)| \nonumber\\
    \le & \frac{2|{\cal E}|^{3/2}\eta \epsilon}{\delta} 
     \Biggl[ \sum_{i=1}^{I} L_{F_i}
     \left(1 + L_{F^{-1}_i} \left(p_{\mathrm{c},i,\max} - p_{\mathrm{c},i,\min}\right)\right)
     + \sum_{j=1}^{J} L_{G_j}
    \left( 1 + L_{G^{-1}_j}
    \left(p_{\mathrm{s},j,\max} - p_{\mathrm{s},j,\min}\right)\right)\Biggr]\nonumber\\
    & + \eta |{\cal E}|^{3/2} \left( \sum_{i=1}^{I} |{\cal E}_{\mathrm{c},i}| ( L_{F_i} + p_{\mathrm{c}, i, \max} ) + \sum_{j=1}^{J}  |{\cal E}_{\mathrm{s},j}| ( L_{G_j} + p_{\mathrm{s}, j, \max} ) \right) + 2 \delta |{\cal E}|^{1/2},
\end{align*}
Similarly, we can obtain
\begin{align*}
    & |\lambda_i^{-}(k+1) - \lambda_i^{-}(k)| \nonumber\\
    \le & \frac{2|{\cal E}|^{3/2}\eta \epsilon}{\delta} 
     \Biggl[ \sum_{i=1}^{I} L_{F_i}
     \left(1 + L_{F^{-1}_i} \left(p_{\mathrm{c},i,\max} - p_{\mathrm{c},i,\min}\right)\right)
     + \sum_{j=1}^{J} L_{G_j}
    \left( 1 + L_{G^{-1}_j}
    \left(p_{\mathrm{s},j,\max} - p_{\mathrm{s},j,\min}\right)\right)\Biggr]\nonumber\\
    & + \eta |{\cal E}|^{3/2} \left( \sum_{i=1}^{I} |{\cal E}_{\mathrm{c},i}| ( L_{F_i} + p_{\mathrm{c}, i, \max} ) + \sum_{j=1}^{J}  |{\cal E}_{\mathrm{s},j}| ( L_{G_j} + p_{\mathrm{s}, j, \max} ) \right) + 2 \delta |{\cal E}|^{1/2},\nonumber\\
    & |\mu_j^{+}(k+1) - \mu_j^{+}(k)| \nonumber\\
    \le & \frac{2|{\cal E}|^{3/2}\eta \epsilon}{\delta} 
     \Biggl[ \sum_{i=1}^{I} L_{F_i}
     \left(1 + L_{F^{-1}_i} \left(p_{\mathrm{c},i,\max} - p_{\mathrm{c},i,\min}\right)\right)
     + \sum_{j=1}^{J} L_{G_j}
    \left( 1 + L_{G^{-1}_j}
    \left(p_{\mathrm{s},j,\max} - p_{\mathrm{s},j,\min}\right)\right)\Biggr]\nonumber\\
    & + \eta |{\cal E}|^{3/2} \left( \sum_{i=1}^{I} |{\cal E}_{\mathrm{c},i}| ( L_{F_i} + p_{\mathrm{c}, i, \max} ) + \sum_{j=1}^{J}  |{\cal E}_{\mathrm{s},j}| ( L_{G_j} + p_{\mathrm{s}, j, \max} ) \right) + 2 \delta |{\cal E}|^{1/2},\nonumber\\
    & |\mu_j^{-}(k+1) - \mu_j^{-}(k)| \nonumber\\
    \le & \frac{2|{\cal E}|^{3/2}\eta \epsilon}{\delta} 
     \Biggl[ \sum_{i=1}^{I} L_{F_i}
     \left(1 + L_{F^{-1}_i} \left(p_{\mathrm{c},i,\max} - p_{\mathrm{c},i,\min}\right)\right)
     + \sum_{j=1}^{J} L_{G_j}
    \left( 1 + L_{G^{-1}_j}
    \left(p_{\mathrm{s},j,\max} - p_{\mathrm{s},j,\min}\right)\right)\Biggr]\nonumber\\
    & + \eta |{\cal E}|^{3/2} \left( \sum_{i=1}^{I} |{\cal E}_{\mathrm{c},i}| ( L_{F_i} + p_{\mathrm{c}, i, \max} ) + \sum_{j=1}^{J}  |{\cal E}_{\mathrm{s},j}| ( L_{G_j} + p_{\mathrm{s}, j, \max} ) \right) + 2 \delta |{\cal E}|^{1/2}.
\end{align*}
Therefore, we have
\begin{align*}
    & |\lambda_i^{+/-}(k+1) - F_i^{-1} (p_{\mathrm{c},i}^{+/-}(k,M))| \le |\lambda_i^{+/-}(k+1) - \lambda_i^{+/-}(k)| + |\lambda_i^{+/-}(k) -  F_i^{-1} (p_{\mathrm{c},i}^{+/-}(k,M))| \nonumber\\
    \le & \frac{2|{\cal E}|^{3/2}\eta \epsilon}{\delta} 
     \Biggl[ \sum_{i=1}^{I} L_{F_i}
     \left(1 + L_{F^{-1}_i} \left(p_{\mathrm{c},i,\max} - p_{\mathrm{c},i,\min}\right)\right)
     + \sum_{j=1}^{J} L_{G_j}
    \left( 1 + L_{G^{-1}_j}
    \left(p_{\mathrm{s},j,\max} - p_{\mathrm{s},j,\min}\right)\right)\Biggr]\nonumber\\
    & + \eta |{\cal E}|^{3/2} \left( \sum_{i=1}^{I} |{\cal E}_{\mathrm{c},i}| ( L_{F_i} + p_{\mathrm{c}, i, \max} ) + \sum_{j=1}^{J}  |{\cal E}_{\mathrm{s},j}| ( L_{G_j} + p_{\mathrm{s}, j, \max} ) \right) + 2 \delta |{\cal E}|^{1/2}\nonumber\\
    & + |\lambda_i^{+/-}(k) -  F_i^{-1} (p_{\mathrm{c},i}^{+/-}(k,M))|\nonumber\\
    \le & \frac{2|{\cal E}|^{3/2}\eta \epsilon}{\delta} 
     \Biggl[ \sum_{i=1}^{I} L_{F_i}
     \left(1 + L_{F^{-1}_i} \left(p_{\mathrm{c},i,\max} - p_{\mathrm{c},i,\min}\right)\right)
     + \sum_{j=1}^{J} L_{G_j}
    \left( 1 + L_{G^{-1}_j}
    \left(p_{\mathrm{s},j,\max} - p_{\mathrm{s},j,\min}\right)\right)\Biggr]\nonumber\\
    & + \eta |{\cal E}|^{3/2} \left( \sum_{i=1}^{I} |{\cal E}_{\mathrm{c},i}| ( L_{F_i} + p_{\mathrm{c}, i, \max} ) + \sum_{j=1}^{J}  |{\cal E}_{\mathrm{s},j}| ( L_{G_j} + p_{\mathrm{s}, j, \max} ) \right) + 2 \delta |{\cal E}|^{1/2}\nonumber\\
    & + 2\epsilon \left(1 + L_{F^{-1}_i} \left(p_{\mathrm{c},i,\max} - p_{\mathrm{c},i,\min}\right)\right),
\end{align*}
where the last inequality uses the induction hypothesis. By the Lipschitz property of $F_i$ in Assumption~\ref{assum:1}, we have
\begin{align*}
    & |F_i(\lambda_i^{+/-}(k+1)) - p_{\mathrm{c},i}^{+/-}(k,M)| \nonumber\\
    \le & L_{F_i} |\lambda_i^{+/-}(k+1) - F_i^{-1} (p_{\mathrm{c},i}^{+/-}(k,M))|\nonumber\\
    \le & 
    \frac{2 \eta \epsilon |{\cal E}|^{3/2} L_{F_i} }{\delta} 
     \Biggl[ \sum_{i'=1}^{I} L_{F_{i'}}
     \left(1 + L_{F^{-1}_{i'}} \left(p_{\mathrm{c},i',\max} - p_{\mathrm{c},i',\min}\right)\right)
     + \sum_{j'=1}^{J} L_{G_{j'}}
    \left( 1 + L_{G^{-1}_{j'}}
    \left(p_{\mathrm{s},j',\max} - p_{\mathrm{s},j',\min}\right)\right)\Biggr]\nonumber\\
    & + \eta |{\cal E}|^{3/2} L_{F_i} \left( \sum_{i'=1}^{I} |{\cal E}_{\mathrm{c},i'}| ( L_{F_{i'}} + p_{\mathrm{c}, i', \max} ) + \sum_{j'=1}^{J}  |{\cal E}_{\mathrm{s},j'}| ( L_{G_{j'}} + p_{\mathrm{s}, j', \max} ) \right) + 2 \delta |{\cal E}|^{1/2} L_{F_i}\nonumber\\
    & + 2\epsilon L_{F_i} \left(1 + L_{F^{-1}_i} \left(p_{\mathrm{c},i,\max} - p_{\mathrm{c},i,\min}\right)\right).
\end{align*}
Recall the definition of $e_{\mathrm{c},i}$ and $e_{\mathrm{s},j}$ in \eqref{equ:def-e-c} and \eqref{equ:def-e-s}. We have
\begin{align}\label{equ:error-term-i}
    |F_i(\lambda_i^{+/-}(k+1)) - p_{\mathrm{c},i}^{+/-}(k,M)| \le e_{\mathrm{c},i}.
\end{align}
Similarly, we can obtain
\begin{align}\label{equ:error-term-j}
    & |G_j(\mu_j^{+/-}(k+1)) - p_{\mathrm{s},j}^{+/-}(k,M)|\nonumber\\
    \le & \frac{2 \eta \epsilon |{\cal E}|^{3/2} L_{G_j} }{\delta} 
     \Biggl[ \sum_{i'=1}^{I} L_{F_{i'}}
     \left(1 + L_{F^{-1}_{i'}} \left(p_{\mathrm{c},i',\max} - p_{\mathrm{c},i',\min}\right)\right)
     + \sum_{j'=1}^{J} L_{G_{j'}}
    \left( 1 + L_{G^{-1}_{j'}}
    \left(p_{\mathrm{s},j',\max} - p_{\mathrm{s},j',\min}\right)\right)\Biggr]\nonumber\\
    & + \eta |{\cal E}|^{3/2} L_{G_j} \left( \sum_{i'=1}^{I} |{\cal E}_{\mathrm{c},i'}| ( L_{F_{i'}} + p_{\mathrm{c}, i', \max} ) + \sum_{j'=1}^{J}  |{\cal E}_{\mathrm{s},j'}| ( L_{G_{j'}} + p_{\mathrm{s}, j', \max} ) \right) + 2 \delta |{\cal E}|^{1/2} L_{G_j}\nonumber\\
    & + 2\epsilon L_{G_j} \left(1 + L_{G^{-1}_j} \left( p_{\mathrm{s},j,\max} - p_{\mathrm{s},j,\min}\right)\right)\nonumber\\
    = &  e_{\mathrm{s},j}.
\end{align}
From \eqref{equ:error-term-i} and \eqref{equ:error-term-j} and the monotonicity of $F_i$ and $G_j$ in Assumption~\ref{assum:1}, we have
\begin{align*}
    \lambda_i^{+/-}(k+1) \in & [F^{-1}_i(p^{+/-}_{\mathrm{c},i}(k,M) + e_{\mathrm{c},i}) - \epsilon, F^{-1}_i(p^{+/-}_{\mathrm{c},i}(k,M) - e_{\mathrm{c},i}) + \epsilon], \text{for all $i$,}\nonumber\\
    \mu_j^{+/-}(k+1) \in & [G^{-1}_j(p^{+/-}_{\mathrm{s},j}(k,M) - e_{\mathrm{s},j})-\epsilon, G^{-1}_j(p^{+/-}_{\mathrm{s},j}(k,M) + e_{\mathrm{s},j}) + \epsilon], \text{for all $j$}.
\end{align*}
Therefore, we have proved \eqref{equ:next-arr-in-interval-3} and \eqref{equ:next-arr-in-interval-4}, and hence \eqref{equ:next-arr-in-interval-1} and \eqref{equ:next-arr-in-interval-2} hold. Then following the same induction argument as that used to prove \eqref{equ:arr-in-interval-1} and \eqref{equ:arr-in-interval-2}, we can show that for all $m=1,\ldots,M$,
\begin{align}
    \lambda_i^{+/-}(k+1) \in & [F^{-1}_i(\bar{p}^{+/-}_{\mathrm{c},i}(k+1,m))-\epsilon, F^{-1}_i(\underline{p}^{+/-}_{\mathrm{c},i}(k+1,m)) + \epsilon], \text{for all $i$,}\label{equ:arr-in-interval-k-1}\\
    \mu_j^{+/-}(k+1) \in & [G^{-1}_j(\underline{p}^{+/-}_{\mathrm{s},j}(k+1,m))-\epsilon, G^{-1}_j(\bar{p}^{+/-}_{\mathrm{s},j}(k+1,m)) + \epsilon], \text{for all $j$}.\label{equ:arr-in-interval-k-2}
\end{align}
Since in each bisection iteration, the algorithm divides the search interval by 2, we have
\begin{align*}
    \bar{p}^{+/-}_{\mathrm{c},i}(k+1,M) - \underline{p}^{+/-}_{\mathrm{c},i}(k+1,M)
    \le & \frac{\bar{p}^{+/-}_{\mathrm{c},i}(k+1,1) - \underline{p}^{+/-}_{\mathrm{c},i}(k+1,1)}{2^{M-1}} = \frac{2e_{\mathrm{c}, i}}{2^{M-1}}
    \le  \frac{p_{\mathrm{c},i,\max} - p_{\mathrm{c},i,\min}}{2^{M-1}},\nonumber\\
    \bar{p}^{+/-}_{\mathrm{s},j}(k+1,M) - \underline{p}^{+/-}_{\mathrm{s},j}(k+1,M)
    \le & \frac{\bar{p}^{+/-}_{\mathrm{s},j}(k+1,1) - \underline{p}^{+/-}_{\mathrm{s},j}(k+1,1)}{2^{M-1}} = \frac{2e_{\mathrm{s},j}}{2^{M-1}}
    \le  \frac{p_{\mathrm{s},j,\max} - p_{\mathrm{s},j,\min}}{2^{M-1}},
\end{align*}
By the Lipschitz property of $F_i^{-1}$ and $G_j^{-1}$ in Assumption~\ref{assum:1}, we have
\begin{align}
    F^{-1}_i(\underline{p}^{+/-}_{\mathrm{c},i}(k+1,M)) - F^{-1}_i(\bar{p}^{+/-}_{\mathrm{c},i}(k+1,M)) 
    \le & L_{F^{-1}_i} \frac{p_{\mathrm{c},i,\max} - p_{\mathrm{c},i,\min}}{2^{M-1}},\label{equ:diff-arrs-k-1}\\
    G^{-1}_j(\bar{p}^{+/-}_{\mathrm{s},j}(k+1,M)) - G^{-1}_j(\underline{p}^{+/-}_{\mathrm{s},j}(k+1,M)) 
    \le & L_{G^{-1}_j}
    \frac{p_{\mathrm{s},j,\max} - p_{\mathrm{s},j,\min}}{2^{M-1}}.\label{equ:diff-arrs-k-2}
\end{align}
Note that $p^{+/-}_{\mathrm{c},i}(k+1,M) \in [\underline{p}^{+/-}_{\mathrm{c},i}(k+1,M), \bar{p}^{+/-}_{\mathrm{c},i}(k+1,M)]$ and $p^{+/-}_{\mathrm{s},j}(k+1,M) \in [\underline{p}^{+/-}_{\mathrm{s},j}(k+1,M), \bar{p}^{+/-}_{\mathrm{s},j}(k+1,M)]$. Hence, by the monotonicity of $F_i$ and $G_j$ in Assumption~\ref{assum:1}, we have
\begin{align}
    F^{-1}_i(p^{+/-}_{\mathrm{c},i}(k+1,M)) 
    \in &  [F^{-1}_i(\bar{p}^{+/-}_{\mathrm{c},i}(k+1,M))-\epsilon, F^{-1}_i(\underline{p}^{+/-}_{\mathrm{c},i}(k+1,M))+\epsilon],\label{equ:arr-in-interval-k-3}\\
    G^{-1}_j (p^{+/-}_{\mathrm{s},j}(k+1,M)) 
    \in & [G^{-1}_j(\underline{p}^{+/-}_{\mathrm{s},j}(k+1,M))-\epsilon, G^{-1}_j(\bar{p}^{+/-}_{\mathrm{s},j}(k+1,M))+\epsilon].\label{equ:arr-in-interval-k-4}
\end{align}
From \eqref{equ:arr-in-interval-k-1}, \eqref{equ:arr-in-interval-k-2}, \eqref{equ:arr-in-interval-k-3}, and \eqref{equ:arr-in-interval-k-4}, we have
\begin{align*}
    \left| \lambda_i^{+/-}(k+1) - F^{-1}_i(p^{+/-}_{\mathrm{c},i}(k+1,M)) \right| 
    \le & F^{-1}_i(\underline{p}^{+/-}_{\mathrm{c},i}(k+1,M)) - F^{-1}_i(\bar{p}^{+/-}_{\mathrm{c},i}(k+1,M)) + 2\epsilon,\\
    \left| \mu_j^{+/-}(k+1) - G^{-1}_j(p^{+/-}_{\mathrm{s},j}(k+1,M)) \right| 
    \le & G^{-1}_j(\bar{p}^{+/-}_{\mathrm{s},j}(k+1,M)) - G^{-1}_j(\underline{p}^{+/-}_{\mathrm{s},j}(k+1,M)) + 2\epsilon.
\end{align*}
Substituting \eqref{equ:diff-arrs-k-1} and \eqref{equ:diff-arrs-k-2} into the above inequalities, we have
\begin{align*}
    \left| \lambda_i^{+/-}(k+1) - F^{-1}_i(p^{+/-}_{\mathrm{c},i}(k+1,M)) \right| 
    \le & L_{F^{-1}_i} \frac{p_{\mathrm{c},i,\max} - p_{\mathrm{c},i,\min}}{2^{M-1}} + 2\epsilon,\nonumber\\
    \left| \mu_j^{+/-}(k+1) - G^{-1}_j(p^{+/-}_{\mathrm{s},j}(k+1,M)) \right| 
    \le & L_{G^{-1}_j}
    \frac{p_{\mathrm{s},j,\max} - p_{\mathrm{s},j,\min}}{2^{M-1}} + 2\epsilon.
\end{align*}
Recall that $M=\log_2 \frac{1}{\epsilon}$. Then we have
\begin{align}
    \left| \lambda_i^{+/-}(k+1) - F^{-1}_i(p^{+/-}_{\mathrm{c},i}(k+1,M)) \right| 
    \le & 2\epsilon \left[1 + L_{F^{-1}_i} \left(p_{\mathrm{c},i,\max} - p_{\mathrm{c},i,\min}\right)\right],\label{equ:price-error-k-1}\\
    \left| \mu_j^{+/-}(k+1) - G^{-1}_j(p^{+/-}_{\mathrm{s},j}(k+1,M)) \right| 
    \le & 2 \epsilon \left[ 1 + L_{G^{-1}_j}
    \left(p_{\mathrm{s},j,\max} - p_{\mathrm{s},j,\min}\right)\right].\label{equ:price-error-k-2}
\end{align}
We have shown that \eqref{equ:arr-in-interval-k-1}, \eqref{equ:arr-in-interval-k-2}, \eqref{equ:price-error-k-1}, and \eqref{equ:price-error-k-2} hold for $k+1$. Hence, by induction, they hold for all $k$, i.e., we have shown that \eqref{equ:lemma-price-error-1} and \eqref{equ:lemma-price-error-2} hold for all $k$.

Next, we will use \eqref{equ:lemma-price-error-1} and \eqref{equ:lemma-price-error-2} to prove the lemma. We have
\begin{align}\label{equ:diff-lambda-prime}
    \biggl|\lambda'_i(t) - \sum_{j\in {\cal E}_{\mathrm{c},i}} x_{i,j}(k(t))\biggr| 
    =& |\lambda'_i(t) - \lambda^{+/-}_i(k(t))|  + \delta \sum_{j\in {\cal E}_{\mathrm{c},i}}  | u_{i,j}(k(t)) |\nonumber\\
    \le & |\lambda'_i(t) - \lambda^{+/-}_i(k(t))|  + \delta \sqrt{|{\cal E}_{\mathrm{c},i}|} \sqrt{ \sum_{j\in {\cal E}_{\mathrm{c},i}}
    u^2_{i,j}(k(t)) }\nonumber\\
    \le & |\lambda'_i(t) - \lambda^{+/-}_i(k(t))|  + \delta \sqrt{|{\cal E}_{\mathrm{c},i}|} \| \boldsymbol{u}(k(t)) \|_2 \nonumber\\
    = & |\lambda'_i(t) - \lambda^{+/-}_i(k(t))|  + \delta \sqrt{|{\cal E}_{\mathrm{c},i}|},
\end{align}
where the first inequality is by Cauchy-Schwarz inequality. Recall that $\lambda'_i(t) = F^{-1}_i(p^{+/-}_{\mathrm{c},i}(k(t),m(t)))$. Then from \eqref{equ:diff-lambda-prime}, we have
\begin{align}\label{equ:diff-lambda-prime-2}
    & \biggl|\lambda'_i(t) - \sum_{j\in {\cal E}_{\mathrm{c},i}} x_{i,j}(k(t))\biggr|
    \le \big|F^{-1}_i(p^{+/-}_{\mathrm{c},i}(k(t),m(t))) - \lambda^{+/-}_i(k(t))\bigr|  + \delta \sqrt{|{\cal E}_{\mathrm{c},i}|}\nonumber\\
    \le & \bigl|F^{-1}_i(p^{+/-}_{\mathrm{c},i}(k(t),m(t))) - F^{-1}_i(p^{+/-}_{\mathrm{c},i}(k(t),M))\bigr|
    + \bigl|\lambda^{+/-}_i(k(t)) - F^{-1}_i(p^{+/-}_{\mathrm{c},i}(k(t),M))\bigr|  + \delta \sqrt{|{\cal E}_{\mathrm{c},i}|}\nonumber\\
    \le & L_{F^{-1}_i} \bigl| p^{+/-}_{\mathrm{c},i}(k(t),m(t)) - p^{+/-}_{\mathrm{c},i}(k(t),M) \bigr|
    + \bigl|\lambda^{+/-}_i(k(t)) - F^{-1}_i(p^{+/-}_{\mathrm{c},i}(k(t),M))\bigr|  + \delta \sqrt{|{\cal E}_{\mathrm{c},i}|},
\end{align}
where the second inequality uses triangle inequality and the last inequality uses the lipschitz property of $F^{-1}_i$ by Assumption~\ref{assum:1}. For $k(t)\ge 2$, i.e., $t\ge 2MN+1$, the length of search interval for customer queue $i$ in the bisection algorithm~\ref{alg:bisection} is $2e_{\mathrm{c},i}$. Hence, from \eqref{equ:diff-lambda-prime-2}, we have
\begin{align}\label{equ:diff-lambda-prime-3}
    \biggl|\lambda'_i(t) - \sum_{j\in {\cal E}_{\mathrm{c},i}} x_{i,j}(k(t))\biggr|
    \le 2 L_{F^{-1}_i} e_{\mathrm{c},i} + \bigl|\lambda^{+/-}_i(k(t)) - F^{-1}_i(p^{+/-}_{\mathrm{c},i}(k(t),M))\bigr|  + \delta \sqrt{|{\cal E}_{\mathrm{c},i}|},
\end{align}
for $t\ge 2MN+1$. Under the event ${\cal C}$, by \eqref{equ:lemma-price-error-1} and \eqref{equ:diff-lambda-prime-3}, we have
\begin{align}\label{equ:diff-lambda-prime-final}
    \biggl|\lambda'_i(t) - \sum_{j\in {\cal E}_{\mathrm{c},i}} x_{i,j}(k(t))\biggr|
    \le 2 L_{F^{-1}_i} e_{\mathrm{c},i} + 2\epsilon \left[1 + L_{F^{-1}_i} \left( p_{\mathrm{c},i,\max} - p_{\mathrm{c},i,\min}\right)\right]  + \delta \sqrt{|{\cal E}_{\mathrm{c},i}|},
\end{align}
for $t\ge 2MN+1$. Similarly, for $t\ge 2MN+1$, under the event ${\cal C}$, we can show that 
\begin{align}\label{equ:diff-mu-prime-final}
    \biggl|\mu'_j(t) - \sum_{i\in {\cal E}_{\mathrm{s},j}} x_{i,j}(k(t))\biggr|
    \le 2 L_{G^{-1}_j} e_{\mathrm{s},j} + 2\epsilon \left[1 + L_{G^{-1}_j} \left( p_{\mathrm{s},j,\max} - p_{\mathrm{s},j,\min}\right)\right]  + \delta \sqrt{|{\cal E}_{\mathrm{s},j}|}.
\end{align}
Lemma~\ref{lemma:bisection-error} is proved.

\section{Proof of Lemma~\ref{lemma:improved-bisection-error}}
\label{app:proof-lemma-improved-bisection-error}

Compared to the original pricing algorithm, the difference of the \textit{balanced pricing algorithm} is that we use the initial vector $\boldsymbol{x}(1)$ and the initial searching intervals $[\underline{p}^{+/-}_{\mathrm{c}, i}(k,1),\bar{p}^{+/-}_{\mathrm{c},i}(k,1)]$, $[\underline{p}^{+/-}_{\mathrm{s}, j}(k,1), \bar{p}^{+/-}_{\mathrm{s},j}(k,1)]$ in Assumption~\ref{assum:5}. From Assumption~\ref{assum:5}, we have
\begin{align}\label{equ:assum-5-ineq}
    \sum_{j\in {\cal E}_{\mathrm{c}, i}} x_{i,j}(1) \in & \left[F^{-1}_i(\bar{p}^{+/-}_{\mathrm{c},i}(k,1))-\epsilon + \sqrt{|{\cal E}_{\mathrm{c}, i}|} \delta,
    F^{-1}_i(\underline{p}^{+/-}_{\mathrm{c},i}(k,1)) + \epsilon - \sqrt{|{\cal E}_{\mathrm{c}, i}|} \delta\right], \mbox{for all } i\nonumber\\
    \sum_{i\in {\cal E}_{\mathrm{s}, j}} x_{i,j}(1) \in & \left[G^{-1}_j(\underline{p}^{+/-}_{\mathrm{s},j}(k,1))-\epsilon + \sqrt{|{\cal E}_{\mathrm{s}, j}|} \delta, 
    G^{-1}_j(\bar{p}^{+/-}_{\mathrm{s},j}(k,1)) + \epsilon - \sqrt{|{\cal E}_{\mathrm{s}, j}|} \delta \right], \mbox{for all } j,
\end{align}
Note that
\begin{align}\label{equ:7-ineq-detail}
    \lambda_i^+(1) = & \sum_{j\in {\cal E}_{\mathrm{c}, i}} \left( x_{i,j}(1) + \delta u_{i,j}(1) \right)\nonumber\\
    \le & \sum_{j\in {\cal E}_{\mathrm{c}, i}}  x_{i,j}(1) + \delta \sum_{j\in {\cal E}_{\mathrm{c}, i}}  |u_{i,j}(1)|\nonumber\\
    \le & \sum_{j\in {\cal E}_{\mathrm{c}, i}}  x_{i,j}(1) + \delta \sqrt{|{\cal E}_{\mathrm{c}, i}|} \sqrt{\sum_{j\in {\cal E}_{\mathrm{c}, i}}  |u_{i,j}(1)|^2}\nonumber\\
    \le & \sum_{j\in {\cal E}_{\mathrm{c}, i}}  x_{i,j}(1) + \delta \sqrt{|{\cal E}_{\mathrm{c}, i}|} \|\boldsymbol{u}(1)\|_2 = \sum_{j\in {\cal E}_{\mathrm{c}, i}}  x_{i,j}(1) + \delta \sqrt{|{\cal E}_{\mathrm{c}, i}|},
\end{align}
where the third line is by the Cauchy-Schwarz inequality and the last line uses the fact that $\boldsymbol{u}(1)$ is a unit vector. Similarly, we have
\begin{align}\label{equ:7-ineq}
    \lambda_i^+(1) \ge & \sum_{j\in {\cal E}_{\mathrm{c}, i}}  x_{i,j}(1) - \delta \sqrt{|{\cal E}_{\mathrm{c}, i}|}\nonumber\\
    \lambda_i^-(1) \le & \sum_{j\in {\cal E}_{\mathrm{c}, i}}  x_{i,j}(1) + \delta \sqrt{|{\cal E}_{\mathrm{c}, i}|}\nonumber\\
    \lambda_i^-(1) \ge & \sum_{j\in {\cal E}_{\mathrm{c}, i}}  x_{i,j}(1) - \delta \sqrt{|{\cal E}_{\mathrm{c}, i}|}\nonumber\\
    \mu_j^+(1) \le & \sum_{i\in {\cal E}_{\mathrm{s}, j}}  x_{i,j}(1) + \delta \sqrt{|{\cal E}_{\mathrm{s}, j}|}\nonumber\\
    \mu_j^+(1) \ge & \sum_{i\in {\cal E}_{\mathrm{s}, j}}  x_{i,j}(1) - \delta \sqrt{|{\cal E}_{\mathrm{s}, j}|}\nonumber\\
    \mu_j^-(1) \le & \sum_{i\in {\cal E}_{\mathrm{s}, j}}  x_{i,j}(1) + \delta \sqrt{|{\cal E}_{\mathrm{s}, j}|}\nonumber\\
    \mu_j^-(1) \ge & \sum_{i\in {\cal E}_{\mathrm{s}, j}}  x_{i,j}(1) - \delta \sqrt{|{\cal E}_{\mathrm{s}, j}|}.
\end{align}
Therefore, from \eqref{equ:assum-5-ineq}, \eqref{equ:7-ineq-detail}, and \eqref{equ:7-ineq}, we have
\begin{align*}
    \lambda_i^+(1), \lambda_i^-(1) \in \left[F^{-1}_i(\bar{p}^{+/-}_{\mathrm{c},i}(k,1))-\epsilon,
    F^{-1}_i(\underline{p}^{+/-}_{\mathrm{c},i}(k,1)) + \epsilon \right], \mbox{for all } i\nonumber\\
    \mu_j^+(1), \mu_j^-(1) \in \left[G^{-1}_j(\underline{p}^{+/-}_{\mathrm{s},j}(k,1))-\epsilon, 
    G^{-1}_j(\bar{p}^{+/-}_{\mathrm{s},j}(k,1)) + \epsilon  \right], \mbox{for all } j.
\end{align*}
The above results verify the base case of \eqref{equ:arr-in-interval-1} and \eqref{equ:arr-in-interval-2} in the induction argument in the proof of Lemma~\ref{lemma:bisection-error} in Appendix~\ref{app:proof-lemma-bisection-error}. Following the same induction argument as that in the proof of Lemma~\ref{lemma:bisection-error} in Appendix~\ref{app:proof-lemma-bisection-error}, we can obtain
\begin{align}
    \left| \lambda_i^{+/-}(k) - F^{-1}_i(p^{+/-}_{\mathrm{c},i}(k,M)) \right| 
    \le & 2\epsilon \left[1 + L_{F^{-1}_i} \left( p_{\mathrm{c},i,\max} - p_{\mathrm{c},i,\min}\right)\right], \mbox{ for all } i,\label{equ:improved-lemma-price-error-1}\\
    \left| \mu_j^{+/-}(k) - G^{-1}_j(p^{+/-}_{\mathrm{s},j}(k,M)) \right| 
    \le & 2 \epsilon \left[ 1 + L_{G^{-1}_j}
    \left( p_{\mathrm{s},j,\max} - p_{\mathrm{s},j,\min}\right)\right], \mbox{ for all } j. \label{equ:improved-lemma-price-error-2}
\end{align}

Note that \eqref{equ:7-ineq-detail} and \eqref{equ:7-ineq} also hold for iteration $k$. Hence, by triangle inequality, we have
\begin{align}\label{equ:improved-diff-lambda-prime}
    \biggl|\lambda'_i(t) - \sum_{j\in {\cal E}_{\mathrm{c},i}} x_{i,j}(k(t))\biggr| 
    \le  & |\lambda'_i(t) - \lambda^{+/-}_i(k(t))| + |\lambda^{+/-}_i(k(t)) - \sum_{j\in {\cal E}_{\mathrm{c},i}} x_{i,j}(k(t))|\nonumber\\
    \le & |\lambda'_i(t) - \lambda^{+/-}_i(k(t))| + \delta \sqrt{|{\cal E}_{\mathrm{c},i}|}.
\end{align}
Then from \eqref{equ:improved-diff-lambda-prime}, we have
\begin{align}\label{equ:improved-diff-lambda-prime-2}
    & \biggl|\lambda'_i(t) - \sum_{j\in {\cal E}_{\mathrm{c},i}} x_{i,j}(k(t))\biggr|
    \le \big|F^{-1}_i(p^{+/-}_{\mathrm{c},i}(k(t),m(t))) - \lambda^{+/-}_i(k(t))\bigr|  + \delta \sqrt{|{\cal E}_{\mathrm{c},i}|}\nonumber\\
    \le & \bigl|F^{-1}_i(p^{+/-}_{\mathrm{c},i}(k(t),m(t))) - F^{-1}_i(p^{+/-}_{\mathrm{c},i}(k(t),M))\bigr|
    + \bigl|\lambda^{+/-}_i(k(t)) - F^{-1}_i(p^{+/-}_{\mathrm{c},i}(k(t),M))\bigr|  + \delta \sqrt{|{\cal E}_{\mathrm{c},i}|}\nonumber\\
    \le & L_{F^{-1}_i} \bigl| p^{+/-}_{\mathrm{c},i}(k(t),m(t)) - p^{+/-}_{\mathrm{c},i}(k(t),M) \bigr|
    + \bigl|\lambda^{+/-}_i(k(t)) - F^{-1}_i(p^{+/-}_{\mathrm{c},i}(k(t),M))\bigr|  + \delta \sqrt{|{\cal E}_{\mathrm{c},i}|},
\end{align}
where the second inequality uses triangle inequality and the last inequality uses the lipschitz property of $F^{-1}_i$ by Assumption~\ref{assum:1}. For all iterations including the first iteration, the length of bisection search interval for customer queue $i$ under the \textit{balanced pricing algorithm} is $2e_{\mathrm{c},i}$ by Assumption~\ref{assum:5}. Hence, from \eqref{equ:improved-diff-lambda-prime-2}, we have
\begin{align}\label{equ:improved-diff-lambda-prime-3}
    \biggl|\lambda'_i(t) - \sum_{j\in {\cal E}_{\mathrm{c},i}} x_{i,j}(k(t))\biggr|
    \le 2 L_{F^{-1}_i} e_{\mathrm{c},i} + \bigl|\lambda^{+/-}_i(k(t)) - F^{-1}_i(p^{+/-}_{\mathrm{c},i}(k(t),M))\bigr|  + \delta \sqrt{|{\cal E}_{\mathrm{c},i}|},
\end{align}
for all $t\ge 1$. Under the event ${\cal C}$, by \eqref{equ:improved-lemma-price-error-1} and \eqref{equ:improved-diff-lambda-prime-3}, we have
\begin{align*}
    \biggl|\lambda'_i(t) - \sum_{j\in {\cal E}_{\mathrm{c},i}} x_{i,j}(k(t))\biggr|
    \le 2 L_{F^{-1}_i} e_{\mathrm{c},i} + 2\epsilon \left[1 + L_{F^{-1}_i} \left( p_{\mathrm{c},i,\max} - p_{\mathrm{c},i,\min}\right)\right]  + \delta \sqrt{|{\cal E}_{\mathrm{c},i}|},
\end{align*}
for all $t\ge 1$. Similarly, for all $t\ge 1$, under the event ${\cal C}$, we can show that 
\begin{align*}
    \biggl|\mu'_j(t) - \sum_{i\in {\cal E}_{\mathrm{s},j}} x_{i,j}(k(t))\biggr|
    \le 2 L_{G^{-1}_j} e_{\mathrm{s},j} + 2\epsilon \left[1 + L_{G^{-1}_j} \left( p_{\mathrm{s},j,\max} - p_{\mathrm{s},j,\min}\right)\right]  + \delta \sqrt{|{\cal E}_{\mathrm{s},j}|}.
\end{align*}
Lemma~\ref{lemma:improved-bisection-error} is proved.

\vfill

\end{document}